\documentclass[10pt]{amsart}
  \usepackage{amsmath,amssymb,amsthm,graphicx}
  \usepackage{pdfsync}
  \usepackage{color}

\usepackage{enumerate}
\usepackage{verbatim}
\def\wl{\par \vspace{\baselineskip}}
\usepackage[margin=1in]{geometry}
\numberwithin{equation}{section}

\newtheorem{theorem}{Theorem}[section]
\newtheorem{cor}[theorem]{Corollary}
\newtheorem{lemma}[theorem]{Lemma}

\newtheorem{prop}[theorem]{Proposition} \theoremstyle{definition}
\newtheorem{definition}[theorem]{Definition}
 \theoremstyle{remark}
\newtheorem{rem}[theorem]{Remark}

\newtheorem*{theorem*}{\textbf{Theorem}}

\newenvironment{lproof}
	{\begin{proof}}
	{\end{proof}}

\def \tP{ \tilde{P}} \def \tA{ \tilde{A}}
\def \dagS{ S^{\dag}}
\def\beq{\begin{equation} } \def\eeq{\end{equation}}

\def\supp{\hbox{supp}}

\def\tw{\widetilde{w}}
\def\idl{\frac{1}{\delta}} \def\cC{\mathcal C}
\def\eql{equilibrium} \def\qed{$\blacksquare$}

\def\cA{ {\mathcal A} } \def\cB{ {\mathcal B} } \def\cC{ {\mathcal
C} } \def\cD{ {\mathcal D} } \def\cE{ {\mathcal E} } \def\cF{
{\mathcal F} } \def\cG{ {\mathcal G} } \def\cH{ {\mathcal H} }

\def\cJ{ {\mathcal J}} \def\cK{ {\mathcal K} } \def\cL{ {\mathcal
L} } \def\cI{ {\mathcal I}}

\def\cM{{\mathcal M}} \def\cN{ {\mathcal N} } \def\cO{ {\mathcal
O} } \def\cP{{\mathcal P}}

\def\cR{{ \mathcal R }}

\def\cS{{\mathcal S} } \def\cT{{\mathcal T}} \def\cV{{\mathcal V}}
\def\cW{{\mathcal W}}
\def\fX{{ \mathfrak{X}}} \def\fZ{{\mathfrak{Z}}}
\def\fY{{ \mathfrak{Y}}}

\def\tA{\tilde A} \def\tc{\tilde c} \def\tp{\tilde p}

\def\del{\delta}

\def\pt{\partial}

\def\RR{\mathbb R}
\def \CC{\mathbb C}
\def\eps{\varepsilon} \def\Qn{\mathbb R^n_{>0}} \def\Qr{\mathbb
R^r_{>0}}

\def\ben{\begin{enumerate} }

\def\een{\end{enumerate} }

\def \R{ {\mathbb R}}
\def\Rn{ {\mathbb R}^n } \def\Rnp{ {\mathbb R}^n_{>0} } \def\Rnnn{
{\mathbb R}^n_{\geq 0} }

\def\bRplus{ {\mathbb R}^n_{ >0} } \def\bRplusc{ {\mathbb R}^n_{
\geq 0} }
\def \tE{ \tilde{E}}     \def \tu{ \tilde{u}}
\def \tD{ \tilde{ \mathcal{D}}}

\def \Id{ \text{Id} } \def \sgn{ \text{sgn}}
\def \div{ \text{div}} \def \grad{ \text{grad}}
\newcommand{\Cdot}[2]{\dot{\mcC}^{\infty}_c(M)^{\bullet}_{[#1,#2]}}
\newcommand{\Hbar}[2]{ \bar{H}_{s,b}^{-1,1}(M)_{[#1,#2]} }
\newcommand{\Hdot}[4]{ \dot{H}_{s,b}^{#1,#2}(M)_{[#3,#4]} }
\def \mcD{ \mathcal{D} } \def \dmcD{ \dot{\mathcal{D}}}
\def \mcF{\mathcal{F}}
\def \tmcD{ \tilde{ \mathcal{D} } }
\def \mcC{ \mathcal{C}} \def \mcA{ \mathcal{A}}
\def \mcM{ \mathcal{M} } \def \mcH{ \mathcal{H}}
\def \mcV{ \mathcal{V} } \def \mcL{ \mathcal{L}}
\def \mcO{ \mathcal{O} } \def \mcI{ \mathcal{I}}

\def \dD{ \dot{\mcD}}
\def \tphi{\tilde{\phi}}
\def \dmcH{ \dot{ \mathcal{H} } } \def \dmcF{ \dot{ \mathcal{F} } }
\def \eTX{ \left.\right.^eT^*X}
\def\ecTX{ \left.\right.^e\bar{T}^*X}
\def\pr{ \partial_{\xi_r}} \def\pth{ \partial_{\xi_{\theta}}} \def\pz{ \partial_{\xi_z}} \def\ptau{ \partial_{\tau}}
\def\p{ \partial }

\def\hxi{ \hat{\xi} } \def\heta{ \hat{\eta} }
\def\htau{ \hat{\tau} } \def\hzeta{ \hat{\zeta} }

\def\bhxi{ \underline{\hat{\xi}} }
\def\bheta{ \underline{\hat{\eta}} }
\def\bhtau{ \underline{\hat{\tau}} }
\def\bhzeta{ \underline{\hat{\zeta}} }

\def\uxi{ \underline{\xi} }
\def\ueta{ \underline{\eta} }
\def\utau{ \underline{\tau} }
\def\uzeta{ \underline{\zeta} }

\def\ePDO{ \left.\right.^e\Psi}

\def\eS{ \left.\right.^eS } \def\eSX{ \left.\right.^eS^*X }

\def\esig{ \left.\right.^e\sigma }

\def\e{ \left.\right.^e} \def\b{ \left.\right.^b}

\def\tP{ \tilde{P} } \def\u{ \tilde{u}}
\def\P{ \tilde{P}}
\def\d{ \delta } \def \g{ \gamma }

\def \teps{ \tilde{ \epsilon } }
\def \db{ \bullet, - } \def \bd{ -, \bullet }

\def \bt{ \textbf{t}}
\def \l{ \langle} \def \r{ \rangle}
\def \a{ \acute{a}}
\def \v{ \mathbf{v}}
\def \mcB{ \mathcal{B}}
\def \mcE{ \mathcal{E}} \def \hE{ \hat{\mathcal{E}}}
\def \dmcE{ \dot{\mathcal{E}}}
\def \be{ \mathbf{e}}
\def \he{ \hat{\mathbf{e}}}
\def \mcS{ \mathcal{S}}
\def \End{ \textnormal{End}} \def\Diff{ \textnormal{Diff}}
\def\Ell{ \text{Ell}} \def\Ran{ \text{Ran}}
\def\tpsi{ \tilde{\psi}} \def\trho{\tilde{\rho}}
\def\I{ \mathbf{I}}
\def\ord{ \text{ord}}
\def \a{ \acute{a}}
\newcommand{ \IH}[4]{ \mathcal{I}^{#1}H^{#2,#3}_e(#4, \mcM^p(#4))}
\newcommand{\commen}[1]{}
\newcommand{ \lsp}[3]{ #1_{(#2,#3)}}
\newcommand{ \ti}[1]{ \tilde{#1}}
\newcommand{ \todo}[1]{
}
\def \WF{ \textnormal{WF}}

\begin{document}

\title{Diffraction of elastic waves by edges}
\author{Vitaly Katsnelson}

\begin{abstract}
We investigate the diffraction of singularities of solutions
to the linear elastic equation on manifolds with edge singularities. Such manifolds are
modelled on the product of a smooth manifold and a cone over a compact fiber. For
the fundamental solution, the initial pole generates a pressure wave (p-wave), and a
secondary, slower shear wave (s-wave). If the initial pole is appropriately situated near
the edge, we show that when a p-wave strikes the edge, the diffracted p-waves and s-waves
(i.e. loosely speaking, are not limits of p-rays which just miss the edge) are weaker in a Sobolev sense than the incident p-wave. We
also show an analogous result for an s-wave that hits the edge, and provide results for more general situations.
\end{abstract}

\maketitle

\section{Introduction}
The purpose of this paper is to investigate the diffractive behavior of singularities of solutions to the linear elastic equation. In elastic theory, if we consider a bounded isotropic elastic medium with a smooth boundary, then a singular impulse in the interior generates two distinct waves called the pressure wave and the slower, secondary wave or shear wave ($p$-wave and $s$-wave for short). In such a situation, Taylor in \cite{Taylor75}, and Yamamoto in \cite{Yam87} showed that when either of these waves hits the boundary transversely, this interaction will generate at least a $p$ or an $s$ wave moving away from the boundary, with the possibility of both a $p$ and $s$ wave being generated, which often occurs in seismic experiments. To be more precise, this means that if the elastic wave solution $u$ has singularities along the ray path of either wave hitting the boundary, then it will have a singularity of the same Sobolev strength along the ray path (called $p$-ray or $s$-ray) of at least one outgoing, reflected $p$ or $s$ wave. That is, if a solution $u$ fails to be in $H^s$ along an incoming $p$-ray hitting the boundary transversely, then it fails to be in $H^s$ along either the reflected $p$ or $s$ ray. Yamamoto in \cite{Yam87} refined this result considerably that fits with seismic data, by showing that if there is an incoming $p$ wave which hits the boundary at a certain time, then even when there is no incoming $s$ wave at that time, both a reflected $p$ and $s$ wave will be generated moving away from the boundary.
An elementary parametrix construction using ``geometric optics'' solutions easily demonstrates such results as done in \cite{TaylorElasticity}. \\
\indent Of considerable interest is what happens when the medium has edges or corners beyond the codimension one case. In particular, one is interested in what happens when a singular impulse generates $p$ and $s$ waves that approach such an edge (indeed, it is known by the results in \cite{Hor} and Dencker \cite{Dencker} that along such waves in the interior of the medium, the solution must have the same Sobolev strength singularity). For a simpler, scalar wave equation Vasy in \cite{VasyCorners} showed that if a solution to the wave equation is singular along an incoming ray (these turn out to be geodesics in the manifold) approaching an edge of codimension $\geq 2$, then it will generally produce singularities along a whole cone-generating family of outgoing rays (i.e. a cone of geodesics moving away from the edge). This result did not distinguish as seen in experimental physics between the stronger ``geometric'' waves versus the weaker, ``diffracted'' waves (when speaking of weaker and stronger waves, we mean in the Sobolev sense where we measure which Sobolev space solutions lie in along certain bicharacteristics corresponding to the wave operator). \\
\indent Over time, many results were obtained describing propagation of singularities on singular manifolds, but they also did not show whether diffracted rays were weaker than the incoming ray as seen in experimental physics. In a remarkable breakthrough, Melrose, Vasy, and Wunsch in \cite{MWCones,MVWEdges,MVWCorners} showed how to distinguish between weaker ``diffracted'' waves and the other waves to show that under a certain ``nonfocusing'' assumption (which in a model case of manifolds with a warped product metric, would mean that in cylindrical coordinates, one is able to smooth out the solution even a little bit beyond its overall regularity by merely smoothing out its angular coordinates), the solution is smoother along the outgoing diffracted front by an amount related to the codimension of the edge being hit. They even confirmed the intuition that the diffracted waves are precisely those that together with the incoming wave, cannot be approximated in the Sobolev sense by waves which just miss the edge by an infinitesimal amount. As a start, our goal is to obtain such results for the linear elastic equation, which is a nonscalar setting. Unfortunately, even though we have some conjectures on what propagation of singularities looks like in this setting, the nonscalar nature of the problem amplifies the complexity by a considerable degree and we do not have a useful result in this direction yet. Nevertheless, distinguishing between regular and ``diffracted'' $p$-waves and $s$-waves is considerably easier, and it is precisely this direction we pursue in this paper. Indeed, under certain semi-global hypotheses, we will show what happens on the diffracted front of a $p$ and $s$ wave hitting an edge transversely.

\subsection{Basic Setup}
The setting will be a $n$-manifold $X$ with boundary equipped with an \emph{edge metric}, which is called an {\it edge manifold}. The way to visualize this is by taking a manifold with corners and then introducing cylindrical coordinates near an {\it edge}, by blowing up the edge and introducing coordinates on this blow up.\footnote{The edge manifolds we work with here are not exactly this type of blowup since the fiber at the boundary of such a blow up would have corners, but our edge manifolds have boundaryless fibers. Nevertheless, when one stays away from the corners of the fibers on the blown-up manifold, it represents a good visualization of the manifolds in the setting of this paper. See section \ref{sec:edge manifolds} for a more precise description.} Precisely, the boundary of $X$ has a fibration
 $$ Z \to \p X \xrightarrow{\pi_0} Y, \text{ with compact fiber }Z, \quad \text{dim}(Z)= f.$$
 Also, $X$ has a boundary defining function $x$, and near $\p X$ the metric is of the form
 $$ g = dx^2 + \ti{\pi}_0^*h + x^2k$$
 with $h \in C^{\infty}([0,\epsilon) \times Y; \text{Sym}^2T^*([0,\epsilon) \times Y))$ and $k \in C^{\infty}(U; \text{Sym}^2T^*M)$; we further assume that $h|_{x=0}$ is a nondegenerate metric on $Y$ and $k|_{x=0}$ is a nondegenerate fiber metric. Here we extended the fibration $\pi_0$ to a fibration $\ti{\pi}_0:U \to [0,\epsilon) \times Y$ on a neighborhood $U$ of $\p X$.
 \\

{\it (A motivational non-example)} Before proceeding further, we want to motivate how edge manifolds will arise naturally in practice by considering a non-example of a manifold with corners. Suppose that near an edge of some manifold with corners $X$, we have the coordinates $x_1, \dots, x_{f+1}, y_1, \dots, y_{n-f-1}$ and the edge is given by the vanishing of $x_1, \dots, x_{f+1}$. Since we are interested in understanding what happens when a wave interacts with the edge, as well as getting extra information along the diffracted waves, we introduce `cylindrical coordinates' as
$$ x = \sqrt{x_1^2 + \cdots x_{k+1}^2}, \ z_j = x_j/x, \ y_i.$$
This will transform the standard Riemannian metric into an edge metric. If one does a real blow-up of this edge, then the above coordinates act as local projective coordinates on the blown-up manifold. The fibers of the blowup have corners given by the vanishing of some of the $z_j$. Away, from such corners in the fibers, this blow up is exactly the setting of edge manifolds we consider in this paper, where the fibers do not have corners.

\indent Since we work with the linear elastic equation, set $M= \R_t \times X$, which is an $n+1$ dimensional edge manifold representing the space-time setting. The boundary of $M$ still has a fibration with compact fiber $Z$ and base $Y \times \RR_t$. Local coordinates on $M$ will be denoted
$$(t,x,y,z) = (t,x,y_1,\dots,y_{n-f-1},z_1, \dots,z_f).$$

\indent We consider distributional solutions $u \in \mcD'(M;TX)$ to the elastic equation
\begin{equation}\label{eq:intro elastic equation}
 Pu = (D_t^2 - L)u = (D_t^2 - \nabla^* \mu \nabla - \div^* (\lambda + \mu) \div+ R_0)u = 0 \text{ on }M
\end{equation}
where $\nabla$ is the Levi-Civita connection on $X$ pulled back to $M$ via the projection $p:M \to X$, $\div$ is the divergence operator on sections of $TX$ pulled back to the manifold $M$ via $p$, $R_0 \in \mcC^{\infty}(M;\text{End}(TX))$, and $\mu, \lambda \in \mcC^{\infty}(X)$ are the Lam\'{e} parameters.

    We shall consider below only solutions of (\ref{eq:intro elastic equation}) lying in some `finite energy space', which plays an analogous role as setting boundary conditions. Thus, let us denote $\mcD_{\alpha}$ as the domain of $L^{\alpha/2}$, where $L$ is the Friedrichs extension of the operator above, also labeled $L$, on the space $\dot{\mcC}^{\infty}(X;TX)$, of smooth functions vanishing to infinite order at the boundary. We require that a solution be {\it admissable} in the sense that it lies in $\mcC(\R; \mcD_{\alpha})$ for
    {\it some} $\alpha \in \R$.

    As described in \cite{MVWEdges}, in terms of adapted coordinates $t,x,y,z$ near a boundary point of $M$, an element of $\mathcal{V}_e(M)$ is locally an arbitrary smooth combination of the basis vector fields
    \begin{equation}
   x\p_t, \ x\partial_x,\ x\partial_{y_j}, \ \partial_{z_k}
    \end{equation}
and so $\mathcal{V}_e(M)$ is equal to the space of all sections of a vector bundle, which is called the {\it edge tangent bundle} and denoted $\e TM.$ This bundle is canonically isomorphic to the usual tangent bundle over the interior (and non-canonically isomorphic to it globally) with a well-defined bundle map
$ \e TM \to TM$ which has rank $f$ over the boundary. As we will justify shortly, we should think of the fiber coordinate, $z_j$, as angular coordinates, with dual coordinates $\zeta_j$ being the angular momentum. The dual bundle is the \emph{edge cotangent bundle}
 $$ \e T^*M;$$
 it is spanned by $\frac{dt}{x}, \frac{dx}{x},\frac{dy_j}{x}, dz_j$, with corresponding dual coordinates
 $$\tau, \xi, \eta_j, \zeta_j.$$

 Such bundles and vector fields show up naturally when studying the wave operator or the elastic operator since in cylindrical coordinates, these operators are shown to be products of vector fields in $x^{-1}\mcV_e(M)$. Nevertheless, Vasy already showed in \cite{VasyCorners} that singularities of solutions to the wave equation should be described by a different bundle called the \emph{$b$-cotangent bundle}, denoted $\b T^*M$ (which is the dual to the $b$-tangent bundle, denoted $\b TM$, whose basis elements are locally described by $x\p_x, \ \p_{y_j}, \ \p_{z_k}$; see Section \ref{sec: edge calculus} for complete definitions). Thus, if we want to describe the propagation of singularities for the elastic equation, this would be the most natural bundle to use as well. However, since $P$ is not a $b$-operator, trying to obtain such results would be very complicated for two reasons: first, the interaction between edge operators and b-operators requires a significant effort to describe. Secondly, $P$ no longer has a scalar principal symbol, so trying to find clever $b$-operators that are positive along the Hamilton flow associated to $P$ have so far been too challenging to pursue. We expect that in the $b$-setting, a $p$-wave hitting $\partial M$, would give rise to a whole cone of singularities as in the scalar wave equation, but should also give rise to $s$-waves as well. A more manageable task that we pursue here is to at least describe the diffractive behavior of an incoming $p$ and $s$ wave. The edge setting is precisely adapted for this purpose.

As commonly known, the characteristic set of the elastic operator $P$, denoted $\Sigma \subset \e T^*M$, can be decomposed into two mutually disjoint sets corresponding to two waves with different wave speeds, called the pressure wave and the shear wave. Indeed, if we denote $\sigma(P)$ as the principle symbol of the elastic operator, then $\text{det}(\sigma(P))$ is the product of principal symbols of two \emph{scalar} wave equations with different sound speeds. These are the $p$ and $s$ waves, and it's precisely the characteristic sets of these two scalar waves which determine the characteristic set of $P$. We will use the notation
$$ \Sigma := \text{det}(\sigma(P))^{-1}(0)= \Sigma_p \cup \Sigma_s$$
to describe $\Sigma$ as the union of characteristic sets for the $p$ and $s$ waves (see Section \ref{sec: principal symbol and Hamilton vector field} and \ref{sec: p/s bicharactersitics} for a precise description of this and the definitions that follow). Hence, the notions of elliptic, glancing, and hyperbolic sets make sense for each of these scalar waves, so we can refer to the elliptic/hyperbolic/glancing set of $P$ in terms of the $p$ and $s$ waves, but we have to be sure to specify which of the two elliptic,hyperbolic, or glancing sets we are referring to. We will use superscripts and subscripts $`p'$ and $`s'$ in the notation for various sets and functions to denote which wave we are referring to; if we do not want to specify, then we will just write $`p/s'$ for such superscripts and subscripts. In order to fix things, lets assume in this introduction that we are working inside $\Sigma_p$, i.e. we are going to work with the bicharacteristic flow of the pressure wave, which means we are inside the elliptic region of the s-waves. Local coordinates on $M$ with their respective dual coordinates provide local coordinates for $\e T^*M$ denoted
$$ (t, x,y,z, \tau, \xi, \eta, \zeta) \in \e T^*M.$$

With the notation of \cite{MVWEdges}, for each normalized point
\begin{equation}\label{eq:intro normalized p point}
 \alpha = (\bar{t},\bar{y}, \bar{z},\bar{\tau} = \pm 1, \bar{\eta}) \in \mcH^p \Leftrightarrow |\bar{\eta}|<c_p^{-1},
 \end{equation}
where $0< c_p \in \mcC^{\infty}(X)$ denotes the speed of a $p$-wave,
it was shown that there are two line segments of `normal' null bicharacteristics in $\Sigma_p$, each ending at one of the two points above $\alpha$ inside $\e T^*_{\p M}M$ given by solutions $\bar{\xi}$ of $\bar{\xi}^2 + \bar{\eta}^2=c_p^{-2}$. These will be denoted
$$ \mcF^p_{\bullet,\alpha},$$
where $\bullet$ is permitted to be $I$ or $O$, for `incoming' or `outgoing', as $\text{sgn}\bar{\xi} = \pm \text{sgn}\bar{\tau}$ (+ for $I$ and $-$ for $O$). Thus, one should view a $p$-bicharacteristic $\mcF^p_{I,\alpha}$ hitting the boundary at a point above $\alpha$ and then immediately exiting the boundary along another $p$-bicharacteristic $\mcF^p_{O,\alpha'}$ where $\alpha'$ lies in the same fiber $Z_{(\bar{t},\bar{y})}$ as $\alpha$, i.e. they only differ by their $z$ coordinate. The exact relation between $\alpha$ and $\alpha'$, and the relation between the Sobolev regularity of $u$ along $\mcF^p_{I,\alpha}$ versus its Sobolev regularity along $\mcF^p_{O,\alpha'}$ is the main interest of this paper in order to describe the diffraction of waves.
The sets $\mcF^p_{\bullet,\alpha}$ are quite explicit when the fibration and metric are of true product form
$$ dx^2 +h(y,dy) + x^2k(z,dz).$$
Then the principal symbol corresponding to the $p$-wave is simply
$$ q_{p} = \frac{ \tau^2 - c^2_{p}(\xi^2 + |\eta|^2_h + |\zeta|^2_k)}{x^2}.$$
When $c_p$ is constant, the bicharacteristics (i.e. the flow curves of the Hamilton vector field $H_{q_p}$ inside $\e T^*M$) hitting the boundary are simply
$$ \mcF^p_{I,\alpha} = \{ t \leq \bar{t}, x = c^2_p(\bar{t} -t)|\bar{\xi}|, y = y(t), z = \bar{z}, \tau = \bar{\tau}, \xi = \bar{\xi}, \eta = \eta(t), \zeta = 0 \}$$
and
$$\mcF^p_{O,\alpha} = \{ t \geq \bar{t}, x = c^2_p(t-\bar{t})|\bar{\xi}|, y = y(t), z = \bar{z}, \tau = \bar{\tau}, \xi = \bar{\xi}, \eta = \eta(t), \zeta = 0 \};$$
where $(y(t), \eta(t))$ evolves along a geodesic in $Y$ with speed $c_p$ which passes through $(\bar{y},\bar{\eta})$ at time $t = \bar{t}$, and where $\bar{\tau}^2 = c_p^2(\bar{\xi}^2 + |\eta|^2) = 1,$ and we have chosen the sign of $\bar{\xi}$ to agree/disagree with the sign of $\bar{\tau}$ in the incoming/outgoing cases. The case of $\mcF^s_{\bullet,\alpha}$ is almost the same, except one has a different wave speed denoted $c_s < c_p$. This is exactly analogous to the example given in the introduction of \cite{MVWEdges}.

\subsection{Past results on the wave equation} We summarize some results taken from \cite[Section 1]{MVWEdges}. When one considers the standard wave operator $\Box = D_t^2 - \Delta_g$, then one has the same definitions and notation above except with $c_p = 1$. In the model case, similar to above, for each normalized point
$$\alpha = (\bar{t},\bar{y}, \bar{z},\bar{\tau} = \pm 1, \bar{\eta}) \in \mcH, \quad |\bar{\eta}|<1,$$
the bicharacteristics are
$$ \mcF_{I,\alpha} = \{ t \leq \bar{t}, x = (\bar{t} -t)|\bar{\xi}|, y = y(t), z = \bar{z}, \tau = \bar{\tau}, \xi = \bar{\xi}, \eta = \eta(t), \zeta = 0 \}$$
and
$$\mcF_{O,\alpha} = \{ t \geq \bar{t}, x = (t-\bar{t})|\bar{\xi}|, y = y(t), z = \bar{z}, \tau = \bar{\tau}, \xi = \bar{\xi}, \eta = \eta(t), \zeta = 0 \}.$$
As it is $Z$-invariant over the boundary, we may write $\mcH$ as the pull-back to $\p M$ via $\pi_0$ of a corresponding set $\dot{\mcH}$. One may therefore consider all the bicharacteristics meeting the boundary in a single fiber, with the same `slow variables' $(t,y)$ and set
$$ \dot{\mcF}_{\bullet,q} = \bigcup_{p \in \pi_0^{-1}(q)} \mcF_{\bullet,p}, \ q\in \dot{\mcH}.$$
These are pencils of bicharacteristics touching the boundary at a given location in the `slow' spacetime variables $(t,y)$, with given momenta in those variables; the union over all $(t,y)$ of such families form smooth coisotropic (involutive) manifolds in the cotangent bundles near the boundary. Then Melrose, Vasy, Wunsch have already shown Snell's Law for the wave equation, stating that tangential momentum and energy is preserved when a wave interacts with an edge, in the form of the following theorem:

\begin{theorem}(\cite[Theorem 1.1]{MVWEdges} For an admissable solution, $u$, to the wave equation and any $q \in \dot{\mcH}$ \footnote{See Section \cite[Section 1]{HVRadial} for the definitions of $\WF^k$},
$$ \dot{\mcF}^o_{I,q} \cap \textnormal{WF}^k(u) = \emptyset \Rightarrow \dot{\mcF}^o_{O,q} \cap \textnormal{WF}^k(u) = \emptyset.$$
\end{theorem}
(The fact that $\eta(q)$ is the same for the incoming and outgoing rays is the preservation of tangential momentum even though $\p \dot{\mcF}_{I,q}$ and $\p \dot{\mcF}_{O,q}$ \footnote{These refer to the endpoints at the boundary of the respective families of bicharacteristics; see Section \ref{sec: p/s bicharactersitics} for a precise description} are different. The fact that the rays stay in the characteristic set shows the energy preservation.) As mentioned already, this is the type of theorem that has remained elusive for the elastic equation since its proof for the wave equation relies heavily on the fact that $\Box$ is an operator with a scalar principal symbol.

The analysis in \cite{MVWEdges} then distinguishes between the ``diffracted'' waves and the ``geometric'' waves. Indeed, let $o \in X$ be near the boundary, and
$$ u_o(t) = \frac{\sin t\sqrt{\Delta}}{\sqrt{\Delta}}\d_o$$
be the fundamental solution. Then $u_0(t) \in H^s_{loc}(M)$, and one has the following theorem

\begin{cor}(\cite[Corollary 1.4]{MVWEdges}) For all $o \in X^o$ let $\mcL_o$ denote the flowout of $SN^*(\{o\})$ along bicharacteristics lying over $X^o$. If $o$ is sufficiently close to $\p X$, then for short time, the fundamental solution $u_0$ is a Lagrangian distribution along $\mcL_0$ lying in $H^{-n/2+1-0}$ together with a diffracted wave, singular only at $\mcF_O$, that lies in $H^{-n/2+1+f/2-0}$, away from its intersection with $\mcL_o$.
\end{cor}

\begin{figure}[ht]
	\centering
	\includegraphics{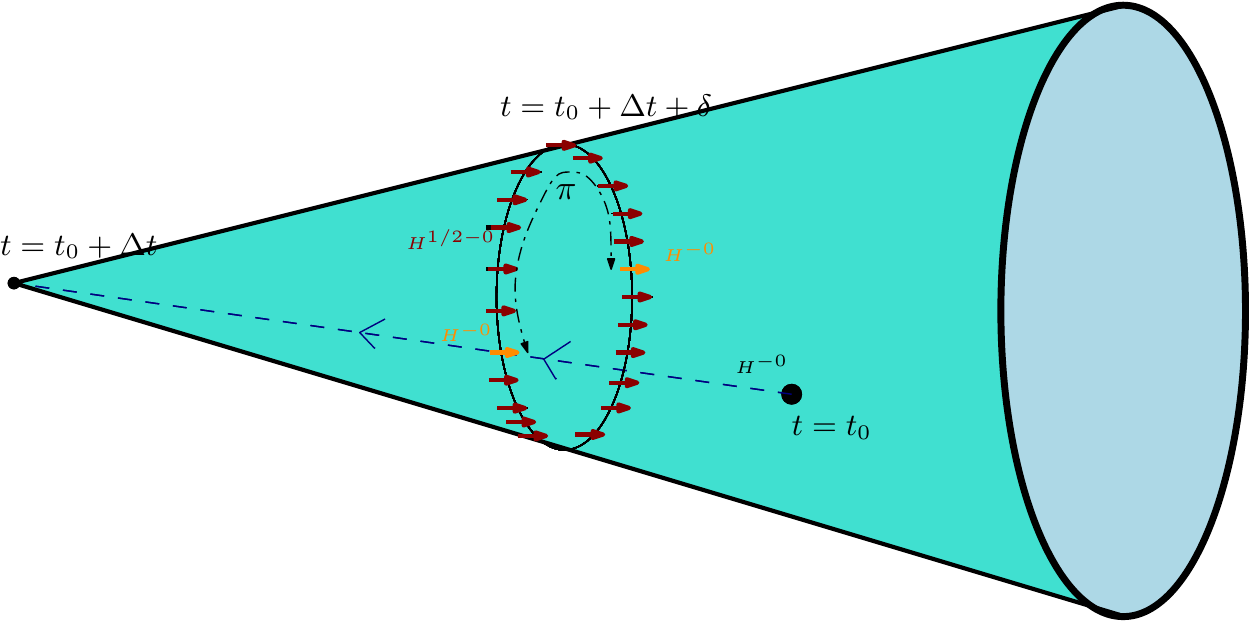}
	\caption{For the scalar wave equation (on a simple cone), an incident wave (navy blue) of Sobolev order $H^{-0}$ starts at time $t=t_0$ and reaches the cone tip at $t=t_0+\Delta t$. An infinite circular family of waves is diffracted, the dark red ones being weaker of Sobolev order $H^{1/2-0}$. However, two rays (orange) are geometrically related to the incident ray having its Sobolev strength as well. They are an angular distance $\pi$ away from the incident ray. }
\end{figure}

\subsection{Sketch of the results}
We will obtain an analogous result for the elastic equation. However, the situation becomes more interesting because there are two waves to consider. Indeed, the unique feature is that
$$ \dot{\mcH}^p \cap \dot{\mcH}^s \neq \emptyset.$$
For example, for the $\alpha$ introduced in (\ref{eq:intro normalized p point}), one also has
$$ |\bar{\eta}| < c_p^{-1} < c_s^{-1}$$
so points that lie above $\alpha$ are solutions $\bar{\xi}$ of
$ \bar{\xi}^2 + \bar{\eta}^2 = c_p^{-2}$ or of $\bar{\xi}^2 + \bar{\eta}^2 = c_s^{-2}$. Thus, for a solution to the homogeneous elastic equation $u$, a singularity of $u$ may enter the boundary along a particular ray in $\dot{\mcF}^p_{I,\pi_0(\alpha)}$ and then exit the boundary along not just rays in $\dot{\mcF}^p_{O,\pi_0(\alpha)}$, but along rays in $\dot{\mcF}^s_{O,\pi_0(\alpha)}$ as well. With $\alpha'$ in the same $Z$ fiber as $\alpha$, a ray $\mcF^{p/s}_{O,\alpha'}$ is a \emph{geometric} $p$-bicharacteristic if this ray, together with $\mcF^p_{I,\alpha}$ is locally a limit of $p$-bicharacteristics lying in $T^*M^o$ that just miss the edge. Otherwise, it is \emph{diffractive}. In the case of the scalar wave equation, the fundamental solution is less singular in a Sobolev sense on the ``diffractive'' bicharacteristics, but has the same Sobolev strength singularity as the incident wave along the geometric bicharacteristics. For the elastic equation however, the incoming $p$-ray $\mcF^p_{I,\alpha}$ together with the outgoing $s$-ray $\mcF^s_{O,\alpha'}$, could never be a limit of $p$-bicharactersics and so $\mcF^s_{O,\alpha'}$ must be diffractive, in which we would expect an improvement in the Sobolev order of $u$ along such an $s$-ray.

 If $(t_0,o) \in M^o$ then the solution to
 $$Pu = \d_{t_0,{o}},$$
vanishing for $t < t_0$, is called the \emph{forward fundamental solution}, where $\d_{t_0,o}$ denotes the delta distribution. Thus, one consequence of our main theorem is

\begin{figure}[ht]
	\centering
	\includegraphics{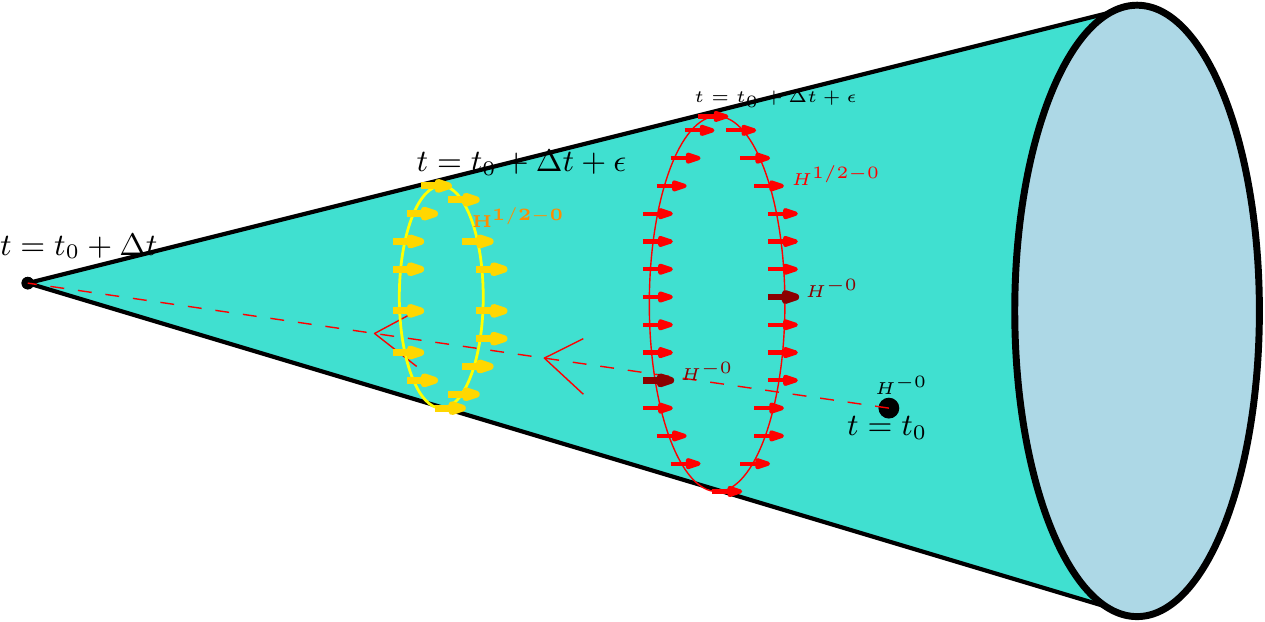}
	\caption{For the elastic wave equation, an incident $p$-wave (red) say, of Sobolev order $H^{-0}$ starts at time $t=t_0$ and reaches the cone tip at $t=t_0+\Delta t$. An infinite circular family of $p$-waves (red) and $s$-waves (yellow) is diffracted, the red ones being weaker of Sobolev order $H^{1/2-0}$. However, two rays (dark red) are geometrically related to the incident ray having its Sobolev strength as well. The $s$-waves are not geometrically related to the incident $p$-wave and are weaker, being of Sobolev strength $H^{1/2-0}$.}
\end{figure}

\begin{theorem}\label{thm:Diffraction}
For all $(t_0,o) \in M^o$ let $\mathcal{L}^p_{t_0,o}$ resp. $\mathcal{L}^s_{t_0,o}$, denote the flowout of $SN^*(\{o\})$ along $p$-bicharacteristics, resp. $s$-bicharacteristics, lying over $M^o$, which lie over $o$ at $t=t_0$. If $o$ is sufficiently close to $\partial X$, then for a short time beyond when the first $p$-wave emanating from $o$ at time $t_0$ hits the edge, the forward fundamental solution $u = u_{t_0,o}$ is a Lagrangian distribution along $\mathcal{L}_{t_0,o} := \mathcal{L}^p_{t_0,o} \cup
\mathcal{L}^s_{t_0,o}$ lying in $H^s$ for all $s < -n/2 +1$ together with diffracted waves, singular only at $\mcF_O^p \cup \mcF_O^s$, that lie in $H^r$ for all $r < -n/2 +1 +f/2$, away from its intersection with $\mathcal{L}_{t_0,o}$.
More precisely, if we consider the first incoming $p$-wave transverse to the boundary, i.e. $u \in H^{-n/2+1-0}$ along $\mcF^p_{I,\alpha}$ but $\textnormal{WF}(u) \cap \mcF^p_{I,\alpha} \neq \emptyset$, then each of the outgoing diffracted $p$ and $s$ waves are weaker in the sense that
$u \in H^r$ along the diffracted $p$ and $s$-bicharacteristics generated by $\mcF^p_{I,\alpha}$ for all $r < -n/2 + 1 + f/2$. Similarly, if we consider the first incoming $s$-wave transverse to the boundary, i.e. $u \in H^s$ along $\mcF^s_{I,\alpha}$ then each of the outgoing diffracted $p$ and $s$ waves are weaker in the sense that
$u \in H^r$ along the diffracted $p$ and $s$-bicharacteristics generated by $\mcF^s_{I,\alpha}$ for all $r < -n/2 + 1 + f/2.$
\end{theorem}

\subsection{Plan for the proof}
To prove the theorem, we will adopt an approach similar to the one presented in \cite[Section 1.3]{MVWEdges}. First, we would like to prove an analog of H\"{o}rmander propagation theorem (see \cite[Section 1.2]{HVRadial} ) for the bicharacteristic flow of $H_{q_{p/s}}$ inside $\e T^*M$ for rays that approach the boundary, which we labeled $\mcF^{p/s}$ earlier. However, such an approach runs into the obstruction presented by manifolds with \emph{radial points}, which occur at the boundary and at which the Hamilton flow vanishes. Hence, such points become saddle equilibria for the Hamilton flow, and $\mcF^{p/s}_{I/O}$ form part of the stable ($I$) or unstable ($O$) manifolds of such equilibria that are transversal to the boundary $x=0$. The other stable/unstable manifolds of the these \emph{critical manifolds} of equilibria are contained in the boundary $x=0$. Thus, the type of propagation of singularities result we want is that a singularity enters the boundary $x=0$ along (say) the stable manifold of one of these critical manifolds, propagating through the critical manifold and out through its unstable manifold; propagating across the boundary to the stable manifold of the other critical manifold; and then through it and back out of the boundary along the corresponding unstable manifold. The key is that the propagation across the boundary leaves the variables $(t,y)$ unaffected, and so this process will show us which bicharacteristics inside $\mcF^{p/s}_O$ will be the geometric continuations of an incident bicharacteristic (say) $\mcF^p_{I,\alpha}$.

The problem is that propagation into or out of a radial point is subject to a threshold amount of regularity that one may propagate. In particular, propagation results into and out of the boundary along bicharacteristics in the edge cotangent bundle up to a given Sobolev order are restricted by the largest power of $x$ by which $u$ is divisible, relative to the corresponding scale of edge Sobolev spaces. Thus, propagating (say) $H^s$ regularity along $\mcF^p_{I}$ into a radial point, will only lead to $H^{s'}$ regularity along certain rays in $\mcF^{p}_O$, but $s'$ may be much smaller than $s$ to provide any useful information.

Hence, we must initially settle for less information. It will turn out that $\dot{\mcF}^p_{I}$ is a \emph{coisotropic} submanifold of the cotangent bundle and as such `coisotropic regularity' with respect to it may be defined in terms of iterated regularity under the application of pseudodifferential operators with symbols vanishing along $\dot{\mcF}^p_{I,\bullet}$. Thus, we begin by showing that coisotropic regularity in this sense, of any order, propagates through the boundary, with a fixed loss of derivatives. Then under the assumption that a solution $u$ lies locally in time in some fixed energy space for the operator $L$, we prove a propagation of semi-global regularity to show that it must lie in such an energy space for all times. By interpolation of such a result combined with the coisotropic regularity propagation, it will follow that coisotropic regularity propagates into, along, and out of the boundary with epsilon derivative loss.

 As we are in a non-scalar setting, we cannot directly adopt the commutator techniques used to prove such results, since such techniques heavily rely on the scalar wave equation setting. Nevertheless if (say) we are trying to propagate along $\dot{\mcF}^p_I$, we may project an elastic wave solution $u$ to the s-wave eigenspace of $p$; the projected distribution will satisfy an elliptic type of equation which will provide simpler elliptic estimates for this piece. This is because even though $\mcH^p \cap \mcH^s$ is not always empty, in the bundle $\e T^*M$, $\Sigma_p$ and $\Sigma_s$ are disjoint, which means that $\sigma_e(P)$ will always have an elliptic eigenvalue (i.e. $q_p \neq 0$ or $q_s \neq 0$ at each point in $\e T^*M$ not in the $0$-section). For the piece $u$ projected to the $p$-wave eigenspace, we'll actually be able to adapt a commutator proof as in the scalar wave equation. Hence, in Section \ref{sec:partial elliptic regularity}, we'll prove partial elliptic estimates for the $`s'$ part of $u$. Sections \ref{sec: edge propagation} and \ref{sec: cosio regularity} will be used to prove a propagation result for the $`p'$-part of $u$. We then combine both of these two results to yield the full propagation of coisotropic regularity of $u$ under certain semi-global hypotheses (see Corollary \ref{Cor:PropOfCoisotropy}). Afterwards, we dualize the argument to obtain the propagation of coinvolutivity (this is analogous to the ``nonfocusing'' condition introduced in \cite{MVWEdges}). In the final section, combining the propagation of coisotropic regularity with the dual notion of coinvolutivity, we will be able to interpolate to prove the main theorem.

\section{Edge and \MakeLowercase{B}-calculus}\label{sec: edge calculus}
The edge calculus of pseudodifferential operators was introduced by Mazzeo in \cite{MazzEdge} and a full summary of their wavefront set and composition properties was given by Wunsch and Melrose in \cite[Section 5]{MWCones}. We will use the exact notation appearing in \cite[Section 3]{MVWEdges} for the calculus so we will avoid repeating it here. The notation uses $\Psi^{*,*}_e(M)$ to denote the bifiltered algebra of pseudodifferential edge operators on $\mcC^{-\infty}(M)$.

As in \cite[Section 3]{MVWCorners}, we also fix a non-degenerate $b$-density $\nu$ on $M$, i.e. $\nu$ is of the form $x^{-1}\nu_0$, $\nu_0$ a non-degenerate $\mcC^{\infty}$ density on $M$, which is a nowhere-vanishing section of the density bundle $\Omega M := |\bigwedge^{n}|(M)$. The density gives an inner product on $\dot{\mcC}^{\infty}(M)$. When below we refer to adjoints, we mean this relative to $\nu$, but the statements listed below not only do not depend on $\nu$ of the stated form, but would even hold for any non-degenerate density $x^{-l}\nu_0$, $\nu_0$ as above, $l$ arbitrary, as the statements listed below imply conjugation by $x^l$ preserves the calculi.

An important feature is that we have principal symbol maps $\sigma_{e,m}$
$$ \sigma_{e,m}: \Psi_e^{m,l}(M) \to x^l \Big[ S^m_{\text{phg}}(\e T^*M) / S_{\text{phg}}^{m-1}(\e T^*M) \Big ];$$
the range space for $\sigma_e$ can be conveniently identified with $\mcC^{\infty}(\e S^*M)$

If $A \in \Psi_e^{m,l}(M)$ and $B \in \Psi_e^{m',l'}(M)$ then
     \begin{equation*}\sigma_{e,m+m'-1,l+l'}(i[A,B]) = \{ \sigma_{e,m,l}(A),  \sigma_{e,m,l}(B) \} := H_{e,\sigma_{e,m,l}(A)}(\sigma_{e,m,l}(B)),\end{equation*}
    where the Poisson bracket is computed with respect to the singular symplectic structure on $\e T^*M$ described above, and $H_{e,\sigma_{e,m,l}(A)}$ is the edge-Hamilton vector field.
If $\mathbf{A} = \{A_{\d}\}_{\{\d \in [0,1]\}}$ is a uniformly bounded family in $\Psi_e^{m,l}(M)$ (sometimes written $A_{\d} \in L^{\infty}([0,1]_{\d},\Psi_e^{m,l})$) then
    $$ q \notin \WF'_e(\mathbf{A}) (\text{ sometimes written } \WF'_{e,L^{\infty}}(A_{\d}))$$
    if there exists a $B \in \Psi_e^{0,0}(M)$ such that $BA_{\d}$ is uniformly bounded in $\Psi_e^{-\infty,l}(M)$.

There is a continuous quantization map (by no means unique)
$$ \text{Op}_e: x^l S^m_{\text{phg}}(\e T^*M) \to \Psi_e^{m,l}(M)$$
which satisfies
\begin{align*}
\sigma_{e,m,l}(\text{Op}_e(a)) = [a] \in
x^{l}&S^m_{\text{phg}}(\e T^*M) / S_{\text{phg}}^{m-1}(\e T^*M)
\ \forall a \in x^lS^m_{\text{phg}}(\e T^*M) \text{ and }\\
& \WF'_e (\text{Op}_e(a)) \subset \text{ess supp}(a).
\end{align*}
\indent Associated with the edge calculus there is a scale of Sobolev spaces. For integral order these may be defined directly. Thus for $k \in \mathbb{N}$ and any $s \in \mathbb{R}$ we set
\begin{align}\label{eq: edge Sobolev space}
&H^{k,s}_e(\mathbb{R} \times X) = \{ u \in x^sL^2_{b,loc}(\mathbb{R} \times X);\\
& Pu \in x^s L^2_{b,loc}( \mathbb{R} \times X) \ \forall P \in \Diff_e^{k}(X) \}, \ k \in \mathbb{N}.
\nonumber
\end{align}
where
$$ L^2_b(M) = \{f: \int |f|^2 \nu < \infty \}.$$
For general orders, the edge Sobolev spaces can be defined using the calculus.
\begin{definition}
$u \in H^{m,l}_e(M) \Leftrightarrow \Psi_e^{m,-l}(M) \cdot u \subset L^2_b(M)$.
\end{definition}
The usual properties for Sobolev spaces and wavefront sets in the standard PsiDO setting carry over to these spaces and a summary may be found in \cite[Section 3]{MWCones}.

The passage of the above calculus to vector bundles is only notational with all the essential properties preserved. For any vector bundle $E$ over a manifold $M$, we denote
$ \Psi_e^{m,l}(M;E)$ as the bi-filtered $\star$-algebra with all the properties described above, except we in addition use trivializations of $E$ to construct the operators locally. Elements of this algebra are now maps
$$ \Psi_e^{m,l}(M;E) \ni A : \dot{\mcC}^{\infty}(M;E) \to \dot{\mcC}^{\infty}(M;E),$$
$$ \Psi_e^{m,l}(M;E) \ni A : \mcC^{-\infty}(M;E) \to \mcC^{-\infty}(M;E),$$
and
$$ \Psi_e^{m,l}(M;E) \ni A : H_e^{m',l'}(M;E) \to H_e^{m'-m,l'+l}(M;E),$$
with $H_e^{m',l'}(M;E)$ defined analogously to the scalar case.
The principal symbol maps are the same, except locally inside a trivialization, $A$ is a matrix of edge operators and $\sigma_e(A)$ is a matrix of symbols. Precisely, we have
$$ \sigma_{e,m,l}: \Psi^{m,l}_e(M;E) \to
x^lS^m_{hom}(\e T^*M \setminus o; \pi^*\text{Hom}(E,E)),$$
where $\pi : \e T^*M \to M$ is the bundle projection, and $S^m_{hom}$ denotes homogeneous degree $m$, $\mcC^{\infty}$ functions on $\e T^*M \setminus o$, while
$$ \sigma_{e,m,l}: \Psi^{m,l}_{e\infty}(M;E) \to
x^l\frac{S^m(\e T^*M \setminus o; \pi^*\text{Hom}(E,E))}
{S^{m-1}(\e T^*M \setminus o; \pi^*\text{Hom}(E,E))}.
$$
are equivalence classes of symbols.

As explained in \cite[Section 3]{VasyForms}, the only addition caveat is that for $B_j \in \Psi^{m_j,l_j}_e(M;E)$, it is not necessarily true that $[B_1,B_2]$ becomes lower order, i.e.\ it does not necessarily lie in the space $\Psi_e^{m_1+m_2-1,l_1+l_2}(M;E)$ since the principal symbols of $B_1$ and $B_2$ may fail to commute. However, suppose $B_1, B_2$ are principally \emph{scalar}, i.e. a multiple of the identity homomorphism:
$$ \sigma_{e,m_j,l_j}(B_j) = x^{l_j} b_j \Id, \ b_j \in S^{m_j}_{hom}(\e T^*M\setminus o),$$
then the principal symbols do commute and their commutator is
$$ [B_1,B_2] \in \Psi^{m_1+m_2-1,l_1 + l_2}_e(M;E)$$
with
$$ \sigma_{e,m_1+m_2-1,l_1+l_2}([B_1,B_2]) = x^{l_1+l_2}i(H_{e,b_1}b_2)\Id.$$

On the other hand, suppose now that only $B_1$ has a scalar principal symbol of above. Then $\sigma_{e,m_1,l_1}(B_1)$ and $\sigma_{e,m_2,l_2}(B_2)$ commute, hence
$$ \sigma_{e,m_1+m_2,l_1+l_2}([B_1,B_2]) = 0$$
so
$$ [B_1,B_2] \in \Psi^{m_1+m_2-1,l_1 + l_2}_e(M;E).$$

The $b$-calculus is now exactly analogous, and a good exposition may be found in \cite[Section 2 and 3]{VasyCorners}. In the next section we will describe the relevant manifolds and bundles where we do our microlocal analysis.

\section{Edge Manifolds and Bundles}\label{sec:edge manifolds}
In this section, we will give a concrete description of edge manifolds and edge metrics, and then give several examples. We will then describe the Hamilton vector fields associated with the elastic operator. This exposition is taken almost verbatim from \cite{MVWEdges}.

\subsection{Edge Manifolds and Edge Metrics}
Let $X$ be an $n$-dimensional manifold with boundary, where the boundary, $\p X$ is the total space of a fibration
$$ Z \rightarrow \p X \xrightarrow{\pi_0} Y,$$
where $Y,Z$ are without boundary. Let $b$ and $f$ respectively denote the dimensions of $Y$ and $Z$ (the `base' and the `fiber').
  As in \cite[Proposition 2.1]{MVWEdges}, we can choose change coordinates to get a convenient form of the edge metric \footnote{We never actually need this simplified form and all arguments go through without it, except it makes calculations simpler in several places.}:
\begin{equation}\label{prop:edgeMan nice edge metric}
g = dx^2+h(x,y,dy)+xh'(x,y,z,dy)+x^2k(x,y,z,dy,dz).
\end{equation}
The essential properties of edge manifolds and metrics are already described in \cite{MVWEdges} so we simply refer the reader there for the basic definitions.

\subsection{Principal symbols and Hamilton vector fields}\label{sec: principal symbol and Hamilton vector field}
In this part, we will use the edge bundles just described to give a nice description of the operator $P$, its principal symbol, and its Hamilton flow. Recall that $g$ denotes the edge metric on $X$, and $\tau, \xi, \eta, \zeta$ are the fiber coordinates on the bundle $\e T^*M$.
From now on, we'll denote the canonical coordinates on $\e T^*M$ as $(t,x,y,z,\tau,\xi,\eta,\zeta) \equiv (w,\tau,\vartheta)$. As a coordinate free description, the elastic operator is given by $P= D_t^2 - L$ where
\begin{align*}
 L  &= \nabla^*\mu\nabla  + \div^*(\lambda + \mu)\div + R_0
\end{align*}
with all operators as described in the introduction.
The upshot of using the edge cotangent bundle is that we now naturally have
$P \in x^{-2}\text{Diff}^2_e(M;TX)$, and
$\sigma_e(P) \in x^{-2}C^{\infty}(\e T^*M \setminus o; \pi^*\text{End}(TX))$, denoting the principal symbol of $P$. In a local coordinate chart, where $TX$ is trivilialized using the coordinate trivialization, we have
\begin{equation}
\sigma_e(P)(w,\tau,\vartheta)= \left(\frac{\tau^2}{x^2} - \mu |\vartheta|^2_g \right) \otimes \Id - (\lambda + \mu) \frac{ \vartheta \otimes \vartheta}{x^2} \in x^{-2}S^2_{hom}(\e T^*M ; \text{End}(TX)),
\end{equation}
keeping in mind that we view $g^{-1}$ as a metric on the fibers of $\e T^*X$.
It will be convenient to denote $p = \sigma_e(P)$ and $\tilde{g} = x^{-2}g$. Then we can easily compute
$$\text{det}(p) = q_p (q_s)^{n-1}$$
where
$$q_p = \frac{\tau^2 - c^2_p |\vartheta|^2_{\tilde{g}}}{x^2}, \ \ q_s = \frac{\tau^2 - c^2_s |\vartheta|^2_{\tilde{g}}}{x^2},$$
 with $c_p = \sqrt{\lambda + 2\mu}$, and $c_s = \sqrt{\mu}$, where $\mu, \lambda + \mu$ are assumed to be strictly positive. These correspond to the principal symbols for the $p$-wave and $s$-wave respectively.

In order to connect with the notation used in the introduction, the characteristic set of $p$, i.e. $\Sigma = \Sigma_{det(p)} = det(x^2p)^{-1}(0)$, can
then be decomposed into two disjoint sets
$$ \Sigma = \Sigma_{q_p} \cup \Sigma_{q_s} := \Sigma_p \cup \Sigma_s,$$
with $\Sigma_{p/s}$ given by the vanishing of $x^2q_{p/s}$.

In order to get a propagation result, we must look at the Hamilton vector field of $q$ (with $q$ being either $q_p$ or $q_s$) as a section of the tangent space of the edge cotangent bundle, i.e. $H_q \in C^{\infty}(\e T^*X; T( \e T^*X))$. By considering two new edge metrics
 \begin{equation}
 g_{p/s} := c_{p/s}^{-2}g,
 \end{equation}
 then $q_{p/s}$ are the principal symbols of the wave operators obtained from these metrics.

 With the notation of the edge metric in Proposition \ref{prop:edgeMan nice edge metric}, let $(H_{ij})$ and $(K_{ij})$ (which are nondegenerate) be defined respectively as the $dy \otimes dy$ and $dz \otimes dz$ parts of $h$ and $k$ at $x=0$. Let $(H^{ij})$ and $(K^{ij})$ denote the inverses, and an $O(x^k)$ term denotes $x^k$ times a function in $\mcC^{\infty}(X)$.
 Hence, we may copy down for later use the computation done in \cite[Equation (2.4)]{MVWEdges} adapted to the edge metrics $g$:
\begin{align}\label{eq:Hamilton vector field computation}
-\frac{1}{2}x^2&H_{q_{p/s}} = -\tau x \p_t + \tau \xi c^2_{p/s} \p_{\tau}
+ c^2_{p/s}\xi x\p_x + (c^2_{p/s}\xi^2 + c^2_{p/s}|\zeta|^2_{\bar{K}})\p_{\xi}
\\
&+
(c^2_{p/s}\zeta_i \bar{K}_{p/s}^{ij} + O(x)) \p_{z_j}
+(-\frac{1}{2}c^2_{p/s}\zeta_i \zeta_j \frac{ \p \bar{K}^{ij}}{\p z_k}
- c_{p/s}|\zeta|^2_{\bar{K}} \frac{\p c_{p/s}}{\p z_k}
+ O(x)) \p_{\zeta_k}
 \nonumber \\
&+ (x c^2_{p/s}\eta_jH^{ij} + O(x^2))\p_{y_i}
+ (c^2_{p/s}\xi \eta_i + O(x))\p_{\eta_i};
\nonumber
\end{align}
where as in \cite[Equation (2.4)]{MVWEdges}, $\bar{K}^{ij}$ denotes a term of the form $K^{-1} + O(x)$.

As usual, it is convenient to work with the cosphere bundle, $\e S^*M$, viewed as the boundary `at infinity' of the radial compactification of $\e T^*M$. Introducing the new variable
$$ \sigma = |\tau|^{-1},$$
we have a lemma taken from \cite[Section 2]{MVWEdges} whose proof is almost verbatim in our setting
\begin{lemma}(\cite[Lemma 2.3]{MVWEdges})
$$ \text{Inside $\e S^*M \cap  \Sigma$, $-\frac{1}{2}x^2\sigma H_{q_{p/s}}$ vanishes exactly at $x=0, \ \hzeta =0$.}$$
\end{lemma}

Let the linearization of $-(1/2)x^2\sigma H_{q_{p/s}}$ at $q \in \Sigma \cap \e S^*M$ (where $x=0, \hzeta =0$) be $A_q$.
We then have the following taken directly from \cite[Lemma 2.3]{MVWEdges} and its proof, but rewritten to include the weights $c_{p/s}$:

\begin{lemma}
For $q \in \dot{\mcH}^{p/s}$, i.e. such that $\hxi(q) \neq 0$, the eigenvalues of $A_q$ are $-c_{p/s}^2\hxi$, $0$, and $c^2_{p/s}\hxi$, with $dx$ being an eigenvector of eigenvalue $c^2_{p/s}\hxi$. Moreover, modulo the span of $dx$, the $-c^2_{p/s}\hxi$-eigenspace is spanned by $d\sigma$ and the $d\hzeta_j$.
\end{lemma}

\begin{rem}(\cite[Remark 2.4]{MVWEdges})\label{rem: eigenvalues for stable/unstable}
This shows in particular that the space of the $d\hzeta_j$ (plus a suitable multiple of $dx$) is invariantly given as the stable/unstable eigenspace of $A_q$ inside $T^*_q\e S^*M$ according to $\hxi>0$ or $\hxi < 0$. We denote this subspace of $T^*_q \e S^*M$ by $T^{*,-}_q( \e S^*M)$.
\end{rem}

 Our main focus will be to understand those bicharacteristics associated to $q_p$ and $q_s$ which approach the boundary $\partial M$ transversely. Even though in our case we only care about the bicharacteristic flow in $\e T^*M$, general broken bicharacteristics are usually defined in $\b T^*M$ so we will adopt the notation in \cite[Section 7]{MVWEdges}, and proceed to write down the relevant concepts adapted to our setting.

 Let $\pi$ denote the bundle map $\e T^*M \to \b T^*M$ given in canonical coordinates by
 $$ \pi(t,x,y,z,\tau,\xi,\eta,\zeta) = (t,x,y,z,\tau, x\xi, \eta,x\zeta).$$ The compressed cotangent bundle is defined by setting
 $$ \b \dot{T}^*M = \pi( \e T^*M)/Z,$$
 $$ \dot{\pi}: \e T^*M \to \b \dot{T}^*M$$
 the projection, where, here and henceforth, the quotient by $Z$ acts only over the boundary, and the topology is given by the quotient topology. The cosphere bundles
 $$\b S^*M, \ \b \dot{S}^*M, \ \e S^*M$$
  are naturally defined in an analogous manner as done in Section \ref{sec: edge calculus}.

  Next, it will be convenient to denote
  $$\ti{q}_{p/s}= x^2 q_{p/s}$$
  so that $\ti{q}_{p/s}$ is smooth up to the boundary.
  Observe that
  $$ \ti{q}_{p/s}|_{x=0} = \tau^2 - c_{p/s}^2(\xi^2 + |\eta|_h^2 + |\zeta|_k^2).$$
  Hence, we have that on $\Sigma$ (which is away from the $0$-section of $\e T^*M$), $\tau \neq 0$ so that restricted to $\Sigma$, non-zero covectors are mapped to non-zero covectors by $\pi$ and $\dot{\pi}$ (that is, $\partial M$ is \emph{non-characteristic}). Thus, $\pi,\dot{\pi}$ define maps, denoted with the same letter:
 $$ \pi: \Sigma \to \b S^*M, \ \dot{\pi}: \Sigma \to \b \dot{S}^*M, \ \Sigma = (det(x^2p))^{-1}(\{0\})/\mathbb{R}^+.$$
 (Note that here we denote $\Sigma,\Sigma_p,\Sigma_s$ as before except now as a subset of $\e S^*M$ when we quotient out by the $\mathbb{R}^+$ action on the fibers.)
 We also set
 $$ \dot{\Sigma} = \dot{\pi}(\Sigma) = \dot{\pi}(\Sigma_p) \cup
 \dot{\pi}(\Sigma_s)= \dot{\Sigma}_p \cup \dot{\Sigma}_s;$$
these are called the \emph{compressed characteristic sets}. As mentioned in the introduction, we can now define the `elliptic', `glancing', and `hyperbolic' sets separately for the $p$ and the $s$ waves:
\begin{align*}
&\mathcal{E}^p = \pi(\e S^*_{\p M} M) \setminus \pi(\Sigma_p) \\
&\mathcal{G}^p = \{ q \in \pi(\e S^*_{\p M}M):
\text{Card}(\pi^{-1}(q) \cap \Sigma_p ) = 1\} \\
& \mathcal{H}^p = \{ q \in \pi(\e S^*_{\p M}M):
\text{Card}(\pi^{-1}(q) \cap \Sigma_p ) \geq 2\}.
\end{align*}
In coordinates, we have
$$ \pi(0,t,y,x,\xi,\hat{\tau},\hat{\eta},\hat{\zeta}) =
(0,t,y,x,0,\hat{\tau},\hat{\eta},0)$$ hence the three sets are defined by $\{ \bhtau^2 < c_p^2 |\bheta|^2\}, \ \{ \bhtau^2 = c_p^2 |\bheta|^2\}, \ \{ \bhtau^2 > c_p^2 |\bheta|^2\}$ respectively inside $\pi(\e S^*_{\p M}M)$, which is given by $x = 0, \ \bhzeta = 0, \ \bhxi = 0$.

Continuing the notation in \cite{MVWEdges} we may define the corresponding set in $\b \dot{S}^*_{\p M}M$( hence quotiented by $Z$ and denoted with a dot):
\begin{align*}
&\dot{\mathcal{E}}^p = \mathcal{E}^p / Z,\\
&\dot{\mathcal{G}}^p = \mathcal{G}^p / Z, \\
&\dot{\mathcal{H}}^p = \mathcal{H}^p / Z.
\end{align*}
Naturally, we have the analogous sets for the $s$-waves:
$$\mathcal{E}^s,\ \mathcal{G}^s,\ \mathcal{H}^s, \ \dot{\mathcal{E}}^s,
\ \dot{\mathcal{G}}^s,\ \dot{\mathcal{H}}^s.$$
Notice now that $\dot{\Sigma}_p \cap \dot{\Sigma}_s \neq \emptyset$. This is precisely why an incoming $p$-wave hitting the boundary may generate both $p$ and $s$ waves traveling away from the boundary. We will now make such notions very precise by discussing bicharacteristics.

\subsection{$p/s$-Bicharacteristics}\label{sec: p/s bicharactersitics}
In order to better understand what we mean by $p$-waves and $s$-waves, let us define the notion of bicharacteristics as done in \cite{MVWEdges}.

\begin{definition}
Let the flow of $H_{q_p}$ inside $\e T^*M^o$ be called a $p$-bicharacteristic, and the flow of $H_{q_s}$ be called an $s$- bicharacteristic.
\end{definition}

 We now explain these notions of incoming/outgoing and make the connection with the notation $\mcF, \dot{\mcF}$ presented in the introduction. Given $\alpha \in \mcH^p$, then $\S$2 in \cite{MVWEdges} shows that there exist unique maximally extended incoming/outgoing $p$-bicharacteristics $\gamma_{I/O}$ (\emph{incoming} just means $\gamma_I$ approaches the boundary as $t$ increases, while \emph{outgoing} means $\gamma_O$ moves away from the boundary as $t$ increases), where $\sgn \ \xi(\alpha) = \pm \sgn \ \tau(\alpha),$ such that $\alpha = \p (\bar{\gamma_{\bullet}});$ we denote these curves
$$ \mcF^p_{\bullet,\alpha} \subset \e T^*M.$$
Likewise, for $\beta \in \dot{\mcH}^p$ we let
$$ \dmcF^p_{\bullet,\beta} = \bigcup_{\alpha \in \pi_0^{-1}\beta}
\mcF^p_{\bullet,\alpha}.$$
As in \cite{MVWEdges} we abuse notation slightly to write
$$ \dmcH^p_{I/O} = \p \dmcF^p_{I/O}$$
for the endpoints of incoming/outgoing hyperbolic $p$-bicharacteristics at the boundary. We also define all these sets for $p$ replaced by $s$ for the $s$-bicharacteristics. We end with a crucial remark to explain our notation throughout the paper.

\begin{rem}
Even though all the sets just defined are natural subsets of $\b T^*M$, nevertheless we will abuse notation and view $\mcH^p_{I/O}, \mcF^p_{I,O},\mcH^s_{I/O}, \mcF^s_{I,O}$ as sitting inside $\e T^*M$ instead. This is because most of our analysis is done on the edge cotangent bundle rather than the b-bundle. Concretely, for $\alpha \in \mcH^p$
one has
$$ \pi^{-1}(\alpha )= (t,x=0,y,z,\tau, \xi = \pm \sqrt{ c_p^{-2}\tau^2 - |\eta|^2}, \eta, \zeta = 0 ) \in \mcH^p_{I/O}.$$
The same goes for the $s$-version of these sets. We do this in order to stay consistent with the notation used in \cite{MVWEdges} and to avoid introducing new notation which offers very little distinction with the notation already introduced.
\end{rem}

We now introduce the analog to forward and backward geodesic flow, which was Definition 7.12 in \cite{MVWEdges}.
\begin{definition}\label{def:geometricFlow}
Let $q,q' \in \mcH^p$, with $\pi_0(q) = \pi_0(q')$. We say that
$$ \mcF^p_{I,q}, \mcF^p_{O,q'} $$
are related under the forward geometric flow (and vice-versa under the backward flow) if there exists a $p$-bicharacteristic in $\e T^*_{\p M} M$ whose limit points are $q, q' \in \e T^*_{\p M}M$ with the identification introduced in the previous remark. In such a case, we sometimes write
$$ q \sim_G q'$$
to signify that they are ``geometrically'' related.
If $a \in \mcF^p_{I,q}$ with $q \in \mcH$, we let the forward flowout of $a$ to be the union of the forward $p$-bicharacteristic segment through $a$ and all the $\mcF^p_{O,q'}$ that are related to $\mcF^p_{I,q}$, under the forward geometric flow (and vice-versa for backward flow). We make the analogous definitions for the $s$-geometric flow.

Also, even though bicharacteristics might sometimes be infinitely extended, the sets $\dot{\mcF}^{p/s}_{I/O}$ are smooth manifolds only for short times near $t(q)$. Thus, when we refer to such sets in our proofs, we are assuming some underlying time interval near $t(q)$ where they are well-defined as manifolds.
\end{definition}

\begin{rem}
By definition $\mcF^p_{I,q}$ and $\mcF^s_{O,q'}$ can never be geometrically related by the definition given. Hence, for any $q \in \mcH^p_I$, then $\mcF^s_{O,q'}$ will always be a nongeometric (i.e. diffracted) ray generated by $\mcF^p_{I,q}$ whenever $\pi_0(q) = \pi_0(q')$.
\end{rem}

\section{Domains}\label{sec: domains}
 In this section, we will describe the Friedrichs form domain for the elastic operator, and then use that to identify the Dirichlet form domain for $P$. This will help us identify some basic properties for solution to the elastic equation to be used later when we prove regularity results.
For edge propagation, it will be essential to identify the domains of powers of
\begin{align*}
L &=  2\nabla_s^* \mu \nabla_s  +  \div^*\lambda\div
\\
&= \nabla^* \mu \nabla  +  \div^*(\lambda + \mu)\div + R_0
\end{align*}
introduced in the introduction, where $R_0 \in x^{-2}\mcC^{\infty}(M, \text{End}(TX))$ is a bundle endomorphism, and $\nabla_s$ is the algebraic symmetrization of $\nabla$ making $\nabla_su$ a symmetric $(1,1)$ tensor for $u \in \mcC^{\infty}$. The equality above relating $\nabla_s$ and $\nabla$ just follows from well known Weitzenb\"{o}ck identities. Precisely, one has
\begin{equation}
\nabla_s^*\nabla_s = \frac{1}{2}\nabla^*\nabla  +  \frac{1}{2}\div^*\div + R_0.
\end{equation}

Now, the metric $g$ on $X$ allows us to define the Hilbert space $L^2_g(X;TX)$, and so we start with the following definition.
\begin{definition}
For $u,v \in \dot{\mcC}^{\infty}(X;TX)$, we define the sesqilinear form associated to $L$
$$ B(u,v):=\int_X \lambda(\text{div}(u))( \text{div}(\bar{v})) dg
+2\int_X \mu(\nabla_s u , \nabla_s \bar{v}) dg,$$
where $(\cdot,\cdot)$ denotes the metric inner product on $TX \otimes T^*X$.
This allows us to define the quadratic form domain
$$ \mathcal{D} = \text{cl}\{ \dot{\mcC}^{\infty}(X;TX) \ w.r.t. \ B(u,u) + ||u||^2_{L^2_g(X;TX)} \}
$$
The Friedrich's form domain is then just
$$\text{Dom}(L_{Fr}) = \{ u \in \mathcal{D}: Lu \in L^2_g(X;TX) \}. $$
We also let $\mcD_s$ denote the corresponding domain of $L^{s/2}$.
\end{definition}

  Notice that we automatically have $H_e^{1,1-(f+1)/2} \subset \mcD$ along with the inequality $$ ||v||_{\mcD} \lesssim ||v||_{H_e^{1,1-(f+1)/2}}$$
 for $v \in \dot{\mcC}^{\infty}(X, TX)$. This is because in a local coordinate chart where $TX$ is trivialized, terms such as $||\text{div}(v)||$ and $||\nabla_s v||$ may be estimated by $||x^{-1}v||, ||D_x v||, ||D_y v||$, and $||x^{-1}D_zv||$. The key point is that the reverse inequality is true as well when we have $f>1$. Indeed, we have
$$||x^{-1}D_z v||^2 + ||D_xv||^2 + ||D_yv||^2 \lesssim ||v||^2+||\nabla v||^2
\lesssim ||v||^2+B(v,v) = ||v||_{\mcD}^2.$$
However, we also have by Hardy's inequality that for $f>1$
$$||x^{-1}v||^2 \lesssim ||D_xv||^2 + ||v||^2 \lesssim ||\nabla v||^2 +||v||^2 \lesssim ||v||^2_{\mcD}.$$
Hence, just as in \cite{MVWEdges} we have
\begin{lemma}
If $f>1$, then $\mcD = H_e^{1,1-(f+1)/2}(X;TX)$.
\end{lemma}

As in \cite[Section 5]{MVWEdges} we remark that multiplication of $\mcC^{\infty}_Y(X)$ (the subspace of $\mcC^{\infty}(X)$ consisting of fiber constant functions on $\p X$) preserves $\mcD$. Thus, $\mcD$ can be characterized locally away from $\p X$, plus locally in $Y$ near $\p X$ (i.e. near $\p X$ the domain does not have a local characterization, but it is local in the base $Y$, so the non-locality is in the fiber $Z$.)

Since we will be working on the manifold $M$ throughout the paper, we will need that the metric $dt^2 + dg$ gives rise to the Hilbert space $L^2_g(M;TX).$  More precisely, we will now describe the notation used for this inner product.
\begin{definition} \label{def:ip} \todo{label(def:ip)}
Let $u,v \in \dot{C}^{\infty}(M;TX)$. Suppose $TX$ is given a local trivialization and $u = (u^i)$, $v=(v^j)$ with respect to the trivialization. Denote $g= (g_{ij})$ as the matrix corresponding to our given edge metric; the fiber inner product takes the form
$$ (u,v)_g = \sum_{ij} g_{ij}u^i \bar{v}^j .$$
 By a standard partition of unity argument, the global inner product of sections of $TX$ takes the form
$$\langle u, v \rangle = \int_M (u,v)_g \ dg dt
 .$$
\end{definition}
This gives rise to a dual pairing between $H_e^{m,l}(M;TX)$ and $(H_e^{m,l}(M;TX))^*$. This convenient
choice of inner product makes $P$ formally self adjoint, so we indeed have
\begin{equation}
P = P^* \ \text{formally}
\end{equation}

Adopting the conventions in \cite{MVWEdges}, we also
write $\tmcD$, etc. for the analogous space on $M$:

\begin{definition}
$||u||^2_{\tmcD} = ||D_tu||^2_{L^2(M)} + \int B(u,u) dt + ||u||^2_{L^2(M)}.$ We also
write $\tmcD([a,b])$ for the space with the same norm on $[a,b] \times X$.
\end{definition}

A further localization of $\mcD$ will be useful.
\begin{definition}
For $u \in \mcC^{-\infty}(X)$, we say $u \in \mcD_{\text{loc}}$ if $\phi u \in \mcD$ for all
fiber constant $\phi \in \mcC_c^{\infty}(X)$. Similarly, for $u \in \mcC^{-\infty}(X)$, we say that
$u \in \mcD'_{\text{loc}}$ if $\phi u \in \mcD'$ for all fiber constant $\phi \in \mcC_c^{\infty}(X)$.
The localized domains on $M$ are defined analogously along with powers of $L$.
\end{definition}

 We will also use Melrose's $b$-calculus since certain \emph{admissable} elastic wave equation solutions will naturally lie in a $b$-based Sobolev space. Their definition and properties are in \cite[Section 6]{MVWEdges} for the scalar case, and in \cite{VasyForms} for the vector bundle case. We will use the notation in those papers such as $$H^m_{b,\mcD,c}(M;TX),H^m_{b, \mcD, loc}, \WF_b, \ etc.$$  Thus, we say a solution to the elastic equation is \emph{admissable} if it lies in some $H^m_{b,\tmcD,loc}(M:TX)$. We will mostly use them to prove finite propagation speed with respect to such spaces as in Section \ref{thm:EnergyEstimate}.

{\bf Also, to avoid cluttering with notation, we will often omit the bundle $TX$ when describing spaces such as $\Psi_e, H^{s,m}_e$,etc... when there's no risk of confusion. In fact, the nonscalar nature of all these spaces will only be relevant in a few key places which we will indicate explicitly.} One such place where it is particularly relevant is when describing adjoint of pseudodifferential operators, which we prepare to do in the next sections.

\section{Adjoints}
An important ingredient in the proof of diffraction will be to place our operator $P$ in a model form, as well as using positive commutator arguments. Such analysis invariably uses the $L^2_g$-based adjoints of pseudodifferential operators, and since our operators are acting on sections of a vector bundle, these are no longer so trivial.

\subsection{Adjoints of Edge Pseudodifferential Operators}
To begin, consider an arbitrary $A \in \Psi^{\infty}_e(M;TX)$.  Picking a trivialization of $TX$, the principal symbol, $a= \sigma_e(A)$ is an $n \times n$ matrix of symbols. As is known (see \cite{Hor} for example), if we denote $A^{\dag}$ as the $L^2(M;TX)$ adjoint of $A$ with the Euclidean inner product on the fibers of $TX$ using a trivialization, and integration with respect to the metric, then
$$ \sigma_e(A^{\dag}) = a^{\dag}, $$
 where in local coordinates, $a^{\dag}$ is the conjugate transpose of the matrix $a$. So with the notation in Definition \ref{def:ip}, we compute
 $$ \langle u, Av \rangle = \int_M (u,Av)_g dtdg = \int_M (gu,Av) dtdg
 = \int_M (g^{-1}A^{\dag}gu,v)_g dtdg. $$
 Hence, we have
 $$ A^* = g^{-1}A^{\dag}g \text{ and }\sigma_e(A^*)= g^{-1} a^{\dag} g.$$
 Thus, if one is not dealing with a principally scalar operator $A$, it's adjoint is more complicated than just being the conjugate transpose of its principal symbol. However, if $A$ is principally scalar, then things are much nicer. Indeed, following Vasy in \cite{VasyForms}, for $A_0 \in \Psi_e(M)$, let $A_0^{\dag}$ denote the $L^2_g(M)$-adjoint of $A_0$ with principal symbol $a_0$, and let $A = A_0 \otimes \Id$. In this case
 \begin{equation}\label{eq:scalarAdjoint} \todo{ label(eq:scalarAdjoint)} \\
 A^* = g^{-1}(A_0^{\dag} \otimes \Id) g \text{ and } \sigma_e(A^*) = \bar{a}_0 \otimes \Id.
\end{equation}

An analogous discussion applies for operators in $\Psi_e^{\infty,-\infty}$.
In particular, (\ref{eq:scalarAdjoint}) implies that when $A$ has a real, scalar principal symbol, then for $A \in \Psi^{m,l}_e(M)$, we in fact have
$$ A^* - A \in \Psi_e^{m-1,l}(M;TX).$$
Thus, we will constantly use this fact that when an operator has a real scalar, principal symbol, then it differs by it's adjoint by an operator of lower order, which is usually nonscalar.

\section{Coisotropic regularity and Coinvolutivity}
In this section, we will make the formal definitions of a distribution being coisotropic or coinvolutive, which was described in only loose terms in the introduction.

\subsection{ Coisotropic Regularity and Coinvolutive Regularity}\label{sec:Coiso Regularity}
First, we cite an important theorem in \cite{MVWEdges} adapted to our notation whose proof is identical as the one presented in that paper:
\begin{theorem}(\cite[Theorem 4.1]{MVWEdges}
Away from glancing rays, the sets $\mcF^p$ and $\mcF^s$ are coisotropic submanifolds of the symplectic manifold $\e T^*M$, i.e. each contains its symplectic orthocomplement.
\end{theorem}

We also make the following definition, taken from \cite{MVWEdges}, but changed only slightly to distinguish between $p$ and $s$ bicharacteristics, and to allow certain nonscalar terms. We first give the definition corresponding to $p$-bicharacteristics since the $s$-rays are similar.

Fix an arbitrary open set $U \subset \e T^*M$ disjoint from rays meeting $x=\xi=0$, i.e. away from the set $\mathcal{G}^p$.

\begin{definition}\label{def:coisotropy}

\begin{enumerate}
 \item[(a)] Let $\Psi_e(U;TX)$ be the subset of $\Psi_e(M;TX)$ consisting of operators $A$ with $\WF'_e(A) \subset U$.
\item[(b)] Let $\mcM^p$ denote the module of pseudodifferential operators in $\Psi_e^1(M;TX)$ given by
    $$\mcM^p = \{ A \in \Psi_e^1(U;TX): \sigma(A)= a_0 \otimes \Id
    \text{ with } a_0 |_{\dmcF^p} = 0\}.$$
\item[(c)] Let $\mcA_p$ be the algebra generated by $\mcM^p$, where we require its elements to have scalar principal symbol, with $\mcA_p^k = \mcA \cap \Psi_e^k(M;TX).$ Let $\bf{H}$ be a Hilbert space on which $\Psi_e^0(M;TX)$ acts, and let $K \subset \e T^*(M)$ be a conic set.
\item[(d)] We say that $u$ has $p$-coisotropic regularity of order $k$ relative to $\bf{H}$ in $K$ if there exists $A \in \Psi^0_e(M)$, elliptic on $K$, such that $\mcA_p^kAu \subset \bf{H}$.
\item[(e)] We say that $u$ is coinvolutive of degree $k$ relative to $\bf{H}$ on $K$ if there exists $A \in \Psi^0_e(M;TX)$, elliptic on $K$, such that $Au \subset \mcA^k\bf{H}$. We say $u$ is coinvolutive relative to $\bf{H}$ on $K$ if it satisfies the condition to some degree.
\end{enumerate}
\end{definition}

We also have the following important lemmas taken from \cite{MVWEdges}, tweaked in order to account apply to the vector bundle case, but whose proofs nevertheless remain the same:

\begin{lemma}(adapted from \cite[Lemma 4.4]{MVWEdges})\label{lemma:Module}
$\mcM^p$ is closed under commutators and is finitely generated in the sense that there exists finitely many $A_i \in \Psi_e^1(M;TX)$, $i=0,1,\dots,N$, with scalar principal symbol, $A_0 = \Id$, such that
$$\mcM^p = \{ A \in \Psi^1_e(U;TX): \exists Q_i \in \Psi^0_e(U) \text{ with scalar principal symbol}, \
A = \sum_{j=0}^N Q_i A_i \}$$
Moreover, we make take $A_N$ to have principal symbol $|\tau|a_N \otimes \Id = |\tau|^{-1}q_p \otimes \Id$,
and $A_i$ to have principal symbol $|\tau|a_i \otimes \Id$ with $da_i(q) \in T^{*,-}_q(\e S^*M)$ for $q \in \p \dot{\mcF}^p$, where we used the notation of Remark \ref{rem: eigenvalues for stable/unstable}.
\end{lemma}

We also have as in \cite{MVWEdges}
\begin{lemma}
If $A_i$, $0 \leq i \leq N$, are generators for $\mcM^p$ in the sense of Lemma \ref{lemma:Module} with $A_0 = \Id$, then
$$
\mcA_p^k = \left\{ \Sigma_{|\alpha| \leq k} Q_{\alpha} \Pi_{i=1}^N A_i^{\alpha_i}, Q_{\alpha} \in
\Psi_e^0(U;TX) \text{ with scalar principal symbol } \right\}
$$
where $\alpha$ runs over multiindices $\alpha:\{1,\dots,N\} \to \mathbb{N}_0$ and $|\alpha| = \alpha_1 + \dots \alpha_N.$
\end{lemma}
\begin{rem}(\cite[Remark 4.6]{MVWEdges}) The notation here is that the empty product is $A_0 = \Id$, and the product is ordered by ascending indices $A_i$. The lemma is an immediate consequence of $\mcM^p$ being both a Lie algebra and a module; the point being that products may be freely rearranged, module terms in $\mcA^{k-1}$.
\end{rem}

As in \cite[Section 6]{MVWCorners} we will use some important facts about coisotropic manifolds. Indeed, away from $\{ x = 0\}$, we may always (locally) conjugate by an FIO to a convenient normal form: being coisotropic, locally $\mcF^{p/s}$ can be put in a model form $\zeta = 0$ by a symplectomorphism $\Phi$ in some canonical coordinates $(y,z,\eta,\zeta)$, see \cite[Theorem 21.2.4]{Hor} (for coisotropic submanifolds one has $k=n-l$, $\text{dim}(S)=2n$, in the theorem). We state the result as a lemma:

\begin{lemma}\label{lemma:module FIO normal form}
We may quantize $\Phi$ to a FIO $T$, elliptic on some neighborhood of $w \in \mcF^{p/s}$ to have the following properties
\begin{enumerate}
\item[(i)] $u$ has coisotropic regularity of order $k$ (near $w$) with respect to $H^s$ if and only if $D_z^{\g}Tu \in H^s$ whenever $|\g|<k$.
    \item[(ii)] $u$ is coinvolutive of order $k$ (near $w$) with respect to $H^s$ if and only if $Tu \in \Sigma_{|\g|\leq k} D_z^{\g}H^s$.
\end{enumerate}
\end{lemma}

The key additional information we need is a lemma in \cite{MVWEdges}.
\begin{lemma}(Lemma 4.7 in \cite{MVWEdges})
For $l = 1, \dots, N-1,$
\begin{equation}\label{eq:module commuting Q_p with module}
x^2i[A_l,Q_p] = \Sigma_{j=0}^N C_{lj}A_j, \ C_{lj} \in \Psi_e^1, \ \sigma(C_{lj})|_{\p \dmcF^p} = 0 \text{ for }j \neq 0.
\end{equation}
(All operators above are principally scalar)
\end{lemma}

The same definition and lemmas apply for the $s$-bicharacteristics. We now introduce the following notation to describe the Hilbert spaces of distributions which have coisotropic regularity of a certain degree. From now on we will use the notation $\mcM$ or $\mcA$ to refer to either the $p$ or $s$ versions of the module, and it will be clear in context which one we are referring to.

\begin{definition}\label{def:coisoSpace}
For the set $U$ introduced in Definition \ref{def:coisotropy}, we can define the space of distributions which have coisotropic regularity of degree $k$ w.r.t. $H^{m,l}_e$ on microlocally inside $U$.
\begin{align*}
 &\mathcal{I}^kH_e^{m,l}(U,\mathcal{M}^{p/s}(U)) := \{ u \in H^{-\infty,l}_e(M):
\mcA^k_{p/s} u \in H^{m,l}_e \}.
\end{align*}
(here, $\mcA^k_{p/s}$ really stands for $\mcA^k_{p/s}(U)$.)
\end{definition}

We will need one final piece of information in order to analyze regularity with respect to $\mcM^{p/s}$ in the later sections. Let $Q_{p/s} \in \Psi_e^{2,-2}(M)$ be a quantization of the edge symbols $q_{p/s}$ corresponding to the $p/s$ waves. Let $W_{m,l} \in \Psi_e^{m,l}(M)$ have symbol
$w_{m,l} = |\tau|^mx^l$ near $\Sigma$ where $\tau \neq 0$, and has an arbitrary smooth extension elsewhere. As shown in \cite[Section 4]{MVWEdges}, we have the following crucial computation
\begin{equation}\label{eq:Modules [W,Q_p] calculation}
i[W_{m,l},Q_{p/s}] = W_{m+1,l-2}C_0 \text{ where }C_0 \in \Psi_e^{0,0}(M), \
\sigma_e(C_0)|_{\p \dot{\mcF}^{p/s}} = -(m+l)c_{p/s}^2\hat{\xi}.
\end{equation}

Finally, we will need some theorems that describe the propagation of coisotropic regularity on $M^o$, where we do not deal with boundaries. Such results are old and well-known. For example Dencker in \cite{Dencker} shows how inside $M^o$ one may find a parametrix to reduce $P$ to a scalar wave operator in order to invoke H\"{o}rmander's theorem (cf. \cite[Theorem 6.1.1]{HorDuis}) for standard propagation of singularities. Thus, one has the following, stated in a similar fashion to \cite[Section 6]{MVWCorners}:

If $K \subset M^o$ is compact, then there is $\d>0$ such that if $p \in S^*_KM^o$ and $\g_{p/s}$ is a $p/s$-bicharacteristic, then for $s \in (-\d,\d), \g(s) \in M^o$. As $s$ is equivalent to $t$ as a parameter along a bicharacteristic, we have the following result similar to \cite[Corollary 6.13]{MVWCorners}.

\begin{cor}\label{cor:module interior coiso for forward solution}
Suppose $K \subset M^o$ is compact. Suppose that $f$ is coisotropic, resp. coinvolutive (on the coisotropic $(\mcF^{p/s})^o$), of order $k$ relative to $H^{m-1}$, supported in $t>T$. Let $u$ be the unique solution of $Pu = f$, supported in $t>T$. Then there exits $\d_0>0$ such that $u$ is coisotropic, resp. coinvolutive ((on the coisotropic $(\mcF^{p/s})^o$)), of order $k$ relative to $H^m$ at $p\in S^*_KM^o$ if $t(p) < T+ \d_0$.

The analogous statements hold if $f$ is supported in $t<T$, and $u$ is the unique solution to $Pu=f$ supported in $t<T$, by virtue of vanishing there, except one needs $t(p) > T - \delta_0$ in the above notation.
\end{cor}

(We remark that when we write $p/s$ in the statement of the proposition as well as the following corollaries, we mean that the statement holds either in the $p$ case where each $p/s$ is replaced by $p$, or in the $s$ case where each $p/s$ is replaced by $s$.)

(We also remark that one could certainly give an alternative proof of this proposition by positive commutator arguments similar to, but much easier than, those used for propagation of coisotropic edge regularity in the following several sections.)

Of course, what happens to coisotropic regularity and coinvolutive regularity when bicharacteristics reach $\p M$ is of considerable interest, and occupy the remainder of this paper.

Having described the notation and definitions to be used, we proceed with the first piece in the proof of the main theorem in the next section.

\section{Inner Products and Trivializations}

Throughout the sections, we will be working in local coordinates of the manifold where the bundle $TX$ is trivialized as well. Now the metric inner product $\langle \cdot, \cdot \rangle_g$ on $L^2(\mathbb{R}_t \times X;TX)$ certainly does not depend on which trivialization we choose for $TX$, but when we express our operators as matrices, we are certainly fixing a local frame for $TX$ in which our matrices are expressed. For example, the matrix we wrote for the principal symbol of the elastic operator $p = \sigma_{e,2,-2}(P)$ was expressed in the coordinate local frame $\{ \p_x, \p_{y_i}, \p_{z_j}\}$ within some coordinate patch. However, it will be computationally convenient to use orthonormal frames to express our operators so that if $Q \in \Psi_e^{m,l}(M;TX)$ is a self-adjoint operator with respect to the metric inner product, then one may find an operator $S \in \Psi_e^{0,0}(M;TX)$ such that $S^*QS$ will have a principal symbol which may be written as a diagonal matrix with respect to an orthonormal trivialization. We will now give the details of such a construction in a general abstract setting.

Let us first describe this process in an abstract simplified setting where $(X,g)$ is a compact, $n$-dimensional Riemannian manifold without boundary with metric $g$.
Let $E$ be a vector bundle of rank $N$, endowed with an inner product $( \cdot, \cdot )_E$ and consider operators $P \in \Psi^m(X;E)$. Again, $\l \cdot, \cdot \r$ denotes the metric inner product where $(\cdot,\cdot)$ is an inner product on $E$ and integration is with respect to $dg$. Observe that the principal symbol $p =\sigma_m(P)$ is an element of $S^m(T^*X; \End(\pi^*E))$ where $\pi : T^*X \to X$ is the projection. Now, fix a coordinate patch $O \subset X$ where we have $x, \xi$ denoting local coordinates on $T^*X$ and an orthonormal trivialization
$$ \hE = \hE_O = \{ \he_1(x), \dots, \he_n(x) \} \text{ for } x \in O.$$
(Orthonormal means $(\he_i(x), \he_j(x)) = \delta_{ij} \ \forall x \in O$.)
Hence, for a distribution $u = \Sigma_k u_k \he_k(x) \in \mcD_c'(O;E)$, $u_k \in \mcD_c'(O)$, the action of $P$ on $u$ is given by
$$ Pu = P( \Sigma_{k=1}^n u_k \he_k(x)) = \Sigma_i\Big(\Sigma_j P_{ij}(u_j)\Big) \he_i(x)$$
for some $P_{ij} \in \Psi^m(O)$ scalar pseudodifferential operators. Thus, in the chart $O$, to say that $P$ is represented by a matrix $(P_{ij})$ of pseudodifferential operators, we are implicitly using a choice of trivialization $\hE$. Thus, it might be more correct to write $P^{\hE} = (P_{ij})$ to make this choice more explicit. Likewise the principal symbol $p$ depends on this trivialization and we may write
$$ [p]_{\hE} = (p_{ij}) $$
where $p_{ij}$ are the principal symbols of $P_{ij}$.
Now suppose that $P$ is self-adjoint so that $p=p(x,\xi)$ is self-adjoint as well in the sense that
$$( p v, w)_E = (v,pw)_E, \ \forall v,w \in E_x .$$
Thus, identifying $E$ with $O \times \mathbb{C}^N$ inside $O$, then by the spectral theorem, we may find a orthonormal basis of eigenvectors (which we assume to be smooth inside this local patch) $$\mcE = \mcE_O = \{ \be_1(x,\xi), \dots, \be_n(x,\xi) \}, \ \be_j \in C^{\infty}(T^*O;\mathbb{C}^N)$$
such that
$$ p(x,\xi) (\be_j) = \lambda_j \be_j \text{ for scalar functions } \lambda_j \in C^{\infty}(T^*O).$$
Indeed, certainly for each fixed $(x,\xi) \in T^*O$, such $\lambda_j(x,\xi)$ will exist, but the smoothness is less obvious and needs additional assumptions. However, in the case of the isotropic elastic operator, the eigenvalues will be smooth, so we just assume here that the $\lambda_j(x,\xi)$ are in $C^{\infty}(T^*O)$.

Now, define the linear operator $k: \pi^*_{x,\xi}E \to \pi^*_{x,\xi}E$ by $k(\he_j) = \be_j$ on the set $T^*O$. Indeed we may view $k$ as an element of $S^0(T^*O; \End(E))$. Here, $\text{End}(E)$ is the vector bundle with fiber at $(x,\xi)$ consisting of linear maps from $E_x$ to $E_x$. Since the $\he_i$ are orthonormal, then $k$ is actually orthogonal. Indeed,
for $ v = \Sigma_i v_i \he_i$ and $ w = \Sigma_i w_i \he_i$
$$ ( kv, kw )_E = v_i w_j ( \be_i, \be_j )_E
= v_i w_j \delta_{ij} = v_i w_j ( \he_i, \he_j )_E
= ( v,  w )_E.$$
Hence, certainly $k^{-1}$ exists and we may diagonalize $p$ by setting
$$ \tilde{p} = k^{-1} p k \Rightarrow \tilde{p}(\he_j) = \lambda_j \he_j
\text{ for each }j.$$ Thus, we have
$$ [\tilde{p}]_{\hE} = \text{diag}( \lambda_1, \dots, \lambda_n).$$

The entire discussion above applies in the exact same way to edge-pseudodifferential operators since the discussion was entirely microlocal. Thus, translating the discussion above to the setting of this paper means our underlying manifold is $(M,g)$, the vector bundle is $E = TX$ and the inner product on $TX$ is the metric inner product $(\cdot,\cdot)_g$. Let us fix any orthonormal frame $\hE$ for $TX$. Since $p = \sigma_{e,2,-2}(P)$ is symmetric with respect to $(\cdot,\cdot)_{g}$, there is a symbol $s \in S^0( \e T^*M; \End(TX))$ such that locally, for any $(t,x,y,z,\tau,\xi,\eta,\zeta) \in \e T^*M$, $$s(t,x,y,z,\tau,\xi,\eta,\zeta) \in \End(TX) \text{ is orthogonal with respect to } (\cdot, \cdot)_g.$$
Then the adjoint $s^*$ of $s$ with respect to this inner product is the inverse of $s$ and $s^*ps$ is a diagonal matrix with respect to the trivialization $\hE$.  {\bf \emph{Thus, for sections \ref{sec:partial elliptic regularity}, \ref{sec: edge propagation}, and \ref{sec: cosio regularity} whenever we trivialize $TX$ we will be using this orthonormal trivialization where all vectors and matrices are expressed with respect to this trivialization without explicitly saying so.}}

We now have the tools necessary to conjugate the elastic operator $P$ into a model form.

\section{Constructing the conjugated elastic operator}

We have arranged things on a principal symbol level, but for our purposes, we will need information beyond the principal symbol. Again we let $O \subset M$ be a neighborhood inside of a local chart where $TX$ is trivialized according to the trivialization described above. Then let $S_0$ be a quantization of $s$. Since $s$ is orthogonal, one has $s^{-1} = s^*$ (with the adjoint taken with respect to $(\cdot,\cdot)_g$ inner product on $TX$). Thus,
\begin{align}\label{eq:trivial R_-1 is symmetric}
&s^*s = ss^*= \Id
\nonumber \\
&\Rightarrow S_0^*S_0 = \Id + R_{-1}, \quad \text{for some }R_{-1} \in \Psi_e^{-1,0}(O)
\end{align}
 by using the edge calculus and since $\sigma_{e,0,0}(S_0^*) = s^*$. Thus, $S^*_0$ is certainly not a parametrix for $S_0$ but if one replaces $S_0$ by an asymptotic sum $S_0 + \sum_j S_{-j}$ where $S_{-j}$ is of order $-j$ then one can solve algebraic equations for the principal symbols of $S_{-j}$ (analogous to the microlocal square root construction argument as in ) so that this sum becomes a true parametrix. We leave out the details and merely state the lemma.

 \begin{lemma}\label{lemma:trivial S^* parametrix for S}
 With the above notation, one may find an operator $S \in \Psi_e^{0,0}(O;TX)$ such that $SS^* - \Id \in \Psi_e^{-\infty,0}(O;TX)$, and $\P := S^*PS \in \Psi_e^{2,-2}(O;TX)$ has principal symbol
 $$ \sigma_{e,2,-2}(\P) = \text{diag}(q_p, q_s \Id_{n-1}).$$
 \end{lemma}

\section{Partial Elliptic Regularity}\label{sec:partial elliptic regularity}
The first step to proving Theorem \ref{thm:Diffraction} is to show that coisotropic regularity is preserved along $\dmcF^{p/s}_I$, and that coisotropic regularity along $\dmcF^p_{I,\alpha}$ implies coisotropic regularity along $\dmcF^p_{O,\alpha}$ (one may look at Corollary \ref{Cor:PropOfCoisotropy} for an exact statement). The key to proving this is to break up a solution $u$ of the equation $Pu=0$ as $u = u_p + u_s$ corresponding to the $q_p$ and $q_s$ eigenspaces of $\sigma_e(P)$ respectively. The point is that along $\mcF^p$, which is inside the $p$ characteristic set $\Sigma_p$, $u_s$ will solve an elliptic equation, and so we will obtain elliptic estimates for $u_s$. We will show in the later sections how this will allow us to analyze the piece $u_p$ separately. We now proceed to give a precise description of $u_s$ and $u_p$ along with the elliptic estimates that follow.

First, we have several remarks regarding notation. We will continue to suppress the manifold and the bundle in the notation of various spaces to avoid cluttering when there's no risk of confusion. Also, unless specifically mentioned, all our operators will be assumed scalar unless mentioned specifically so that symbols $a$ are identified with $a \otimes \Id$ and scalar operators $A$ with $A \otimes \Id$.

To prove a propagation result, we will employ a positive commutator argument as done in \cite{MVWEdges}. One of the main advantages to working on $\e T^*M$ is that $P$ is naturally an edge operator, and we can put it into a model form by conjugating $P$ to be a diagonal operator. The principal symbol $p=\sigma_{e,2,-2}(P)$ is symmetric with respect to the metric inner product so we may find $s,s^* \in \eS^{0,0}, ss^* = \Id $, matrices of symbols with $s^*$ the adjoint of $s$ with respect to the metric inner product on $TX$, so that
$$ s^*ps = \text{diag}(q_p,q_s\Id_{n-1} ) :=\ti{p}$$
as explained in the previous section when we use an orthonormal trivialization of $TX$.
 Recall $q_p, q_s \in \e S^{2,-2}$ are the principal symbols corresponding to the pressure and shear waves.
 Now we let $S,S^*$ be quantizations of $s,s^*$ respectively in $\Psi_e^{*,*}$.
 As done in Lemma \ref{lemma:trivial S^* parametrix for S}, we may locally arrange that $S,S^*$ are inverses of each other modulo a low order error term we denote by $E \in \Psi_e^{-\infty,0}.$

 For now, let us work with the conjugated operator
 \begin{equation}\label{eq:conjOp}\todo{label(eq:conjOp)}
 \tilde{P} = S^*PS,
 \end{equation}
  such that
\[  \sigma_{e,2,-2}(\tilde{P})=
\left[ \begin{array}{cc}
q_p & 0  \\
0 & q_s\Id_{n-1}
\end{array} \right]. \]

 Now, let $\Pi_p,\Pi_s \in \Psi_e^{0,0}(M;TX)$, denote the projections to the $q_p,q_s$-eigenspaces of $\tilde{P}$ respectively. In fact, inside a local chart where all bundles are trivialized, we can write these down explicitly for $v \in \mcC^{-\infty}(M;TX)$:
 \[ \Pi_pv = \Pi_p  \left( \begin{array}{c}
v_1 \\
\vdots \\
v_n
\end{array} \right)
= \left( \begin{array}{c}
v_1 \\
0 \\
\vdots \\
0
\end{array} \right), \ \text{and} \ \
 \Pi_s  \left( \begin{array}{c}
v_1 \\
\vdots \\
v_n
\end{array} \right)
= \left( \begin{array}{c}
0 \\
v_2 \\
\vdots \\
v_n
\end{array} \right).
\]
With this explicit form, we clearly see that $[\P, \Pi_{p/s}]$ is lower order since the principal symbol of $\P$ is diagonal, that is
\begin{equation}
[\P, \Pi_{p/s}] \in \Psi_e^{1,-2}(M).
\end{equation}
For a distribution $v$, denote
 $$ v_p = \Pi_p v, \ \ v_s = \Pi_s v.$$
Now that we have introduced these projections, let us introduce another piece of notation that will keep  equations from becoming cluttered. In local charts where bundles are trivialized, we will view elements of $\mcC^{-\infty}(M;TX)$ as vectors with two components corresponding to the $p$ and $s$ projections above. That is, for $ v \in \mcC^{-\infty}(M;TX)$ we will write
$v = \begin{bmatrix} v_I \\ v_{II} \end{bmatrix}$ where $v_{I}$ is a vector with $1$ component corresponding to the $p$-eigenspace and $v_{II}$ is a vector with $n-1$ components corresponding to the $s$-eigenspace. Correspondingly, we may write operators $E \in \Psi^{\infty,-\infty}_e(M;TX)$ within a local chart as
$$ E = \begin{bmatrix} E_I & E_{II} \\ E_{III} & E_{IV} \end{bmatrix}, $$
where $E_I$ is a $(1 \times 1)$ matrix, $E_{II}$ a $(1 \times n-1)$ matrix,
$E_{III}$ a $(n-1 \times 1)$ matrix, and $E_{IV}$ a $(n-1 \times n-1)$ matrix. We will again use the convention that if $E_{j}$ is scalar, then we will write it as $E_j$ rather than $E_j \otimes \Id$ when it is clear in the context.

To proceed, take $\alpha \in \mathcal{H}^p$ and let $O$ be a local chart containing the projection of $\alpha$ to $M$ where all bundles are trivialized. Suppose we have $\tilde{P}u \in H^{m-1,l-2}_e$ and $u \in H^{m,l}_e$ microlocally near $\alpha$. Thus, for any $\tilde{G}' \in \Psi_e^{0,0}$ that is elliptic at $\alpha$, has Schwartz kernel that is compactly supported in $O \times O$, and is microlocally supported close to $\alpha$, we have $\tilde{G}'u \in H_e^{m,l}$. Observe that
 $$ \tilde{G}'\tilde{P}u_s = \tilde{G}'\Pi_s \tilde{P}u + \tilde{G}'[\tilde{P},\Pi_s]u.$$
 Since $\tilde{G}'[\tilde{P},\Pi_s] \in \Psi^{1,-2}_e$ and it is microlocally supported near where $u\in H_e^{m,l}$, then $\tilde{G}'[\tilde{P},\Pi_s]u \in \Psi_e^{m-1,l-2}$.
 We thus conclude that
 \begin{equation}\label{eq:ellEst}
 \tilde{G}'\tilde{P}u_s \in \Psi_e^{m-1,l-2}.
 \end{equation}
This proves part of the following proposition, which gives the main semi-elliptic estimate:

 \begin{prop}\label{prop:semiElliptic}
 Suppose $\alpha \notin \textnormal{WF}_e^{m,l}(u)$, $\tilde{P}u \in H_e^{m-1,l-2}$ microlocally near $\alpha$, and $\alpha \in \Sigma_p \setminus \Sigma_s$. Let $G \in \Psi_e^{0,0}$ elliptic near $\alpha$ with a Schwartz kernel compactly supported in $O \times O$, such that $\WF'(G) \cap \Sigma_s = \textnormal{WF}'(G)\cap \textnormal{WF}_e^{m,l}(u) = \emptyset$. Then $\alpha \notin \textnormal{WF}^{m-1,l-2}_e(\P u_s)$ and $Gu_s \in H_e^{m+1,l}.$ Moreover, the following estimate holds
 $$||Gu_s||_{H^{m+1,l}_e} \leq C ( ||G'\tilde{P}u||_{H^{m-1,l-2}_e} + ||G''u||_{H_e^{m,l}} + ||u||_{H_{e,loc}^{-N,l}}),$$
 for some $G',G'' \in \Psi_e^{0,0}$ elliptic on $\alpha$, microsupported in a neighborhood of $\alpha$, and whose Schwartz kernels are supported in $O \times O$.
 \end{prop}
 \begin{rem}
 This is the crucial place where we need $\Sigma_p$ and $\Sigma_s$ to be disjoint, since otherwise, we could never get such a semielliptic estimate. Without such an estimate, none of the theorems that we prove later would go through.
 \end{rem}

 \proof This is a symbolic exercise using that $q_s(\alpha) \neq 0$, together with the usual microlocal elliptic regularity. We also work in a local chart near the projection of $\alpha$ to the manifold where all bundles are trivialized. Indeed, take a parametrix $\tilde{P}_0^- \in \Psi_e^{-2,2}$ such that $\tP_0^-q_s(w, \e D_w)\otimes\Id - \Id \in \Psi_e^{0,0},$ and $\alpha \notin \WF'(\tP_0^-q_s(w, \e D_w) \otimes \Id - \Id)$. In fact, since $\WF'(G)$ is compact and disjoint from $\Sigma_s$, a parametrix may be chosen so that
 $$ \WF'(G) \cap \WF'(\tP_0^-q_s(w, \e D_w)\otimes\Id - \Id) = \emptyset.$$
 Then set $\tP^- \in \Psi_e^{-2,2}$ as
 \[
 \tP^- =\left[ \begin{array}{cc}
0 & 0  \\
0 & \tP_0^-
\end{array} \right]. \]
We then have
\[ \tP^- \tP - \Id = \left[ \begin{array}{cc}
-\Id & 0  \\
R' & R''
\end{array} \right] := R
 \in \Psi_e^{0,0} \text{ and } \WF'(G) \cap \WF'( R'' ) = \emptyset,\]
 where $R' \in \Psi^{-1,0}_e(M)$ and $R'' \in \Psi^{0,0}_e(M)$.

We now compute
\begin{equation}\label{eq:semiElliptic Gu_s term estimate}
Gu_s = G\tP^- \tP u_s + GR u_s = G\tP^- \tP u_s + GR'' u_s \in H_e^{m+1,l}
\end{equation}
since $\WF'(\tilde{G}\tP^-) \bigcap \WF^{m-1,l-2}(\tP u_s) = \emptyset$ by (\ref{eq:ellEst}), and $\WF'(GR'') \bigcap \WF^{m,l}(u_s) = \emptyset$.
 Now let $\tilde{G}',G'\in \Psi_e^{0,0}$ with $\tilde{G'}$ elliptic on $\WF'(G\tP^-)$ and $G'$ elliptic on the microsupport of $\tilde{G'}$,
 and such that one still has
 $$\WF'(G')\cap \WF^{m-1,l-2}_e(\tP u) =
 \WF'(G')\cap \WF^{m,l}_e(u)= \emptyset,$$
 with $\tilde{G}'$ having the same property.
 Then by the ellipticity of $G', \tilde{G}'$ in the aforementioned regions together with microlocal elliptic regularity, and equation (\ref{eq:ellEst}) gives
 \begin{align} \label{eq:semiElliptic tildePu_s estimate}
  ||G\tP^-\tP u_s||_{H^{m+1,l}_e} &\lesssim  ||\tilde{G'}\tilde{P}u_s||_{H^{m-1,l-2}_e} + ||u||_{H_{e,loc}^{-N,l}}
  \nonumber \\
  &\lesssim ||\tilde{G}'\Pi_s \tilde{P}u||_{H^{m-1,l-2}_e} + ||\tilde{G'}[\tilde{P},\Pi_s]u||_{H^{m-1,l-2}_e} + ||u||_{H_{e,loc}^{-N,l}}
  \nonumber \\
   &\lesssim ||G'\tilde{P}u||_{H^{m-1,l-2}_e} + ||G'u||_{H^{m,l}_e} + ||u||_{H_{e,loc}^{-N,l}},
\end{align}
where in the last inequality we again use the ellipticity of $G'$ on $\WF'(\tilde{G}')$ combined with microlocal elliptic regularity. We also have a similar estimate using microlocal elliptic regularity for the other term
\begin{equation}\label{eq:semiElliptic GR''u_s term}
 ||GR''u_s||_{H^{m+1,l}_e} \lesssim  ||G'u||_{H^{m,l}_e}+||u||_{H_{e,loc}^{-N,l}}.
 \end{equation}
Thus, combining (\ref{eq:semiElliptic Gu_s term estimate}) with the inequalities (\ref{eq:semiElliptic tildePu_s estimate}), and (\ref{eq:semiElliptic GR''u_s term}) gives the result of the proposition. $\Box$

The essential point in the proof was which characteristic set the point $\alpha$ belonged to. Thus, with essentially no changes except notation in the above proof, we get the analogous semi-elliptic estimates for $\alpha \in \Sigma_s$. We record it here for later use:

\begin{prop}\label{prop:semiElliptic s wave version }
 Suppose $\alpha \notin \textnormal{WF}_e^{m,l}(u)$, $\tilde{P}u \in H_e^{m-1,l}$ microlocally near $\alpha$, and $\alpha \in \Sigma_s \setminus \Sigma_p$. Let $G \in \Psi_e^{0,0}$ elliptic near $\alpha$ with a Schwartz kernel compactly supported in $O \times O$, such that $\WF'(G) \cap \Sigma_p = \textnormal{WF}'(G)\cap \textnormal{WF}_e^{m,l}(u) = \emptyset$. Then $\alpha \notin \textnormal{WF}^{m-1,l-2}_e(\P u_p)$ and $Gu_p \in H_e^{m+1,l}.$ Moreover, the following estimate holds
 $$||Gu_p||_{H^{m+1,l}_e} \leq C ( ||G'\tilde{P}u||_{H^{m-1,l-2}_e} + ||G''u||_{H_e^{m,l}} + ||u||_{H_{e,loc}^{-N,l}}),$$
 for some $G',G'' \in \Psi_e^{0,0}$ elliptic near $\alpha$, microsupported in a neighborhood of $\alpha$, and whose Schwartz kernels are supported in $O \times O$.
 \end{prop}

In the next section, we will discuss propagation into and out of the edge, which will rely on our semi-elliptic estimates.

 \section{ Edge Propagation  }\label{sec: edge propagation}
In this section, we describe the propagation of edge regularity into and out of the edge. First, let us state our main theorem for propagation to/away from the edge.

 \begin{theorem}\label{thm:propInOutEdge}
 Let $u \in H_e^{-N,l}$ be a distribution.
 \begin{enumerate}
 \item Let $m> l + f/2$. Given $\alpha \in \mathcal{H}^p_I$, if $u \in H^m$ microlocally on
 $\mathcal{F}^p_{I,\alpha} \setminus \partial M$ and $Pu \in H_e^{m-1,l-2}$ microlocally on $\bar{\mcF}^p_{I,\alpha}$, then $u \in H_e^{m,l'}$ microlocally at $\alpha, \forall l' <l.$

 \item  Let $m < l+f/2.$ Given $\alpha \in \mcH^p_O$, suppose $U$ is a neighborhood of $\alpha$ in $\left.\right.^eS^*|_{\partial M}M$ such that
 $\textnormal{WF}_e^{m,l}(u) \cap U \subset \partial \mathcal{F}^p_O$ and $\textnormal{WF}_e^{m-1,l-2}(Pu) \cap U = \emptyset$, then $\alpha \notin \textnormal{WF}_e^{m,l}(u).$

 \end{enumerate}
 \end{theorem}

 \begin{rem}\label{rem:coiso outgoing U_1 improvement}
 In $(2)$, although the theorem is stated with $ U \subset \left.\right.^eS^*_{\partial M}M$, such that
 $ (U \setminus \partial \mathcal{F}^p_O) \cap \WF^{m,l}(u) = \emptyset $ we may actually enlarge $U$ as follows. Since $\WF_e^{*,*}(u)$ is closed, we can find a small neighborhood $U_1 \subset \left.\right.^eS^*M$ of $U$ so that
 $ (U_1 \setminus \overline{\mathcal{F}}^p_O) \cap \WF_e^{m,l}(u) = \emptyset.$ We will refer to $U_1$ in the proof.
 \end{rem}

To prove this theorem, we will use the conjugated operator $\P$ introduced in equation (\ref{eq:conjOp}) of the previous section. Let $u$ be as in the statement of the above theorem. Since $S^*$ is a parametrix for $S$, we transform the problem by letting $\u = S^*u$ so that with Lemma \ref{lemma:trivial S^* parametrix for S}, one has
 \begin{equation}\label{eq:transformedPu}
  \P\u = S^*PS(S^*u) = S^*Pu + S^*PRu, \text{ for some } R \in \Psi_e^{-\infty,0}
  \end{equation}
  $$ \Rightarrow \WF_e^{m-1,l-2}(\P \u) = \WF_e^{m-1,l-2}(Pu)$$
by the ellipticity of $S$ and $S^*$.  Hence, $\P \u$ satisfies the same assumptions as $Pu$ in the statement of the above theorem. We will work with this transformed equation from now on.

We will prove Theorem \ref{thm:propInOutEdge} by a positive commutator argument for `radial points' as done in \cite[Section 11]{MVWEdges}. Thus, via an inductive argument which we justify later, we want to show that if $\alpha \notin \WF_e^{m-\teps,l}(u)$ then in fact $ \alpha \notin \WF_e^{m,l}(u)$ where $\alpha$ and $u$ satisfy the assumptions of the above theorem, and $0 < \teps \leq 1/2$. By the ellipticity of $S^*$, $\tu$ satisfies the same property, so we're essentially going to improve the Sobolev order of $\u$ microlocally by $\teps$ at each step. To proceed, let $\{A_{\d}\}_{\d \in [0,1]}$ be a scalar, uniformly bounded family in $\Psi_{e}^{m',l'}(M)$, and microlocally supported in a set disjoint from $\WF_e^{m-\teps,l}(u)$.
  If $m'$ is picked so that all the following pairings are finite, we compute
 \begin{align}\label{eq:coisoBreaking up the commutator into s and p parts}
 \langle A_{\d}\P\u,A_{\d}^*\u\rangle &- \langle A_{\d}\u, A_{\d}\tilde{P}\u \rangle
 \nonumber \\
 &= \langle [A^2_{\d},\tilde{P}]\u,\u \rangle
 \nonumber \\
 &= \langle [A^2_{\d},\tilde{P}]( \u_p + \u_s), \u_p + \u_s \rangle
 \nonumber \\
 &= \langle [A^2_{\d},\tilde{P}]\u_p, \u_p \rangle + \langle [A^2_{\d},\tilde{P}]\u_p, \u_s \rangle + \langle [A^2_{\d},\tilde{P}]\u_s, \u_p \rangle
 +\langle [A^2_{\d},\tilde{P}]\u_s, \u_s \rangle
 \end{align}
 Our strategy will be to estimate the first term on the right using a standard positive commutator estimate, while
 the other three terms will be estimated using Proposition \ref{prop:semiElliptic} derived in the previous section.
 This first term will be called the $pp$ term while the others will be called the $ps,sp,ss$ terms.
 We will estimate the terms containing $\u_s$ first since those are elementary now that we have Proposition \ref{prop:semiElliptic}.

 \subsection{ Estimating the $sp,ps,ss$ terms}
 The key is that the terms in (\ref{eq:coisoBreaking up the commutator into s and p parts}) involving $\u_s$ will be controlled by elliptic estimates. The exact orders are very important here.
  Since $A_{\d}$ is scalar, then $[A^2_{\d},\tilde{P}] \in \Psi_e^{2m'+1,2l'-2}(M)$ so in order to eventually show that $\tu$ is in $H_e^{m,l}$ microlocally near $\alpha$, we must have
 $$ 2m'+1 = 2m,\ \ \ 2l'-2 = 2(-l - \frac{f+1}{2})$$
 $$\Rightarrow m' = m-1/2, \ \ \ l' = -l- (f+1)/2.$$

 To justify pairings, we also have that
 $$\tilde{P}A^2_{\d} \in \Psi_e^{2m+1,-2l-(f+1)} \Rightarrow \tilde{P}A^2_{\d}:  H_e^{m-\teps,l} \to  H_e^{-m-1-\teps,-l-(f+1)} =  (H_e^{m+1+\teps,l})^*$$
 and
 \begin{equation}\label{eq:coisoMapping property of [P,A]}
 [\tilde{P},A^2_{\d}] \in \Psi_e^{2m,-2l-(f+1)} \Rightarrow [A^2_{\d},\tilde{P}]:  H_e^{m-\teps,l} \to  H_e^{-m-\teps,-l-(f+1)} =  (H_e^{m+\teps,l})^*
 \end{equation}

We'll state our bounds on the $ps,sp,ss$ terms as a proposition.
\begin{prop} \label{prop:coiso uniform boundedness of u_s terms}
Suppose that $K \subset U \subset \eS^*M$ with $K$ a compact neighborhood of $\alpha \in \mcH^p$, $U$ open, and let $O$ be a coordinate patch containing the projection of $\alpha$ to $M$ where all bundles are trivialized. Let $\mathbf{A} = \{ A^2_{\d} \in \Psi_e^{2m-1,2l'}(M): \d \in (0,1] \}$ a family with $\WF'_{e,L^{\infty}}(\mathbf{A}) \subset K$ which is bounded in $\Psi^{2m}_{e}(M)$ and has Schwartz kernels uniformly supported in $O \times O$.
Suppose that $\WF_e^{m-\teps,l}(u) \cap U = \emptyset$ and $\WF_e^{m-1,l-2}(Pu) \cap U = \emptyset$ for $0<\tilde{\epsilon} \leq 1/2$. Then the following pairings
are justified and remain uniformly bounded even as $\d \to 0$:
$$ \langle [A^2_{\d},\tilde{P}]\u_p, \u_s \rangle, \langle [A^2_{\d},\tilde{P}]\u_s, \u_p \rangle,
 \langle [A^2_{\d},\tilde{P}]\u_s, \u_s \rangle $$
\end{prop}

\proof All operators we mention here are assumed to have Schwartz kernels compactly supported in $O \times O$.
First, note that the assumption on $u$ and the ellipticity of $S$ imply $\u \in H^{m-\tilde{\epsilon},l}_e$ microlocally near $U$ as well. Thus, we trivially have $\u_p \in H^{m-\teps,l}_e$ microlocally near $U$ by microlocality of $\Pi_p$ and since $\Pi_p$ is $0$'th order.
 Pick $G_1 \in \Psi_e^{0,0}$ elliptic on $\text{WF}_{e,L^{\infty}}'(A_{\d})$, microsupported where $\u$ is in $H_e^{m-\teps,l}$. Hence, $G_1$ is elliptic on $\text{WF}_{e,L^{\infty}}'([A^2_{\d},\P]\Pi_p)$ as well. By microlocal elliptic regularity of $G_1$, the mapping property in (\ref{eq:coisoMapping property of [P,A]}), and continuity, we may estimate
 $$||[A^2_{\d},\tilde{P}]\u_p||_{H_e^{-m-\teps,-l-(f+1)}} \leq C( ||G_1\u||_{H_e^{m-\teps,l}} + ||\u||_{H_e^{-N,l}}).$$

 Also, if we let $G \in \Psi_e^{0,0}(M)$ elliptic at $\alpha$, microsupported inside $U$ such that
 \begin{equation}\label{eq:coiso I-G supportProp}
 \text{WF}'_e(I- G) \cap \text{WF}_{e,L^{\infty}}'(\textbf{A}) = \emptyset,
  \end{equation}
  then Proposition \ref{prop:semiElliptic}
implies $G\u_s \in H^{m-\tilde{\epsilon}+1,l}_e \subset H^{m+ \tilde{\epsilon},l}_e $ (the inclusion of Hilbert spaces is due to $\teps \leq 1/2$) and
 $$||G\u_s||_{H_e^{m + \teps,l}}  \leq ||G\u_s||_{H_e^{m - \teps +1,l}}
  \leq C( ||G'\tilde{P}\u_s||_{H_e^{m - \teps - 1,l-2}} +
 ||G''\u||_{H_e^{m-\teps,l}} + ||\u||_{H_e^{-N,l}}),$$
for some $G',G'' \in \Psi_e^0$ microsupported inside $U$ and elliptic at $\alpha$.
 To proceed, the microsupport property of $G$ in (\ref{eq:coiso I-G supportProp}) implies $(I-G)[A^2_{\d},\tilde{P}]$ is a uniformly bounded family in $\Psi_e^{-\infty,l}$, so we get the following uniform estimate (since $\Pi_{p/s}$ are clearly bounded operators between on edge Sobolev spaces)
 $$ |\langle (I-G^*)[A^2_{\d},\tilde{P}]\u_p, \u_s \rangle| \leq C||\u||^2_{H_e^{-N,l}}.$$

 Combining these estimates and using the $L^2_g$-dual pairing with Cauchy-Schwartz inequality, we obtain
 \begin{align*}
|\langle [A^2_{\d},\tilde{P}]\u_p, \u_s \rangle| &= |\langle (I-G^*)[A^2_{\d},\tilde{P}]\u_p, \u_s \rangle|
 + |\langle [A^2_{\d},\tilde{P}]\u_p, G\u_s \rangle|
 \\
 &\lesssim |\langle (I-G^*)[A^2_{\d},\tilde{P}]\u_p, \u_s \rangle|
 \\
 &\qquad \qquad \qquad + ||[A^2_{\d},\tilde{P}]\u_p||_{H_e^{-m-\teps,-l-(f+1)}} ||G\u_s||_{H_e^{m+\teps,l}}
 \\
 &\lesssim  ||G'\u||^2_{H_e^{m-\teps,l}} +  ||G_1\u||^2_{H_e^{m-\teps,l}} + ||G''\tilde{P}\u||^2_{H_e^{m - \teps - 1,l-2}} + ||\u||^2_{H_e^{-N,l}}.
 \end{align*}
The other terms are bounded analogously. $\Box$

 \subsection{ Reducing $\langle [A^2_{\d},\tilde{P}]\u_p, \u_p \rangle$ to the case of a scalar wave equation}

 In this part, since $A_{\d}$ is a scalar operator, we will express it as
 $$ A_{\d} = A_{0,\d} \otimes \Id$$
 to make some calculations more transparent, where $A_{0,\d}$ has all the properties already mentioned for $A_{\d}$ and is an honest edge pseudodifferential operator on scalar distributions.
 First, a careful calculation shows that
  \begin{equation}\label{eq:coiso writing our [A,P]}
  [A^2_{\d},\P] = \left[ \begin{array}{cc}
[A^2_{0,\d}, Q_p] & F'_{\d}  \\
F''_{\d} & F'''_{\d}
\end{array} \right] =
\left[ \begin{array}{cc}
[A^2_{0,\d}, Q_p ] & 0  \\
0 & F'''_{\d}
\end{array} \right] + \tilde{F}_{\d},
\end{equation}
with $\sigma(Q_p) = q_p \in x^{-2}S^2(\e T^*M)$, $F'_{\d},F''_{\d},\tilde{F}_{\d} \in \Psi_{e}^{m',l'-2}$ and $F'''_{\d} \in \Psi_{e}^{m'+1,l'-2}$ uniformly bounded families. We won't have control over $F'''_{\d}$ but that's irrelevant since we have
\begin{equation}\label{eq:coiso breaking up <[A,P]u_p,u_p>}
  \langle [A^2_{\d},\tilde{P}]\u_p, \u_p \rangle = \langle [A^2_{0,\d},Q_p]\u_p,\u_p \rangle + \langle \tilde{F}_{\d}\u_p, \u_p \rangle.
\end{equation}

Thus, since $\tilde{F}_{\d}$ is lower order and will be handled by inductive assumptions as we show later, we are reduced to a commutator with $Q_p$, whose principal symbol is the same as that of a scalar wave operator. Now we proceed to construct $A_{0,\d}$ in the same fashion as done in \cite{MVWEdges}. However, since $Q_p$ does not commute with $D_t$, we will need more care for the regularization argument.

\subsection{Constructing the family of test operators $A_{0,\d}$ }
Notice that the principal symbol of $Q_p$ is that of a scalar wave operator associated to the metric $g_p$ so that the construction of the test operator in \cite[Section 11]{MVWEdges} goes through almost verbatim.
For regularization, we define
\begin{equation}\label{eq:coiso Regularizer construction}
\varphi_{\delta}(y) = (1+\delta y)^{-\teps} \Rightarrow \varphi'_{\delta}(y) = -\teps\delta(1+\delta y)^{-1}\varphi_{\delta}(y),
\end{equation}
where $0< \teps <1/2$,

  We then have the following lemma whose proof is almost identical to what is done in \cite[Section 11]{MVWEdges} where no regularization is done, and in \cite[Section 8]{MWCones} where a regularization is done albeit a slightly different setting.

\begin{lemma}\label{eq:cosio Positive commutator}
With the notation above, and setting $m'= m-1/2$ and $l' = -l-(f-1)/2$, then for $\alpha \in \mcH^p_{I/O}$ and assuming either $m'+l'-\teps > 0$ or $m' + l' < 0$, we have may find an edge PsiDO $A_{\d,0}$ with the following properties:
\begin{equation}
 i[A^2_{\delta,0},Q_p] = \pm B_{\delta}^*B_{\delta} \pm \sum_j B^*_{\delta,j}B_{\delta,j} + E_{\delta} + C_{\delta} + K_{\delta} + F_{\delta}
\end{equation}

with

\begin{enumerate}
\item $A_{\delta,0} \in \Psi_e^{m-1,l'-1}$ for $\delta>0$, elliptic at $\alpha$, uniformly bounded in $\Psi_e^{m,l'-1}$, and $A_{\d,0} \to A_0$ in $\Psi_e^{m+\epsilon,l'-1}$ for any $\epsilon>0$.
    Moreover, given any conic neighborhood $U$ of $\alpha$, the family $\mathbf{A} = \{A_{\d,0}\}_{\d \in [0,1]}$ may be chosen so that $\WF_{e,L^{\infty}}'(\mathbf{A})$ is contained in $U$.
\item $B_{\delta},B_{\delta,j} \in \Psi_e^{m,l'-1}$ elliptic at $\alpha$, uniformly bounded,
$B_{\delta},B_{\delta,j} \in \Psi_e^{m-\teps,l'-1}, \delta > 0$ and
$B_{\delta} \to B, B_{\delta,j} \to B_j$ in $\Psi_e^{m + \epsilon,l'-1}$ for any $\epsilon >0$.

\item $E_{\delta} \in \Psi_e^{2m,l'-2}$ uniformly bounded, $E_{\delta} \in \Psi_e^{2m-2\teps,l'-2}$ for $\delta >0$,
and $\text{WF}_{e,L^{\infty}}'(E_{\delta}) \subset T^*M^o$ in the case $m'+l' - \teps > 0$. In the case $m' + l' < 0$, if $u$ satisfies the hypothesis of (2) in Theorem \ref{thm:propInOutEdge} with $U_1$ as in Remark \ref{rem:coiso outgoing U_1 improvement}, then one may arrange that
$$\text{WF}_{e,L^{\infty}}'(E_{\delta}) \subset U_1 \setminus
\overline{\mcF}^p_{O,\alpha} \text{ uniformly }$$
and such that $u \in H_e^{m,l}$ microlocally on $\WF'_{e,L^{\infty}}(E_{\d})$.

\item  $C_{\delta} \in \Psi_e^{-\infty,l'-2}$ uniformly bounded.

\item $K_{\delta} \in \Psi_e^{2m,l'-2}$ uniformly bounded, \\
$K_{\delta} \in \Psi_e^{2m-2\teps,l'-2}$ for $\delta >0$,
and $\text{WF}_{e,L^{\infty}}'(K_{\delta}) \bigcap \Sigma_p = \emptyset$

\item $F_{\delta} \in \Psi_e^{2m-1,l'-2}$ uniformly bounded.

\end{enumerate}
\end{lemma}

What is noteworthy for us is that
\begin{equation}\label{eq:coiso a' definition for in}
\sigma_e(A_{\d,0}) = (\pm \tau)^{m-1}x^{l'-1}\phi_{\d}(|\tau|^2)a'
\end{equation}
$a' \in S^0(\e T^*M)$ is a nonnegative symbol that is elliptic at $\alpha$, whose support may be made arbitrarily close to $\alpha$, and is supported near $\Sigma(P)$ so that $\tau \neq 0$ on $\supp(a')$. Also,
$$\sigma_e(B_{\delta}) = b_{\delta} = \sqrt{\pm(m'+l'-\teps\delta(1+\delta|\tau|^2)^{-1})c_0}w_{m,l'-1}\varphi_{\delta}
(|\tau|^2)\acute{a}$$
which is a well-defined symbol in $x^{l'}S_e^{m'+1/2}$, and is non-negative since $\pm c_0 > 0$ and $m' + l' + \teps >0$ on $\supp(a')$.

\subsection{ Proving propagation in/out of the edge}
We now have all the pieces to prove Theorem \ref{thm:propInOutEdge}.
First, we prove the key lemma which is a ``baby'' version of the theorem.
\begin{lemma}
\begin{enumerate}\label{lemma:intoTheEdge}

 \item
 Suppose $m > l + f/2 + \teps$ for some $0<\teps \leq 1/2 ,$ and $\alpha \in \mcH^p_I$.
\begin{equation}\label{lemma:coiso baby1}
 \alpha \notin \WF_e^{m-\teps,l}(u),\ \ \WF_e^{m,l}(u) \cap (\mathcal{F}^p_{I,\alpha} \setminus \partial M) = \emptyset,
 \ \ \alpha \notin \WF_e^{m-1,l-2}(Pu),
 \end{equation}
$$ \Rightarrow \alpha \notin \WF_e^{m,l}(u). $$

\item
Suppose $m < l + f/2 $ and $\alpha \in \mcH^p_O$, and $U$ is a neighborhood of $\alpha$ in $\e S^*_{\p M} (M)$.
\begin{equation}\label{lemma:cosio baby2}
 \alpha \notin \WF_e^{m-\teps,l}(u),\ \ \WF_e^{m,l}(u) \cap U \subset \p \mathcal{F}^p_{O,\alpha}, \ \ \alpha \notin \WF_e^{m-1,l-2}(Pu),
 \end{equation}
$$ \Rightarrow \alpha \notin \WF_e^{m,l}(u). $$
\end{enumerate}
\end{lemma}

\begin{lproof}
We will prove both parts simultaneously and point out the relevant differences. Again, we let $O$ be a local coordinate neighborhood of the projection of $\alpha$ to $M$ where all bundles are trivialized, and we assume all operators constructed here have Schwartz kernels compactly supported in $O \times O$. As already shown, $\u$ and $\P \u$ satisfy the same assumptions as $u$ and $Pu$ in this lemma. It will be convenient to let
$$ m' = m-1/2 \text{ and } l' = -l - (f-1)/2,$$
so the hypothesis of the lemma would say $m' + l' -\teps >0$ for $(1)$ and $m'+l'<0$ for $(2)$.
Also, as shown in the proof of Proposition \ref{prop:coiso uniform boundedness of u_s terms}, one trivially has $\alpha \notin \WF_e^{m-\teps,l}(\u_p)$ as well. Let $\mathbf{A} = \{A_{\d}\}_{\d \in [0,1]}$ be as in Lemma \ref{eq:cosio Positive commutator} where the microsupport of $A_{\d}$ may be made sufficiently close to $\alpha$ such that
\begin{equation}
\WF'_{e,L^{\infty}}(\mathbf{A}) \cap \WF_e^{m-\teps,l}(\u_p) = \emptyset
\text{ and } \WF'_{e, L^{\infty}}(\mathbf{A}) \cap \WF_e^{m-1,l-2}(\P \u) = \emptyset
\end{equation}
so that $A_{0,\d} \u_p$ remains bounded in $L^2_g(M)$ for $\d >0$.
By (\ref{eq:coisoBreaking up the commutator into s and p parts}) and (\ref{eq:coiso breaking up <[A,P]u_p,u_p>}), for $\d > 0$ so that the integration by parts and the pairings are justified, we have
\begin{align*}
 &\langle A_{\delta}\tilde{P}\u,A_{\delta}^*\u\rangle - \langle A_{\delta}\u, A_{\delta}^*\tilde{P}\u \rangle
 = \langle [A_{\delta}^2,\tilde{P}]\u, \u \rangle \\
 &= \langle [A_{0,\delta}^2,Q_p]\u_p, \u_p \rangle + \langle \tilde{F}_{\delta}\u_p, \u_p \rangle
  + \langle [A_{\delta}^2,\tilde{P}]\u_p, \u_s \rangle + \langle [A_{\delta}^2,\tilde{P}]\u, \u_s \rangle
\end{align*}
In fact, it is already shown in \cite[Equation (7.17)]{MVWCorners} along with the proof in that paper, that indeed, the first equality above does hold.

Thus, applying
Lemma \ref{eq:cosio Positive commutator} to the computation above, where the $B_{\delta,j}$ terms may be ignored since they have the same sign in front of them as the $B_{\d}$ term, we get the following estimate:
\begin{align}\label{eq:coiso estimating Bu_p}
 ||B_{\delta}\u_p||^2 &\leq |\langle F_{\delta}\u_p,\u_p \rangle| + |\langle A_{\delta}\tilde{P}\u,A^*_{\delta}\u \rangle|
 + |\langle A_{\delta}\u,A^*_{\delta}\tilde{P} \u \rangle| \\
 &+ |\langle E_{\delta}\u_p,\u_p \rangle| + |\langle C_{\delta}\u_p,\u_p \rangle| + |\langle K_{\delta}\u_p,\u_p \rangle|
 + |\langle \tilde{F}_{\delta}\u_p,\u_p \rangle|
 \nonumber \\
 &+ |\langle [A^2_{\delta},\tilde{P}]\u_p, \u_s \rangle| +|\langle [A^2_{\delta},\tilde{P}]\u_s, \u_p \rangle|
 + |\langle [A_{\delta}^2,\tilde{P}]\u_s, \u_s \rangle|
 \nonumber
\end{align}
We now justify why each term in the RHS of the above inequality remains uniformly bounded. First observe that for an operator $G \in \Psi_e^{2m-1,2l'-2}$ one has
\begin{equation}\label{eq:cosio mapping to the dual space}
G: H_e^{m-\teps,l} \to H_e^{-m+1-\teps,-l+(f-1)} \subset H_e^{-m+\teps,-l+(f-1)} = (H_e^{m-\teps,l})^*.
\end{equation}
The operators $F_{\delta}$, $\tilde{F}_{\delta}$, $C_{\delta}$
are uniformly bounded in $\Psi_e^{2m-1,2l'-2}$ i.e. of lower order than the main term. Since $\u_p$ is microlocally in $H^{m-\teps,l}_e$ on their microsupports (which are contained in $\WF'_{e,L^{\infty}}(\mathbf{A})$), then (\ref{eq:cosio mapping to the dual space}) implies $|\langle \tilde{F}_{\delta}\u_p,\u_p \rangle|$, $|\langle F_{\delta}\u_p,\u_p \rangle|$, $|\langle C_{\delta}\u_p,\u_p \rangle|$ are valid dual pairings and remain uniformly bounded even as $\d \to 0$. \par
 Let us turn to the term with $E_{\delta}$, where $E_{\d}$ differs depending on whether we are proving $(1)$ or $(2)$ of this lemma. If the hypothesis of $(1)$ in the lemma are satisfied, then as stated in Lemma \ref{eq:cosio Positive commutator}, we may arrange that $\WF'_{e,L^{\infty}}(E_{\d})$ is uniformly bounded away from $\p M$, and that $\WF'_{e,L^{\infty}}(E_{\d}) \cap \WF^{m,l}(\u_p) = \emptyset$ using (\ref{lemma:coiso baby1}) (recall that $\u_p$ satisfies the same hypothesis as $\u$). Thus, the term with $E_{\d}$ remains uniformly bounded in this case. If instead we are proving $(2)$ of the lemma, then $(3)$ of Lemma \ref{eq:cosio Positive commutator} shows $\u \in H_e^{m,l}$ microlocally on $\WF'_{e,L^{\infty}}(E_{\d})$ (so $\u_p$ will have the same property as explained before), where in this case
 $$ \WF'_{e,L^{\infty}}(E_{\d}) \subset U_1 \setminus \overline{\mcF}^p_{O,\alpha} \text{ uniformly }.$$
 Thus, the term with $E_{\d}$ remains uniformly bounded just like the `incoming' radial point case we just described. \par

 For the $K_{\delta}$ term, we have $q_p$ elliptic on $\text{WF}'_{e,L^{\infty}}(K_{\delta})$. Using the edge pseudodifferential calculus, since $\P$ and $\begin{bmatrix}Q_p & 0 \\ 0 & Q_s \end{bmatrix}$ have the same principal symbol,
 $$ \P = \begin{bmatrix}Q_p & 0 \\ 0 & Q_s \end{bmatrix} + P_1$$
 for some $P_1 \in \Psi_e^{1,-2}(M;TX)$. Thus, inside $O$ where bundles are trivialized,
 $\P \u_p = Q_p \u_p + P_1 \u_p$ which implies $Q_p \u_p \in H_e^{m-1-\teps,l-2}$ microlocally on $\WF'_{e,L^{\infty}}(\mathbf{A})$.
 Hence, by microlocal elliptic regularity of $Q_p$, we then have $u_p \in H_e^{m+1-\teps,l}$ microlocally on
$\text{WF}'_{e,L^{\infty}}(K_{\delta})$ so $K_{\d}\u_p \in H_e^{-m+1-\teps,l-2l'}\subset H_e^{-m+\teps,l-2l'}$ even as $\delta \to 0$. \par
Next, all the terms in (\ref{eq:coiso estimating Bu_p}) involving $\u_s$ are uniformly bounded by Proposition \ref{prop:coiso uniform boundedness of u_s terms}. \par
Lastly, we must justify the uniform boundedness of the terms containing $\P \u$, which requires a closer analysis of the principal symbol of $A_{\d}$. For this part, $\alpha$ being an incoming radial point versus outgoing is irrelevant so we suppress this distinction. As this type of argument is standard for positive commutator estimates, we merely refer the interested reader to \cite[page 19 proof of Theorem 1.5]{HVRadial} and \cite[proof of Lemma 9.6.1]{VKThesis} for the details since the only relevant feature is having a scalar principal symbol for $A_{\delta}$.
\end{lproof}

This lemma will be all we need in order to prove Theorem \ref{thm:propInOutEdge}. \par

 \begin{lproof}[(Proof of Theorem \ref{thm:propInOutEdge}).] The following proof is taken directly from \cite{MVWCorners} with only minor notational changes, and we provide extra details to enhance the transparency of the proof. We will only prove $(1)$ of the theorem, as the proof of $(2)$ would only require some trivial sign changes. By assumption, $m>l + f/2$ implies there exists $\teps \in (0,1/2]$ such that
 $$ m > l + f/2 + \teps.$$
First, observe that if we have $u \in H_e^{-N,l}$ and $-N > l + f/2$, then $-N + \teps > l+f/2 + \teps$, and so by applying the first part of Lemma \ref{lemma:intoTheEdge}
iteratively (with $m$ replaced first by $-N+\teps$), improving by $\teps$ in edge-Sobolev order at each step, we obtain $\alpha \notin \WF_e^{m,l}(u)$ (however at the last step of the iteration we might only need to improve edge-Sobolev order by an amount less than $\teps$).

To consider the other case, suppose $ l \geq -N - f/2$. Now define
$$l_0 = \sup \{ r : \alpha \notin \WF_e^{m,r}(u) \}.$$
First notice that the set over which the supremum is taken is non-empty since we can always find an $r_0 >> 0$ such that $ l+f/2-r_0 < -N,$ so that an analogous iterative procedure as in the first case shows $ \alpha \notin \WF_e^{m,l-r_0}(u)$. Next observe that if we can show that $l_0 = l$, then the theorem will be proved. The details are in \cite{MVWEdges}.
\end{lproof}

In the next section, we will improve this last theorem by showing coisotropic regularity into and out of the edge, that is, regularity at $\alpha$ under application of elements in $\mcA^k_p$ to $u$ under certain assumptions.

\section{ Propagation of coisotropic regularity }\label{sec: cosio regularity}

In this section, we get the coisotropic improvement by building up from the theorem in the previous section. The main result is at the end of this section, which is Theorem \ref{thm:CoisoVersion1}.

The first theorem we will prove is an analogue of Theorem \ref{thm:propInOutEdge} but with an improvement incorporating the propagation of regularity of $u$ under application of elements in $\mcA^k$. Since we will first prove propagation along $p$-bicharacteristics, we will often suppress the $p/s$ distinction in the modules by assuming that $\mcA$ and $\mcM$ will refer to the $p$-versions i.e. $\mcM_p$ and $\mcA_p$ introduced in Section \ref{sec:Coiso Regularity}; we note however that all results here will hold for the $s$-bicharacteristics as well, and we will provide more details of this at the end of the section.

\begin{theorem} \label{theorem:Coisopropagation in/out of edge}
 Let $u \in H_e^{-N,l}$ be a distribution.
 \begin{enumerate}
 \item  Let $m> l + f/2+1/2$. Given $\alpha \in \mcH^p_I$, if $\mathcal{A}_p^ku \subset H^m$ microlocally in
 $\mathcal{F}^p_{I,\alpha} \setminus \partial M$ and $\mcA_p^kPu \subset H_e^{m-1,l-2}$ microlocally at $\alpha$, then $ \mathcal{A}_p^ku \subset H_e^{m,l'}$ microlocally at $\alpha,\ \forall l' <l.$

 \item  Let $m < l+f/2.$ Given $\alpha \in \mathcal{H}^p_O$, if there exists a neighborhood
 $U$ of $\alpha$ in $\left.\right.^eS^*|_{\partial M}M$ such that
 $\WF_e^{m,l}(A_{\gamma}u) \cap U \subset \partial \mathcal{F}^p_O$ for all $A_{\gamma} \in \mathcal{A}_p^k$ and $\mcA^k_p Pu \in H_e^{m-1,l-2}$ microlocally at $\alpha$,
then $ \mathcal{A}_p^ku \subset H_e^{m,l'}$ microlocally at $\alpha,\ \forall l' <l.$
 \end{enumerate}
 \end{theorem}
 \begin{rem}
 One should notice that there is a loss of order $1/2$ in the case of $(1)$ of the theorem compared to the theorem in the previous section. This is because regularization is not free as we saw in the proof of the previous theorem, yet one may only improve coisotropic regularity by positive integer powers $\mcA$. Hence, it is not enough to regularize by just $\teps$ as in the proof of Lemma \ref{lemma:coiso baby1}. If one could microlocalize in such a way as to make sense of $\mcA^k$ for non-integer $k$, then indeed we would not have this loss of $1/2$ edge-derivatives.
 \end{rem}

 \begin{rem}
 This is a remark taken from \cite[Remark 11.2]{MVWEdges}. For any point $\alpha \in \e T^*_{\p M}M \setminus \{ \zeta = 0\}$ there is an element of $\mcA^k_p$ elliptic there, hence $(1)$, with $k = \infty$, shows that solutions with (infinite order) coisotropic regularity have no wavefront set in $ \e T^*_{\p M}M \setminus \{ \zeta = 0\}$. Indeed this result holds microlocally in the edge cotangent bundle. Note that the set $\{x =0, \zeta = 0\}$ is just the set of radial points for the Hamilton vector field $H_{q_p}$.
 \end{rem}

We will again work with $\u$ and $\P \u$ as introduced in the last section. As before, the ellipticity of $S$ and $S^*$ imply that $\u$ and $\P \u$ satisfy the same assumptions as $u$ and $Pu$ of the above theorem. Hence, it will suffice to prove coisotropic regularity at $\alpha$ for $\u$. The case of $k=0$ is exactly Theorem \ref{thm:propInOutEdge} already proven in the last section.
To get the coisotropic improvement, we give the following heuristic to show what we plan to do.

First, since the theorem is local in nature, let us fix some small neighborhoods to make our constructions here more explicit. Let $O$ be a neighborhood of the projection of $\alpha$ to $M$ inside a coordinate patch where all bundles are trivialized. We will assume all operators constructed have Schwartz kernels compactly supported in $O \times O$. Next, let
 \begin{equation}\label{eq:coiso2 K and U local sets}
 U_1 \subset \e S^*M
  \end{equation}
  be a precompact neighborhood of $\alpha$ away from the glancing rays such that
  $$
  \WF_e^{m,l}(\mcA^{k-1} \u) \cap U_1 = \emptyset \text{ and }
   \WF_e^{m-1,l}(\mcA^k \P \u) \cap U_1 = \emptyset.$$
   To avoid cluttered notation we will write $\mcM$ and $\mcA$ when we really mean $\mcM(U_1)$ and $\mcA(U_1)$ as in Definition \ref{def:coisotropy} often without explicitly clarifying. \par

Let $A_{\g} \in \mcA^k$ be a generator with multiindex $\g$ as introduced in Section \ref{sec:Coiso Regularity}, and
\begin{align}\label{eq:coiso2 introducing the regularizer}
A = A_0 \otimes \Id \in \Psi_e^{0,0}(U_1), \ A_{\gamma,m',l'} = W_{m',l'}A_{\gamma}, \text{ and }A_{\gamma,m',l',\d} = \Lambda_{\d}A_{\gamma,m',l'}
\end{align}
 for $\{ \Lambda_{\d}\}_{\d \in [0,1]}$ a scalar, uniformly bounded family of operators in $\Psi_e^{0,0}$ which will serve as a regularizer. Assuming $m'$ is chosen such that all the following quantities are bounded and that the integration by parts is valid, we have
\begin{align}\label{eq:coiso2 expanding <AA_mPu,AA_m u>}
\langle AA_{\gamma,m',l',\d}&\tilde{P}\u, A^*A_{\gamma,m',l',\d}\u \rangle - \langle AA_{\gamma,m',l',\d}\u, A^*A_{\gamma,m',l',\d} \tilde{P}\u \rangle
\nonumber \\
&= \langle [A^*_{\gamma,m',l',\d}A^2A_{\gamma,m',l',\d},\tilde{P}]\u,\u \rangle
\nonumber \\
&= \langle [A^*_{\gamma,m',l',\d}A^2A_{\gamma,m',l',\d},Q_p]\u_p,\u_p \rangle
+ \langle [A^*_{\gamma,m',l',\d}A^2A_{\gamma,m',l',\d},F]\u_p,\u_p \rangle
\nonumber \\
&\qquad + \langle [A^*_{\gamma,m',l',\d}A^2A_{\gamma,m',l'},\tilde{P}]\u_p,\u_s \rangle
 +\langle [A^*_{\gamma,m',l',\d}A^2A_{\gamma,m',l'},\tilde{P}]\u_s,\u \rangle,
\end{align}
where in the last equality we used (\ref{eq:coiso writing our [A,P]}) and (\ref{eq:coiso breaking up <[A,P]u_p,u_p>}) for some $F \in \Psi_e^{1,-2}(M;TX)$.

As a first step, for $A_{\gamma} \in \mathcal{A}^{k}$, to first prove coisotropic regularity of order $k$ with respect to $H^{m-1/2,l}_e$ on $\text{ell}(A)$, we must bound the quantities
$$\langle AA_{\gamma,m',l',\d}\tilde{P}\u_p, AA_{\gamma,m',l',\d}\u_s \rangle \text{ and } \langle AA_{\gamma,m',l',\d}\tilde{P}\u_s, AA_{\gamma,m',l',\d}\u \rangle$$
 appearing in (\ref{eq:coiso2 expanding <AA_mPu,AA_m u>}). To do this, we will obtain elliptic estimates for $A_{\gamma}\u_s$ using Proposition \ref{prop:semiElliptic} to directly bound these quantities. Afterwards, will do a careful commutator computation to bound the term with $Q_p$.

 In order to do commutator estimates, it will be convenient first to obtain a model form for commutators involving $A_{\gamma}$ with the following lemma:

\begin{lemma}\label{lemma:coiso2 commutator with module elements}
Let $G \in \Psi_e^{r',s'}$ and $A_{\gamma} \in \mathcal{A}^k$. Then
$$[G,A_{\gamma}] \in \Psi_e^{r',s'}\mathcal{A}^{k-1} :=  \{\sum_{j=1}^N Q_j \tilde{A}_j: Q_j \in \Psi_e^{r',s'}, \tilde{A}_j \in \mathcal{A}^{k-1}, N \in \mathbb{N} \} $$
\end{lemma}
(Note: No special properties of our particular module $\mathcal{A}$ are used here, so any such module suffices)
\proof This follows by induction and a tedious, explicit computation of the commutator. The details may be found in \cite[proof of Lemma 10.1.4]{VKThesis}.
$\Box$ \par
This lemma actually provides us with a very useful corollary.
\begin{cor}\label{cor:coiso2 I^kH^m mapping property}
Let $A \in\Psi^{m',l'}_e\mcA^{k'}(U_1;TX)$ for some $U_1 \subset \e S^*M$ a precompact open set. Then $A$ has the following mapping property for $k' \leq k$.
$$ A: \IH{k}{m}{l}{U_1} \to \IH{k-k'}{m-m'}{l+l'}{U_1}$$
\end{cor}
\proof
Let $u \in \IH{k}{m}{l}{U_1}$ and take any $A_{\gamma} \in \mcA^{k-k'}(U_1)$. Then
$$ A_{\gamma} Au = A A_{\gamma}u + [A_{\gamma},A]u.$$
Observe that $AA_{\gamma} \in \Psi_e^{m',l'}\mcA^{k}(U_1;TX)$ so that $AA_{\gamma}u \in H^{m-m',l+l'}_e$. Also, by the previous lemma $[A_{\gamma},A] \in \Psi_e^{m',l'}\mcA^{k-1}(U_1;TX)$ and so $[A_{\gamma},A]u \in H^{m-m',l+l'}_e$ as well. This completes the proof. $\Box$

With the aid of the above lemma, we may finally obtain elliptic regularity of $\mathcal{A}^k\u_s$ under certain regularity assumptions for $\mcA^k \u$ and
$\mcA^k \P \u$. This will be essential for bounding the terms in (\ref{eq:coiso2 expanding <AA_mPu,AA_m u>}) containing $\u_s$.
\begin{lemma}Let $\alpha \in \mcH^p$. Then
\begin{align*}
\alpha &\notin \WF_e^{m,l}(\mcA^{k}\u) \text{ and }
\alpha \notin \WF_e^{m-1,l-2}(\mcA^{k} \P \u)
\\
&\Rightarrow
\alpha \notin \WF_e^{m-1,l}(\mcA^k \P \u_s) \text{ and }\alpha \notin \WF_e^{m+1,l}(\mcA^k  \u_s).
\end{align*}
(The point here is that $\u_s$ is microlocally one edge-derivative smoother that $\u$. Also, $\mcA$ here really means $\mcA_p(U)$ where $U \subset \e S^*M$ is a neighborhood of $\alpha$ away from glancing rays. For details, see Section \ref{sec:Coiso Regularity}.)
\end{lemma}

\proof First let $O$ be a neighborhood of $\pi_M(\alpha)$ contained in a local chart where all bundles are trivialized, and we assume that all operators constructed here have Schwartz kernels supported on $O \times O$. Let $A_{\gamma} \in \mcA^k$ be arbitrary.
Then
$$ A_{\gamma}\P \u_s = A_{\gamma} \Pi_s \P \u + A_{\gamma}[\Pi_s,\P] \u.$$
Since $A_{\gamma}\Pi_s \in \Psi^{0,0}_e\mcA^k$ by Lemma \ref{lemma:coiso2 commutator with module elements}, then $A_{\gamma} \Pi_s \P \u \in H_e^{m-1,l-2}$ microlocally at $\alpha$ by the assumption on $\P \u$. Likewise, we have shown in Section \ref{sec:partial elliptic regularity} that $[\Pi_s, \P] \in \Psi_e^{1,-2}(M;TX)$ so $A_{\gamma}[\Pi_s, \P] \in \Psi^{1,-2}\mcA^k$ by Lemma \ref{lemma:coiso2 commutator with module elements}; this implies $A_{\gamma}[\Pi_s,\P] \u \in H^{m-1,l-2}_e$ microlocally at $\alpha$ by the assumption on $\u$. This proves the first part of the lemma, that $\mcA^k\P \u_s \subset H^{m-1,l-2}_e$ microlocally at $\alpha$.

Next let $\P^- \in \Psi_e^{-2,2}(M;TX)$ be the parametrix for $P^-$ as constructed in Proposition \ref{prop:semiElliptic}, where we also showed
$$ \u_s = \P^-\P \u_s + R'' \u_s $$
where $R'' \in \Psi_e^{0,0}(M)$ and $\alpha \notin \WF_e'(R'')$. Thus,
$$ A_{\gamma}\u_s = A_{\gamma}\P^-\P \u_s + A_{\gamma}R''\u_s.$$
Then $A_{\gamma}\P^- \in \Psi_e^{-2,2}\mcA^k$ by Lemma \ref{lemma:coiso2 commutator with module elements} so $A_{\gamma}\P^-\P \u_s \in H^{m+1,l}_e$ microlocally at $\alpha$. Proceeding, since $\alpha \notin \WF'_e(R'')$ then by microlocality of $R''$ then $A_{\gamma}R''\u_s \in H^{m+1,l}_e$ microlocally at $\alpha$ (in fact $\alpha \notin \WF^{\infty,l}_e(A_{\gamma}R''\u_s)$).
This shows, $A_{\gamma}\u_s \in H_e^{m+1,l}$ microlocally at $\alpha$, which concludes the proof of the lemma.
$\Box$

The proof of the above lemma actually shows something a little bit stronger where we can replace $\alpha$ by a small neighborhood of $\alpha$:
\begin{lemma}\label{lemma:coiso2 coiso elliptic regularity}
Let $\alpha \in \mcH^p$ and $U \subset \e S^*M$ be a precompact neighborhood of $\alpha$. Then
\begin{align*}
&U \cap \WF_e^{m,l}(\mcA^{k}\u) = \emptyset \text{ and }
U \cap \WF_e^{m-1,l-2}(\mcA^k \P \u) = \emptyset
\\
&\Rightarrow
U \cap \WF_e^{m-1,l-2}(\mcA^k  \P \u_s) = \emptyset \text{ and } U \cap \WF_e^{m+1,l}(\mcA^k  \u_s) = \emptyset.
\end{align*}
\end{lemma}

We also need to understand adjoints in $\mathcal{A}$. Since the operators in $\mcA^k$ and their adjoints are principally scalar with real principal symbols, if $A_{\gamma} \in \mathcal{A}^k$ then one has \emph{a priori} $A_{\gamma}^* - A_{\gamma} \in \Psi_e^{k-1}$, but
in fact, we can do better:
\begin{lemma}\label{lemma:coiso2 adjoints in A^k}
We have
$$ A - A^* \in \Psi_e^{0,0}\mathcal{A}^{k-1} \text{ for any } A \in \mathcal{A}^k. $$
\end{lemma}

\begin{rem}
We must put $\Psi_e^{0,0}\mathcal{A}^{k-1}$ rather than just $\mcA^{k-1}$ in the above lemma since we are allowing $A -A^*$ to not be principally scalar.
\end{rem}

\proof We use induction. The case $k=1$ is trivial as mentioned before the lemma since elements of $\mcA$ have real, scalar principal symbols. Thus, suppose the result holds for $k-1$.
It suffices to prove this for the generators of $\mathcal{A}^k$, and so we let $A_{\gamma} \in \mathcal{A}^k$ be a generator.
Then $A_{\gamma} = A_{\gamma'}A_j$ for some $A_j$ a generator in $\mathcal{M}$ and $A_{\gamma'} \in \mathcal{A}^{k-1}.$ By the induction hypothesis applied to $A_{\gamma'}$, one has
$A_{\gamma}^* = A_j^*A_{\gamma'}^*=A_j^*(A_{\gamma'} + E_{k-2})$ for some $E_{k-2} \in \Psi^{0}_e\mathcal{A}^{k-2}$. Notice $A_j^* - A_j \in \Psi_e^0$ and
$[A_j,A_{\gamma'}] \in \mathcal{A}^{k-1}$, so we have
\begin{align*}
A_{\gamma}^*
&= A_jA_{\gamma'} + A_j E_{k-2} + (A_j^*-A_j)A_{\gamma'} + (A_j^*-A_j)E_{k-2}
\nonumber \\
&=  A_{\gamma} +[A_j,A_{\gamma'}] + A_j E_{k-2} + (A_j^*-A_j)A_{\gamma'} + (A_j^*-A_j)E_{k-2}.
\end{align*}
By Lemma \ref{lemma:coiso2 commutator with module elements}, $A_j E_{k-2} \in \Psi_e^{0,0}\mcA^{k-1}$, and similarly for all the other terms on the RHS of the above equation besides for $A_{\gamma}$. This gives the desired conclusion. $\Box$

We now turn to the main commutator proof for the term in (\ref{eq:coiso2 expanding <AA_mPu,AA_m u>}) involving $Q_p$. First, we have another crucial lemma will allow us to only worry about only those generators $A_{\gamma}$ which do not contain $A_N$. This
is no longer trivial, as $\tilde{P}\u$ having a certain amount of regularity does not imply automatically that $Q_p\u_p$ has the same amount of coisotropic regularity, which is what's needed
for proving coisotropic regularity of $\u_p$. Nevertheless, we've arranged everything so that the following lemma, similar to \cite[Lemma A.4]{MVWEdges} still remains true.

\begin{lemma}\label{lemma:coiso2 making A_N irrelevant}
Suppose $\u$ is coisotropic order $k-1$ on $O$ relative to $H_e^{m,l}$ and $\P \u \in \IH{k}{m-1}{l}{O}$. Then for $O' \subset O$, $u$ is coisotropic of order $k$ on $O'$
relative to $H_e^{m,l}$ if for each multiindex $\gamma$, with $|\gamma|=k, \gamma_{N}=0$, there exists $B_{\gamma} \in \Psi_e^{0,0}$, elliptic
on $O'$ such that $B_{\gamma}A_{\gamma,m,l}\u \in L^2$.
\end{lemma}

\proof
First, inside a local chart where bundles are trivialized, we represent $\tilde{P}$ as
\[  \tilde{P}=
\left[ \begin{array}{cc}
Q_p & 0   \\
0 & Q_s
\end{array} \right]
+ \left[ \begin{array}{cc}
F_1 & F_2   \\
* & *
\end{array} \right]
 \]
for $F_1,F_2 \in \Psi_e^{1,-2}$ and the $*$ also refer to operators in this space. By a partition of unity, we may suppose WLOG that $O$ is contained inside a local chart where bundles are trivialized. Then $\tilde{P}\u \in$ $ \IH{k}{m-1}{l}{O}$ implies
\begin{align}
 &Q_p\u_p \in \IH{k}{m-1}{l-2}{O} -F_1\u_p - F_2\u_s \\
 &\Rightarrow \mathcal{A}^{k-1}Q_p\u_p \subset\IH{1}{m-1}{l-2}{O} -\mathcal{A}^{k-1}F_1\u_p - \mathcal{A}^{k-1}F_2\u_s.
\end{align}
 One then uses Lemma \ref{lemma:coiso2 coiso elliptic regularity} to estimate $\u_s$, and the rest of the proof proceeds analogously as in \cite[proof of Lemma A.3]{MVWEdges}, with full details in our setting in \cite[Lemma 10.1.11]{VKThesis}.
$\Box$

To proceed with the commutator proof, we recall Lemma A.4 in \cite{MVWEdges} appropriately adapted to our setting where
we view our operators as principally scalar operators acting on a vector bundle. Due to the previous lemma, $\gamma,\beta$ will stand for \emph{reduced multiindices}, with $\gamma_N = 0$, $\beta_N =0$. Also, from now on, the choice of regularizer $\Lambda_{\d}$, and operator $A_0$ will play a crucial role. Indeed we let $A_0$ to be the quantization of $a'$ in (\ref{eq:coiso a' definition for in}), and
$\Lambda_{\d}$ constructed in (\ref{eq:coiso Regularizer construction}), to be the regularizer whose principal symbol is
\begin{equation}\label{eq:coiso2 regularizer principal symbol}
\sigma_e(\Lambda_{\d}) = \sqrt{\phi_{\d}(|\tau|^2)} = (1+\d|\tau|^2)^{-\teps/2}
\text{ on } \supp(a').
\end{equation}
Thus, $\{\Lambda\}_{\d \in (0,1]} \subset \Psi_e^{-\teps,0}(M)$, is uniformly bounded in $\Psi_e^{0,0}(M)$ down to $\d \to 0$, and $\Lambda_{\d} \to \Id$ in
$\Psi_e^{\epsilon,0}(M)$ for any $\epsilon>0$.

\begin{lemma}\label{lemma:coiso2 commuting Q_p through}
Suppose $A_0 \in \Psi_e^{0,0}$ is as described above. Then
\begin{align}\label{eq:coiso2 commuting Q_p through}
\sum_{|\gamma| = k} &i[A^*_{\gamma,m',l',\d}A^2_0A_{\gamma,m',l',\d},Q_p]
\\
&=\sum_{|\gamma|,|\beta| = k} A^*_{\gamma,m'+1/2,l'-1,\d}A^*_0 C'_{\gamma \beta,\d}A_0A_{\beta,m'+1/2,l'-1,\d}
\nonumber \\
&+ \sum_{|\gamma|=k}( A^*_{\gamma,m'+1/2,l'-1,\d}A_0E_{\gamma,m'+1/2,l-1,\d} + E^*_{\gamma,m'+1/2,l-1,\d}A_0A_{\beta,m'+1/2,l'-1,\d})
\nonumber\\
&+ \sum_{|\gamma| = k} A^*_{\gamma,m',l',\d}i[A_0^2,Q_p]A_{\gamma,m',l',\d}
\nonumber
\end{align}
where
$$E_{\gamma,m',l',\d} = W_{m',l'}E_{\gamma,\d}, \quad E_{\gamma,\d} \in \Psi_e^{0,0}\mathcal{A}^{k-1} + \Psi_e^{0,0}\mathcal{A}^{k-1}A_N,$$
$$
\WF'_{e,L^{\infty}}(E_{\gamma,\d}) \subset \WF'(A_0)
\text{ uniformly and for all $\gamma,\beta$},$$
\begin{equation}\label{eq:coiso2 C' is non zero}
C'_{\gamma\beta} \in \Psi_e^{0,0} \text{ uniformly bounded, } \ \sigma(C'_{\gamma\beta})|_{\partial \dot{\mathcal{F}}} = -(m'+l'-\teps \d(1+\d |\tau|^2)^{-1})\hxi c_p^2\delta_{\gamma\beta}
\end{equation}
\end{lemma}
\begin{rem}This is an important remark taken from \cite[Remark A.5]{MVWEdges}. The first term on the right hand side of (\ref{eq:coiso2 commuting Q_p through}) is the principal term in terms of $\mcA^p$ order; both $A_{\gamma, m,l}$ and $\mcA_{\beta,m,l}$ have $\mcA$ order $k$. Moreover, (\ref{eq:coiso2 C' is non zero}) states that it has non-zero principal symbol near $\p \dmcF^p$ depending on the sign of $m'+l'$ and $\hat{\xi}$. The terms involving $E_{\gamma,m',l'}$ have $\mcA$ order $k-1$, or include a factor of $A_N$, so they can be treated as error terms. On the other hand $i[A_0^2,Q_p]$ must be arranged to be positive, which will come from Lemma \ref{eq:cosio Positive commutator} as we show below.
\end{rem}

\proof The proof proceeds exactly as in \cite[Lemma A.4]{MVWEdges} with full details in \cite[Lemma 10.1.12]{VKThesis}.
$\Box$

We now have all the pieces to prove Theorem \ref{theorem:Coisopropagation in/out of edge}.

{\it(Proof of Theorem \ref{theorem:Coisopropagation in/out of edge})} We do the proof by induction. The case $k=0$ was precisely Theorem \ref{thm:propInOutEdge}. We may prove both parts of the theorem at once and point out the relevant differences. Let us assume then the theorem holds for $k-1$, i.e.
$$ \alpha \notin \WF_e^{m,\tilde{l}}(\mcA^{k-1}u).$$
For notational convenience, we will suppose $\tilde{l} =l$ since the $\tilde{l}<l$ condition only came from the interpolation argument in Theorem \ref{thm:propInOutEdge}, but does not affect any of the arguments here.
It will be convenient to denote
$$ \bar{\xi} = \xi(\alpha).$$
Then by the closedness of $\WF_e$, there exists a neighborhood $O_k$ of $\alpha$ over which $\xi$ has a fixed sign which is that of $\bar{\xi}$, such that
\begin{enumerate}
\item[$\bullet$]$u \in \IH{k-1}{m}{l}{O_k}$
\item[$\bullet$] $Pu \in \IH{k}{m-1}{l-2}{O_k}$
\end{enumerate}
and the projection of $O_k$ to $M$ is inside a local coordinate patch where all bundles are trivialized, and all operators constructed have Schwartz kernels compactly supported in $\pi_M(O_k) \times \pi_M(O_k)$.
The ``first step'' is to show $\u$ has coisotropic regularity of order $k$ relative to $H^{m-1/2,l}_e$ microlocally over a compact subset of $O_k$, and then do an analogous ``second step'' to show $\u$ has coisotropic regularity of order $k$ relative to $H^{m,l}_e$ at $\alpha$. Since the proof of the ``first step'' is almost identical and even easier than the ``second step'', we only prove the ``second step'' and afterwards comment on why the ``first step'' is easier. Thus, let us instead assume
\begin{enumerate}
\item[$\bullet$]$\u \in \IH{k}{m-1/2}{l}{O_k}$
\item[$\bullet$] $\P \u \in \IH{k}{m-1}{l-2}{O_k}$
\end{enumerate}
so by Lemma \ref{lemma:coiso2 coiso elliptic regularity} (applying it twice first replacing $m$ by $m-1/2$ and then replacing $k$ by $k-1$) one also has
\begin{align}\label{eq:coiso2 u_s regular in baby proof}
 \u_s &\in \IH{k}{m+1/2}{l}{O_k}\cap \IH{k-1}{m+1}{l}{O_k}\text{ and }
 \\
 \P \u_s &= \Pi_s \P \u + [\P, \Pi_s] \u \in \IH{k}{m-3/2}{l-2}{O_k}
 \cap \IH{k-1}{m-1}{l}{O_k}
 \nonumber
\end{align}
by Corollary \ref{cor:coiso2 I^kH^m mapping property} and since $[\P,\Pi_s]$ is in $\Psi_e^{1,-2}$.
We need only show that $\u$ has coisotropic regularity of order $k$ relative to $H^{m,l}_e$ over some neighborhood of $\alpha$.

Now let $A_{\gamma}$ and $A_{\gamma,m',l',\d}$ be as described in (\ref{eq:coiso2 introducing the regularizer}) with $U_1$ there replaced by $O_k$, and
 $$ m' = m-1/2 \qquad \qquad l' = -l- (f-1)/2 .$$
 Note that for the regularizing term of $A_{\gamma,m',l',\d}$ explained in (\ref{eq:coiso2 regularizer principal symbol}), we take
 $$ \teps = \frac{1}{2}.$$

 Also let $A \in \Psi_e^{0,0}(O_k;\mathbb{C}^n)$ be a scalar operator with principal symbol $a'$ as described in
(\ref{eq:coiso a' definition for in}). For clarity, we also write
$$ A = A_0 \otimes \Id.$$
As, shown there, the microsupport of $A$ may be made arbitrarily close to $\alpha$ so that one indeed has
$$ \WF'_{e,L^{\infty}}(AA_{\gamma,m',l',\d}) \subset O_k.$$

Thus, as shown in (\ref{eq:coiso2 expanding <AA_mPu,AA_m u>}), for $\d >0$ (where the integration by parts is justified by the proof in Lemma \ref{lemma:coiso baby1})
\begin{align}\label{eq:coiso2 mainExpansion}
  \langle A&A_{\gamma,m',l',\d}\u, A^*A_{\gamma,m',l',\d}\P \u \rangle
- \langle AA_{\gamma,m',l',\d}\tilde{P}\u, A^*A_{\gamma,m',l',\d}\u \rangle
\\
&=\langle AA_{\gamma,m',l',\d}\u, A^*A_{\gamma,m',l',\d}\P \u_s \rangle
- \langle AA_{\gamma,m',l',\d}\tilde{P}\u_s, A^*A_{\gamma,m',l',\d}\u \rangle
\nonumber \\
&\qquad+ \langle AA_{\gamma,m',l',\d}\u_s, A^*A_{\gamma,m',l',\d}\P \u_p \rangle
- \langle AA_{\gamma,m',l',\d}\tilde{P}\u_p, A^*A_{\gamma,m',l',\d}\u_s \rangle
\nonumber \\
&\qquad+ \langle [A^*_{\gamma,m',l',\d}A_0^2A_{\gamma,m',l',\d},Q_p]\u_p,\u_p \rangle
+ \langle [A^*_{\gamma,m',l',\d}A^2A_{\gamma,m',l',\d},F]\u_p,\u_p \rangle,
\nonumber
\end{align}
with $F \in \Psi_e^{1,-2}(M;TX)$.

 So with the notation in Lemma \ref{lemma:coiso2 commuting Q_p through}, as done in \cite[Appendix]{MVWEdges}, we let
 $$C'_{\d} = \text{diag}((C'_{\gamma\beta,\d}), \dots, (C'_{\gamma\beta,\d}))\
 \text{($n$ copies) such that } |\gamma| = |\beta| = k,$$
 as a block diagonal matrix (with $n$ blocks) of
operators, or rather as an operator on a trivial vector bundle with fiber $\mathbb{R}^{n|M_k|}$ over
a neighborhood of $\mathcal{F}_{\alpha}^p$, where $|M_k|$ denotes the number of elements of the set $M_k$ of
multiindices $|\gamma|=k$. Also, let
$$c'_{\d} = \sigma(C'_{\d})|_{\partial \dot{\mathcal{F}_{\alpha}^p}} = -(m'+l'-\teps \d(1+\d |\tau|^2)^{-1})c_p^2\hat{\xi} \Id_{n|M_k|}.$$
Thus, $c'_{\d}$ is positive or negative definite with the sign of $-(m'+l'-\teps)\bar{\xi}$, so the same
is true microlocally near $\dot{ \mathcal{F}_{\alpha}^p }.$ Thus we have
$$
\sgn(\bar{\xi})\sigma(C'_0)|_{\partial \dot{ \mathcal{F}^p}}
\begin{cases}
 < 0 &
\mbox{ if } m'+l' -\teps >0  \\
 > 0 &
\mbox{ if } m'+l' < 0
\end{cases}
$$
The first case happens when assuming $(1)$ of the theorem, while the second case happens when $(2)$ is assumed.
 Then shrinking $O_k$ if necessary, we may find $B \in \Psi_e^{0,0}$, $G \in \Psi_e^{-1,0}$, with $\sigma(B)>0$ on $O_k$ such that
\begin{equation}\label{eq:coiso2 A_0 C' A_0 positive}
\sgn(\bar{\xi})A^*_0C_{\d}'A_0 = A^*_0(\mp B^*B + G)A_0
\end{equation}
with the $(-)$ in the case of $(1)$ of the lemma and the $(+)$ in the case of $(2)$ of the lemma.
Also, shrinking $O_k$ if necessary, we have from the proof of Lemma \ref{eq:cosio Positive commutator} (where we do not use the regularizer $\phi_{\d}$ there, and only consider the Hamilton derivative of $a'$ there) that
\begin{align}\label{eq:coiso2 [A_0^2,Q_p]}
&\sgn(\bar{\xi})i[A_0^2, Q_p] = \sum \mp\tilde{B}_j^*\tilde{B}_j + E + C + K + \tilde{F}\\
& \tilde{B}_j \in \Psi_e^{1/2,-1}(O_k),E,K \in \Psi_e^{1,-2}, C,\tilde{F} \in \Psi_e^{0,-2}(O_k),
\nonumber
\end{align}
where $\u_p$ is coisotropic of order $k$ with respect to $H_e^{m,l}$ in a neighborhood of $\WF'_e(E)$ as explained in the proof of Lemma \ref{lemma:coiso baby1}, and $\WF'_e(K) \cap \Sigma_p = \emptyset$. Note that the $(\mp)$ we have in (\ref{eq:coiso2 [A_0^2,Q_p]}) is different from the signs in that lemma since there, the sign of $\bar{\xi}$ was taken into account rather than the sign of $m'+l'$.

Now, let
$$\hat{A}_\d \u_p = (A_0A_{\gamma,m'+1/2,l'-1,\d}\u_p)_{|\gamma| = k },$$ regarded as a column vector
of length $n|M_k|$. Now, using Lemma \ref{lemma:coiso2 commuting Q_p through} and substituting (\ref{eq:coiso2 A_0 C' A_0 positive}),(\ref{eq:coiso2 [A_0^2,Q_p]}) into (\ref{eq:coiso2 commuting Q_p through}) we obtain
\begin{align*}
\sum_{|\gamma|= k} &\sgn(\hat{\xi})\langle i[A^*_{\gamma,m',l',\d}A^2_0A_{\gamma,m',l',\d}, Q_p]\u_p, \u_p \rangle
\\
&= \mp||B\hat{A}_{\d}\u_p||^2 \mp||\tilde{B}_j\hat{A}_{\d}\u_p||^2
 \\
&\qquad +\sum_{|\gamma|=k} \Big(  \langle E_{\gamma,m'+1/2,l'-1,\d}\u_p, A_0A_{\gamma,m'+1/2,l'-1,\d}\u_p \rangle
\\
&\qquad \qquad \qquad + \langle A_0A_{\gamma,m'+1/2,l'-1,\d}\u_p, E_{\gamma,m'+1/2,l'-1,\d}\u_p \rangle
\Big) \\
&\qquad + \sum_{|\gamma|= k} \langle (E+C+K+\tilde{F})A_{\gamma,m',l',\d}\u_p, A_{\gamma,m',l',\d}\u_p \rangle
+ \langle \hat{A}\u_p, G\hat{A}\u_p \rangle
\end{align*}
If we substitute the above equation into (\ref{eq:coiso2 mainExpansion}),
drop the terms involving $\tilde{B}_j$ and apply the Cauchy-Schwartz inequality to the terms with $E_{\gamma,m'+1/2,l'-1,\d}$, we have for any $\epsilon_1>0$
\begin{align}\label{eq:coiso2 main BAu_p estimate}
||B&\hat{A}_{\d}\u_p||^2
\\
&\leq \sum_{|\g|=k}\Big(||\langle AA_{\gamma,m',l',\d}\u, A^*A_{\gamma,m',l',\d}\P \u \rangle|
+ |\langle AA_{\gamma,m',l',\d}\tilde{P}\u, A^*A_{\gamma,m',l',\d}\u \rangle|
\nonumber \\
&\qquad +|\langle AA_{\gamma,m',l',\d}\u, A^*A_{\gamma,m',l',\d}\P \u_s \rangle|
+ |\langle AA_{\gamma,m',l',\d}\tilde{P}\u_s, A^*A_{\gamma,m',l',\d}\u| \rangle
\nonumber \\
&\qquad+ |\langle AA_{\gamma,m',l',\d}\u_s, A^*A_{\gamma,m',l',\d}\P \u_p| \rangle
 + |\langle AA_{\gamma,m',l',\d}\tilde{P}\u_p, A^*A_{\gamma,m',l',\d}\u_s \rangle|
\nonumber \\
&\qquad + |\langle [A_{\gamma,m',l',\d}A^2A_{\gamma,m',l',\d},F]\u_p,\u_p \rangle| \Big)
\nonumber \\
&\qquad + \epsilon_1||\hat{A}_{\d}\u_p||^2 +
\epsilon_1^{-1}\sum_{|\gamma|=k} ||E_{\gamma,m'+1/2,l'-1,\d}\u_p||^2
\nonumber \\
&\qquad + \sum_{|\gamma|= k} |\langle (E+C+K+\tilde{F})A_{\gamma,m',l',\d}\u_p, A_{\gamma,m',l',\d}\u_p \rangle|
+ |\langle \hat{A}_{\d}\u_p, G\hat{A}_{\d}\u_p \rangle|.
\nonumber
\end{align}
Before proceeding, let us first make a remark regarding the intuition of the proof.
\begin{rem}
The goal is to uniformly bound the rest of the terms on the right of the above inequality. The term containing $\P\u$ will remain bounded due to \emph{a priori} assumptions on $\P \u$, which put it in a better space microlocally at $\alpha$ then what would be dictated by the space $\u$ is in. The terms with $\u_s$ are also considered ``error'' terms since $\u_s$ satisfies better elliptic estimates. All the remaining terms are either of lower order so that they may be bounded by the inductive assumption, or are terms with $E$ or $K$, which are microsupported in regions where we assume $\u$ is better to begin with, so may be treated as ``error'' terms as well. Hence, in the ensuing proof we are merely justifying why all the terms may be treat as such error terms that are bounded due to the inductive step or the initial assumptions, as is standard in all positive commutator proofs.
\end{rem}
Choosing $\epsilon_1>0$ small enough, the $||\hat{A}_{\d}\u_p||^2$ term on the right can be absorbed in the left hand side as done in \cite[Proof of Proposition A.6]{MVWEdges}.

Let us justify that the terms with $E_{\gamma,m'+1/2,l'-1,\d},E,K,C,\tilde{F},G,\P \u$ remain uniformly bounded as well.
\begin{enumerate}
\item[$\bullet$]
First, we turn to the term with $E_{\gamma,m'+1/2,l'-1,\d}$. Writing $\u_p = \u - \u_s$ and using (\ref{eq:coiso2 u_s regular in baby proof}) gives
\begin{equation}\label{eq:coiso2 u_p is in coisotropic space}
\u_p \in \IH{k-1}{m}{l}{O_k}.
\end{equation}
Next, the proof of Lemma \ref{lemma:coiso2 making A_N irrelevant} then shows
\begin{align}\label{eq:coiso2 Q_pu_p is better}
 &A_N \u_p \in \IH{k-1}{m}{l}{O_k}
  \\
 &\Rightarrow \Psi_e^{m,l'-1}\mcA^{k-1}(O_k)A_N \u_p \subset H_e^{0,-(f+1)/2}
\nonumber
\end{align}

   Thus,
    \begin{align*}
   &E_{\g,m,l'-l,\d}\u_p \in H^{0,-(f+1)/2}_e \text{ uniformly as } \d \to 0
    \\
    &\Rightarrow
    ||E_{\g,m,l'-l,\d}\u_p|| \text{ remains uniformly bounded as } \d \to 0.
    \end{align*}
\item[$\bullet$]
For the terms involving $C, \tilde{F}$, $G$,$\tilde{F}' \in \Psi_e^{0,-2}(O_k')$, $G' \in \Psi_e^{-1,0}(O_k')$ are lower order than the principal term and so these remain uniformly bounded exactly as in the proof of Proposition \ref{prop:coiso uniform boundedness of u_s terms} when we bounded the term with $F$ there. Exact details are already in \cite[Proposition A.6]{MVWEdges} and \cite[Theorem 10.1.1]{VKThesis}.
\item[$\bullet$]
We now turn to the terms involving $K$ and $E$. Now $Q_p$ is elliptic on $\WF'_{e,L^{\infty}}(\Lambda_{\d}K)$ so by microlocal elliptic regularity, there exists a $Q^- \in \Psi_e^{-2,2}$ such that
$$ I = Q^- Q_p + R_1 \text{ for some }R_1 \in \Psi_e^{0,0}(O_k),
\WF'(R_1) \cap \WF'_{e,L^{\infty}}(\Lambda_{\d}K) = \emptyset.$$
So using Corollary \ref{cor:coiso2 I^kH^m mapping property}, (\ref{eq:coiso2 Q_pu_p is better}) to give us the order of $Q_p \u_p$, and this microlocal elliptic regularity of $Q_p$ implies $\u_p$ is coisotropic of order $k$ on $\WF'_{e,L^{\infty}}(\Lambda_{\d}K)$ with respect to $H^{m,l}_e(M)$.
Similarly, by the hypothesis of the theorem, we similarly have $\u_p$ is coisotropic of order $k$ microlocally
 in a neighborhood of $\WF'_{e,L^{\infty}}(\Lambda_{\d}E)$ with respect to $H_e^{m,l}$. By looking at the orders of the operators $E,K$ and the microlocal edge-Sobolev regularity of $\u_p$ over $\WF'_{e,L^{\infty}}(E+K)$ gives uniform boundedness of these terms exactly as in \cite[Proposition A.6]{MVWEdges}.

\item[$\bullet$]
 We turn to the term involving $\P \u$. Estimating this term is a standard argument in commutator proofs since $\P \u$ is of a better order than \emph{a priori} expected by the Sobolev order of $\u$. We refer the interested reader to \cite[Proof of Theorem 1.5]{HVRadial} and \cite[Proof of Theorem 10.1.1]{VKThesis}.

\item[$\bullet$] Proceeding, $F'$ may also be estimated as a lower order error term exactly as the term with $E_{\gamma, m'+1/2,l'-1,\d}$ as shown in the proof \cite[Theorem 10.1.1]{VKThesis}.

\item[$\bullet$] Finally, we look at the terms containing $\u_s$. Let us first consider
    \begin{align}\label{eq:coiso2 commuting P through}
    \langle [A^*_{\gamma,m',l',\d}A^2&A_{\gamma,m',l',\d},\P]\u_s,\u \rangle
    \end{align}
    Observe that
    $$ \ord(\u_s,\u)= (2k, 2m),$$
    while
    $$
    [A^*_{\gamma,m',l',\d}A^2A_{\gamma,m',l',\d},\P]
    \in \Psi_e^{2m,2l'-2}\mcA^{2k} + \Psi_e^{2m+1,2l'-2}\mcA^{2k-1}.
    $$
    Thus, writing $ [A^*_{\gamma,m',l',\d}A^2A_{\gamma,m',l',\d},\P] = C_1+ C_2$ with $C_1 \in \Psi_e^{2m,2l'-2}\mcA^{2k}$ and $C_2 \in \Psi_e^{2m+1,2l'-2}\mcA^{2k-1}$, the term $\l C_1 \u_s, \u \r$ remains uniformly bounded, while the term $\l C_2 \u_s, \u \r$
     may be dealt with as the $F'$ term above using an analogous estimate as for the $\P \u$ term.
     The other $`\u_s'$ term
     $$\l [A^*_{\g,m',l',\d}A^2A_{\g,m',l',\d},\P]\u_p, \u_s \r$$
      appearing on the right of (\ref{eq:coiso2 main BAu_p estimate}) is dealt with in a similar fashion since all that is relevant are the edge-Sobolev spaces $\u_s$ and $\u_p$ belong to.
\end{enumerate}

Thus, after absorbing all the terms mentioned into the left of (\ref{eq:coiso2 main BAu_p estimate}) and then letting $\d \to 0$ shows $\u_p$ has coisotropic regularity of order $k$ on $\text{ell}(A)$ with respect to $H_e^{m-1/2,l}$ due to Lemma \ref{lemma:coiso2 making A_N irrelevant}. Since $\u_s$ has the same property, we have shown $\u$ has coisotropic regularity of order $k$ on $\text{ell}(A)$ with respect to $H_e^{m-1/2,l}$.

To show the first part of the argument, one assumes $\u \in \IH{k-1}{m}{l}{O_k'}$ for any open $O_k'$ whose closure is compactly contained in $O_k$, and then shows $\u$ has coisotropic regularity of order $k$ with respect to $H^{m-1/2,l}_e$ microlocally over $\alpha$. This is easier than the proof above since we only need to improve coisotropic order, but not Sobolev order. The proof is exactly analogous to the proof of \cite[Proposition A.6]{MVWEdges}, and the full details in our setting are in \cite[Theorem 10.1.1]{VKThesis}.
$\Box$

Notice that all our proofs only relied on the diagonal form
of the principal symbol of $\tilde{P}$, so interchanging $p$ and $s$ in all the
above proofs, but instead looking at $\alpha \in \Sigma_s$, gives us the parallel theorem:
\begin{theorem}
 Let $u \in H_e^{-N,l}$ be a distribution.
 \begin{enumerate}
 \item  Let $m> l + f/2+1/2$. Given $\alpha \in \mcH^s_I$, if $\mathcal{A}_s^ku \subset H^m$ microlocally in
 $\mathcal{F}^s_{I,\alpha} \setminus \partial M$ and $\mcA_s^kPu \subset H_e^{m-1,l-2}$ microlocally at $\alpha$, then $ \mathcal{A}_s^ku \subset H_e^{m,l'}$ microlocally at $\alpha,\ \forall l' <l.$

 \item  Let $m < l+f/2.$ Given $\alpha \in \mathcal{H}^s_O$, if there exists a neighborhood
 $U$ of $\alpha$ in $\left.\right.^eS^*|_{\partial M}M$ such that
 $\WF_e^{m,l}(A_{\gamma}u) \cap U \subset \partial \mathcal{F}^s_O$ for all $A_{\gamma} \in \mathcal{A}_s^k$ and $\mcA^k_s Pu \in H_e^{m-1,l-2}$ microlocally at $\alpha$,
then $ \mathcal{A}_s^ku \subset H_e^{m,l'}$ microlocally at $\alpha,\ \forall l' <l.$
 \end{enumerate}
 \end{theorem}

We may finally combine all our results to prove the main propagation of coisotropic regularity theorem with an important corollary, following closely \cite[Section 8]{MVWCorners}.
\begin{theorem}
\label{thm:CoisoVersion1}
Let $\beta \in \dot{\mcH}^p$, $k \in \mathbb{N}$ and $\d>0$. Suppose that
$u \in H_e^{-N,l}$.
Suppose also that
\begin{enumerate}
\item[(i)] $\alpha \in \mcH_O^p$ and $\alpha$ is projected to $\beta$ in the fiber.
    \item[(ii)]
 $u$ has coisotropic
regularity of order $k$ relative to $H^m$ (on the coisotropic $\dot{\mathcal{F}}^{p}_I$
(resp. $\dot{\mathcal{F}}^{s}_I$)) in an open set containing all points in $\dot{\mathcal{F}}^{p}_{I,\beta}\cap \{ 0 < x < \delta\}$( resp.
$\dmcF^{s}_{I,\beta} \cap \{ 0 < x < \d\}$)
that are geometrically related to $\alpha$.
   \item[(iii)]
 $Pu$ has coisotropic
regularity of order $k$ relative to $H^{m-1}$ in a neighborhood of $(\mcF^p_{O,\alpha})^o$ (resp. $(\mcF^s_{O,\alpha})^o$). In addition, there exists a neighborhood of $\{ x = 0\}$ such that $Pu$ has coisotropic regularity of order $k$ relative to $H^{m-1,l-2}_e$ on this neighborhood.
\end{enumerate}
Then $u$ has
coisotropic regularity of order $k$ relative to $H^{m'}$ for all
$$ m' < \text{min}(m-1/2, l + f/2)$$
(on the coisotropic $\dot{\mathcal{F}}^{p}_O$ (resp. $\dot{\mathcal{F}}^{s}_O$)) in a neighborhood of
$\mathcal{F}^{p}_{O,\alpha}$ (resp. $\mathcal{F}^{s}_{O,\alpha}$) strictly away from $\partial M$.
\end{theorem}

\begin{rem}
The assumption on $Pu$ near $\{ x= 0\}$ is there not just to apply Theorem \ref{theorem:Coisopropagation in/out of edge}, but to also allow us to propagate regularity along the edge $\{ x = 0\}$ away from $\mcH^p_{I/O}$. Indeed, since the module $\mcM^p$ has elliptic elements at such points, the assumption on $Pu$ along this region just means $Pu \in H^{m+k-1,l-2}_e$ microlocally at such points, which will allow us to use standard propagation of singularities as we see in the proof.
\end{rem}

As an immediate consequence, one also gets the following corollary where we don't distinguish among geometric rays:
\begin{cor}
\label{thm:CoisoVersion2}
Let $\alpha \in \dot{ \mathcal{H} }$, and $k \in \mathbb{N}$. Suppose that
$u \in H_e^{-N,l}$ and $Pu$ is coisotropic of order $k$ relative to $H_e^{m-1,l-2}$ in a neighborhood of $\{ x = 0\}$ and on a neighborhood of $\dot{\mcF}^p_{O,\alpha}$. Suppose also that $u$ has coisotropic
regularity of order $k$ relative to $H^m$ (on the coisotropic $\dot{\mathcal{F}}^{p}_I$
(resp. $\dot{\mathcal{F}}^{s}_I$)) near $\dot{\mathcal{F}}^{p}_{I,\alpha}$( resp.
$\dot{\mathcal{F}}^{s}_{I,\alpha}$) strictly away from $\partial M$. Then $u$ has
coisotropic regularity of order $k$ relative to $H^{m'}$ for all
$$ m' < \text{min}(m-1/2, l + f/2)$$
(on the coisotropic $\dot{\mathcal{F}}^{p}_O$ (resp. $\dot{\mathcal{F}}^{s}_O$)) near
$\dot{\mathcal{F}}^{p}_{O,\alpha}$ (resp. $\dot{\mathcal{F}}^{s}_O$) strictly away from $\partial M$.
\end{cor}

{\it( Proof of Theorem \ref{thm:CoisoVersion1} )} We follow closely a clever argument used in \cite{MVWEdges} and \cite{MVWCorners}.
First, set $\tilde{l} = \min(l, m-f/2 -1/2- 0)$ so that $u \in H_e^{-N,\tilde{l}}$.
Notice that if
$l < m- f/2 -1/2 - 0, \text{ then } m > f/2 +l +1/2 = \tilde{l}+f/2 + 1/2$. On the other hand, if
$ m- f/2 -1/2- 0 \leq l, \text{ then } \tilde{l} + f/2 +1/2= m-f/2-1/2-0+1/2+f/2 = m-0 <m$ so this shows
$m > \tilde{l} + f/2 +1/2$.
 Hence, Theorem \ref{theorem:Coisopropagation in/out of edge} (1) is applicable (applied to each $\alpha'$ geometrically related to $\alpha$, with $\alpha$ in the theorem replaced by $\alpha'$), and one may then propagate along the edge and back out of the edge by the same argument as in \cite{MVWEdges}.
$\Box$

As in \cite[Corollary 8.3]{MVWCorners}, we prove that the regularity with respect to which coisotropic regularity is gained in the above results is, in fact, independent of the weight $l$:

\begin{cor}
\label{Cor:PropOfCoisotropy} Assume $f>1$.
Let $\beta \in \dot{ \mathcal{H} }^p$, $\epsilon>0$ and $k \in \mathbb{N}$. Suppose that $u \in H^{s-1}_{b,\tmcD,loc}(I\times X)$ when restricted to $I\subset (\infty,t(\beta)) $, a precompact time interval, and $Pu \in H^s_{b,\tmcD',loc}(M)$. Then there
is a $k'$ (depending on $k$ and $\epsilon$ ) such that if
 \begin{enumerate}
 \item[(i)]$\alpha \in \mcH_O^p$ and $\alpha$ is projected to $\beta$ in the fiber.
     \item[(ii)]
 $u$ has coisotropic regularity of order $k'$ relative to $H^s$ (on the coisotropic $\mathcal{F}^{p}_I$
(resp. $\mathcal{F}^{s}_I$) in an open set containing all points in $\dot{\mathcal{F}}^{p}_{I,\beta}$( resp.
$\dot{\mathcal{F}}^{s}_{I,\beta}$) strictly away from $\partial M$ that are geometrically related to $\alpha$,
\item[(iii)] $Pu$ has coisotropic regularity of order $k'$ relative to $H^{s-1}$ (on the coisotropic $\mathcal{F}^{p}$
(resp. $\mathcal{F}^{s}$)) in an open set containing all points in $\mathcal{F}^{p}_{O,\alpha}$( resp.
$\mathcal{F}^{s}_{O,\alpha}$) strictly away from $\partial M$. Also, there exists a neighborhood of $\{ x=0 \}$ such that
$Pu \in H^{N}_{b,\mcD'}(M)$ for any $N -1> s + k'$,
\end{enumerate}
 then $u$ has
coisotropic regularity of order $k$ relative to $H^{s-\epsilon}$ (on the
coisotropic $\mathcal{F}^{p}_O$ (resp. $\mathcal{F}^{s}_O$)) in a neighborhood of
$\mathcal{F}^{p}_{O,\alpha}$ (resp. $\dot{\mathcal{F}}^{s}_{O,\alpha}$) strictly away from $\partial M$.
\end{cor}

We also record an immediate consequence of the above corollary where we do not distinguish between geometric and diffractive rays.

\begin{cor}Assume $f>1$.
Let $\alpha \in \dot{ \mathcal{H} }$, $\epsilon>0$ and $k \in \mathbb{N}$.  Suppose $u \in H^{s-1}_{b,\tmcD,loc}(I\times X)$ when restricted to $I\subset (\infty,t(\alpha)) $, a precompact time interval, and $Pu \in H^{s}_{b,\tmcD',loc}(M)$.  Then there
is a $k'$ (depending on $k$ and $\epsilon$ ) such that if
\begin{enumerate}
\item[(i)]
$u$ has coisotropic regularity of order $k'$ relative to $H^s$ (on the coisotropic $\dot{\mathcal{F}}^{p}_I$
(resp. $\dot{\mathcal{F}}^{s}_I$)) in a neighborhood of $\dot{\mathcal{F}}^{p}_{I,\alpha}$( resp.
$\dot{\mathcal{F}}^{s}_{I,\alpha}$) strictly away from $\partial M$
\item[(ii)]
$Pu$ has coisotropic regularity of order $k'$ relative to $H^{s-1}$ on a neighborhood of $\dot{\mcF}^p_{O,\alpha}$ strictly away from $\partial M$. Also, there exists a neighborhood of $\{ x=0 \}$ such that
$Pu \in H^{N}_{b,\mcD'}(M)$ for any $N -1> s + k'$ on this neighborhood,
\end{enumerate}
then $u$ has coisotropic regularity of order $k$ relative to $H^{s-\epsilon}$ (on the coisotropic $\dot{\mcF}^p_{O}$ (resp. $\dot{\mcF}^s_{O}$) in a neighborhood of $\dot{\mcF}^p_{O,\alpha}$ ( resp. $\dot{\mcF}^s_{O,\alpha}$) strictly away from $\p M$.
\end{cor}

\begin{rem}
Note that the assumption $u \in H^{s-1}_{b,\tmcD, loc}(I \times X)$ is essential for these corollaries. Indeed, Theorem \ref{thm:CoisoVersion1} only allows us to propagate coisotropic regularity into and out of the edge only with respect to a Sobolev space of low order. Hence, this assumption on $u$ will allow us to improve this order by means of an interpolation argument. Such an assumption would not be needed if we had a $b$-propagation result (i.e. a version of Snell's Law telling us that the angular momentum is preserved along outgoing rays generated by the incident ray) that would already tell us that $u$ is microlocally in $H^s$ along $(\mcF^p_{O,\alpha})^o$.
\end{rem}

{\it (Proof of Corollary \ref{Cor:PropOfCoisotropy}).} \ The key is to identify a weighted edge space that $u$
lies in. The assumption that $u \in H_{b,\tmcD,loc}^{s-1}(I \times X)$ and $Pu \in H_{b,\tmcD',loc}^{s}(M)$ implies by Theorem \ref{cor:speed semiglobal propagation}
\begin{equation}\label{eq:coiso2 u in H_b,D^(s-1)}
 u \in H_{b,\tmcD,loc}^{s-1}(M),
 \end{equation}
i.e. $u$ is in such a space for all times and not just restricted to the interval $I$. One then shows that
$$ u \in H^{s-1,-(f+1)/2}_{e,loc}(M) $$ to be able to apply the previous corollaries to prove a low edge-regularity propagation result along the outgoing $p$-ray, followed by an interpolation argument. The details are almost identical as in the proof of \cite[Theorem 12.1]{MVWEdges}, \cite[Corollary 8.3]{MVWCorners}, and \cite[Corollary 10.1.20]{VKThesis}.
$\Box$

The next several sections will be devoted to establishing a duality result of the previous theorem. The first step will be to establish a vital energy estimate, which we do in the next section.

\section{Energy Estimates}
In this section, we state the crucial energy estimates that are at the heart of dualizing the propagation of coisotropic regularity argument. The proofs here are standard arguments for obtaining energy estimates for hyperbolic equations, so we do not go into the proofs here as not much novelty would be offered. We simply refer the interested reader to \cite[Chapter 11]{VKThesis} for full, comprehensive proofs.

\begin{theorem}(Finite Speed of Propagation \cite[Corollary 11.4.2]{VKThesis})\label{cor:speedFiniteSpeedOfProp}
Let $w_0 \in \p X$ and $K :=\bar{B}_r(w_0)$ a closed geodesic ball of radius $r>0$ around $w_0$ using the metric. Denote $\kappa := [\text{sup}_K(\lambda + 2\mu)]^{-1}$ and fix any sufficiently small $\delta,\epsilon_0 >0$ and times $T_0^- <T_0 <T_1-\epsilon_0 <T_1$ such that $0 < \kappa^{-1}( (T_1-T_0^-)-\epsilon^2_0 ) < r$ and $0 < \kappa^{-1}( (T_1-T_0^-)-\epsilon^2_0 -\d )$. If
$$u \in H^{-N}_{\tD,b}, \ Pu \in H^{m+1}_{\tD',b}
\text{ on the set} $$
$$\{(t, p):\ d^2(p,w_0) \leq \kappa^{-1} (T_1 -t)^2, \ T_0^- \leq t \leq T_1 - \epsilon_0 \}$$
and
$$u \in H^m_{\tD,b} \text{ on } \{ (t,p) :\ d^2(p,w_0) \leq \kappa^{-1} (T_1 -t)^2 , \ T_0^- \leq t \leq T_0 \},$$
then for any compact interval $I \subset (T_0^-,T_1-\epsilon_0)$,
$$u \in H^m_{\tD,b} \text{ on }
\{ (t,p)  :\ d^2(p,w_0) \leq \kappa^{-1}( (T_1 -t)^2 -\delta), t \in I \}$$
Moreover, the following estimate holds
\begin{align*}
&||u||^2_{H^m_{\tD,b}( \{ (t,p) \in I \times X : \ d^2(p,w_0) \leq \kappa^{-1}( (T_1 -t)^2 -\delta)\} ) }
\\
&\leq C\Big( ||Pu||^2_{H^{m+1}_{\tD',b}( \{ (t,p) \in [T_0^-,T_1] \times X : \ d^2(p,w_0) \leq \kappa^{-1} (T_1 -t)^2\} ) }
\\
&\qquad \qquad
+ ||u||^2_{H^m_{\tD,b}( \{ (t,p) \in [T_0^-,T_0] \times X : \ d^2(p,w_0) \leq \kappa^{-1} (T_1 -t)^2 \} ) }
\\
&\qquad \qquad
+||u||^2_{H^{-N}_{\tD,b}( \{ (t,p) \in [T_0^-,T_1] \times X : \ d^2(p,w_0) \leq \kappa^{-1}(T_1 -t)^2 \} ) }
\Big )
\end{align*}
\end{theorem}

We also have a useful estimate giving us a semi-global propagation of $b$ regularity.
\begin{theorem}(\cite[Corollary 11.4.3]{VKThesis})\label{cor:speed semiglobal propagation}
Fix times $T_0^- < T_0 < T_1$. Suppose that
\begin{enumerate}
\item[(i)]
$u \in H^{-N}_{\tD,b}((T_0^-,T_1)\times X)$
 \item[(ii)]
 $Pu \in H^{m+1}_{\tD',b}((T_0^-,T_1)\times X)$
  \item[(iii)]
  $u \in H^{m}_{\tD,b}((T_0^-,T_0)\times X)$ when restricted to $T^-_0 \leq t \leq T_0$.
\end{enumerate} Then for any time interval $I \subset (T_0^-,T_1)$ we in fact have $u \in H^{m}_{\tD,b}(I \times X)$, and moreover, the following estimate holds
\begin{align*}
&||u||^2_{ H^{m}_{\tD,b}((T_0^-+ \d, T_1 -\epsilon)\times X)}
\\
&\qquad\leq C \Big(||u||^2_{H^{-N}_{\tD,b}((T_0^-,T_1)\times X)} + ||Pu||^2_{H^{m+1}_{\tD',b}((T_0^-,T_1)\times X)}
 + ||u||^2_{H^{m}_{\tD,b}((T_0^-,T_0)\times X)} \Big)
\end{align*}
with a constant $C$ that depends on $\epsilon, \d$.
\end{theorem}

Next, using Theorem \ref{cor:speed semiglobal propagation} to propagation $b$-regularity in time, one may deduce the following theorem exactly as done in \cite[Theorem 8.12]{VasyAds} and in \cite[Volume 3]{Hor}:
\begin{theorem}(\cite[Theorem 11.2.9]{VKThesis})\label{thm:energyFundSolutionExistence}Suppose $dim(Z)>1$, and $t_0 < t_1' < t_1$.
Let $m' \in \RR$ and $$f \in  H^{m+1}_{b,\tmcD'}(M;TX)^{\bullet}_{(t_0,\infty)}$$ for some $m \geq m'$. Then there exists a unique $u \in H^{m'}_{b,\tmcD}(M;TX)^{\bullet}_{[t_0,\infty)}$ solving $Pu =f$.
Moreover, one in fact has
$u \in H^{m}_{b,\tmcD}(M;TX)$, and the estimate
 \begin{equation}\label{eq:bReg main existence estimate}
 ||u||_{H^{m}_{b,\tmcD}(M)|_{[t_0,t'_1]}} \lesssim ||Pu||_{H^{m+1}_{b,\tmcD'}(M)|_{(t_0,t_1]}}.
 \end{equation}
\end{theorem}

Finally, we may prove the global result:
\begin{theorem}(\cite[Theorem 11.4.4]{VKThesis})\label{thm:EnergyEstimate}
Let $m \in \RR$ and suppose $\kappa = \sup_X(\lambda + 2\mu)< \infty$.
Given $f \in H^{m+1}_{b,\tD',loc}(M)^{\bullet}_{(t_0,\infty)}$,
 there exists a unique forward solution $u \in H^{m}_{b,\tD,loc}(M)^{\bullet}_{[t_0,\infty)}$ such that $Pu =f$.
\end{theorem}
\proof
The proof is as done in \cite[Volume 3]{Hor} with full details in our setting in \cite[Theorem 11.4.4]{VKThesis}.
 $\Box$

By use of these energy estimates, we will prove propagation of coinvolutive regularity and our main result on diffraction in the next section.

\section{Dualization and Diffraction of Elastic Waves}
In this section, we dualize the results in Section \ref{sec: cosio regularity} to obtain the propagation of coinvolutive regularity through the edge. Throughout this section, we assume that all operators constructed here are scalar unless mentioned otherwise, and that $f>1$.\\
\indent Next, it will be convenient to pick out certain subsets of $\dot{\mcF}$ that restrict the flow of bicharacteristics to certain time intervals, so we introduce the following definition.

\begin{definition}
Let $\alpha \in \dot{\mcH}^{p/s}$ and $\I$ a compact time interval either contained in $(t(\alpha), \infty)$ or $(-\infty,t(\alpha))$. Denote by
 \begin{equation*}
 \dot{\mcF}^{p/s}_{\bullet,\alpha}(\I) = (\dot{\mcF}^{p/s}_{\bullet,\alpha})^o \cap t^{-1}(\I) \subset (\dot{\mcF}^{p/s}_{\bullet,\alpha})^o , \qquad
\text{ with $\bullet = $ $I$ or $O$}
\end{equation*}
which is the set of all points in $ \dot{\mcF}^{p/s}_{\bullet,\alpha}$ whose $t$ coordinate is an element of $\I$. We always assume that $\I$ is picked close enough to $t(\alpha)$ so that $\mcF^{p/s}$ exist as smooth coisotropic manifolds over $\I$.
\end{definition}

\begin{theorem}
\label{thm:propagation of coinvolutive}
Let $u \in H^{-\infty}_{b,\tD,loc}(M)$ be a solution to the elastic equation. Take $\alpha \in \dot{\mcH}^p$.
Assume there is a time interval $\I \subset (-\infty, t(\alpha))$, $k \in \mathbb{N}$, and $B_0 \in \Psi_c^0(M^o)$ elliptic on $\dot{\mcF}^p_{I,\alpha}(\I)$
such that
\begin{equation}
 u \in H^{s+1}_{b,\tmcD',loc}(\I \times X) + B_0v \text{ where } B_0 v \in  \mcA^k_p(H^{s+1}_{b,\tmcD'}(\I \times X))
\end{equation}
 when $u$ is restricted to $\I\times X$ (as a distribution). Then there exists a $k' \in \mathbb{N}$ such that $u$ is coinvolutive on $(\dot{\mathcal{F}}^p_{O,\alpha'})^o$ with respect to $H^r$ for any $r < s$.
\end{theorem}

\begin{rem}
The assumption on $u$ just means that during some time interval before the $p$-bicharacteristics in $\dot{\mcF}^p_{I, \alpha}$ reach $\alpha$, $u$ is globally nice, i.e.\ $H^s$ everywhere, except that it is slightly worse on a subset of $(\dot{\mcF}^p_{I,\alpha})^o$ by being coinvolutive there.
\end{rem}

\begin{rem}
The essential idea of the proof is as in \cite[Remark 9.3]{MVWCorners}, which we will adapt to our setting. Our previous results show that under certain assumptions, coisotropic regularity entering the edge along $\dot{\mcF}^p_{I,\alpha}$ imply coisotropic regularity along $\dot{\mcF}^p_{O,\alpha}$. In other words, regularity under the application of $A^{\gamma}$ along
$(\dot{\mcF}^p_{I,\alpha})^o$ with respect to $H^{-s}$ combined with assumed semi-global regularity by being in $H^{-s-1}_{b,\tD}(\I \times X)$, yields regularity under $A^{\gamma'}$ along $(\dot{\mcF}^p_{O,\alpha})^o$ w.r.t. $H^{-s-\epsilon}$. The dual condition to this is lying in the range of the operators $A^{\gamma}$ in the relevant regions. Thus, using time reversal and duality, the condition of coinvolutive regularity along $\dot{\mcF}^p_{I,\alpha}$ w.r.t. $H^{s}$, i.e. lying in the range of $A^{\gamma}$ microlocalized there, combined with being in $H^{s+1}_{b,\tmcD'}$ elsewhere for some short time interval projected from $(\dot{\mcF}^p_{I,\alpha})^o$, leads to coinvolutive regularity along $(\dot{\mcF}^p_{O,\alpha})$ w.r.t. $H^{s-\epsilon}$, i.e. lying in the range of $A^{\gamma'}$ microlocalized there.
\end{rem}

\proof We follow closely the proof of \cite[Theorem 9.2]{MVWCorners}. We first assume $s \leq 0$ to simplify notation; we'll return to the general case at the end of the argument. Fix an $\epsilon >0$. Let $T = t(\alpha)$, and choose $T_0 < T < T_1$ so that $\I \subset (T_0,T)$ and $T_1$ is close to $T$ as will be specified later. Let $\chi$ be smooth step function such that $\chi \equiv 0$ on a neighborhood of $(-\infty, T_0]$ and $\chi \equiv 1$ on a neighborhood of $(T, \infty)$, so that one also has $\supp(d\chi) \subset \I$.
We find that
$$ v \equiv \chi u$$
satisfies
$$Pv = f \text{ with }f = [P,\chi]u.$$
 Since $d\chi$ is supported in $I$, our initial hypothesis implies that
 $$f= [P,\chi]u_0 + [P,\chi]u_1 := f_0 + f_1,$$
$$\text{for some }u_0 \in H^{s+1}_{b,\tmcD'}(\I \times X),$$
 $$\text{ and }u_1 \text{ is coinvolutive of order $k$ w.r.t. $H^{s+1}_{b,\tmcD'}$ on $\dot{\mcF}^p_{I,\alpha}(\I)$}.$$
 Also, notice that $v$ vanishes on a neighborhood of $(-\infty, T_0] \times X$. We will then write
$$ v = P^{-1}_+f$$
as the unique forward solution to the equation $P\phi = f$.
Observe that since $\p_t$ is a $b$-differential operator, then $[P,\chi] \in \Diff_b^{1}(M;TX)$
 $$\Rightarrow f_0 = [P,\chi]u_0 \in H_{b,\mcD'}^{s}(M)$$
 since the support condition on $d\chi$ allows us to extend $f_0$ to all of $M$. Then by Theorem \ref{thm:EnergyEstimate}, there exists a unique forward solution
 $$P^{-1}_+f_0 \in H^{s-1}_{b,\mcD}(M;TX), \text{ satisfying }P(P^{-1}_+f_0)=f_0.$$

 In particular, $P^{-1}_+f_0 \in H^s$ microlocally at $w$ since $w$ is away from $\p M$. By definition, this is certainly a stronger condition than $P^{-1}_+f_0$ being coinvolutive relative to $H^{s-\epsilon}$ at $w$; hence we conclude that $ P^{-1}_+f_0 $ is coinvolutive of order $k'$ for some $k'$ to be determined relative to $H^{s-\epsilon}$ at $w$. Since $v \equiv u$ for $t>T$, all we need to show is that the \emph{unique} forward solution, denoted $P^{-1}_+f_1$, of $Pv_1 = f_1$ is coinvolutive of order $k'$ relative to $H^{s-\epsilon}$ at $w$. The proof of this proceeds identically to the proof of \cite[Theorem 9.2]{MVWCorners}. All that is necessary is propagation of coisotropic regularity and an energy estimate, both of which have been proven, followed by a duality argument using the Hahn-Banach theorem. $\Box$

We may now combine the result on propagation of coisotropic regularity with the coinvolutive propagation to show that regularity of the solution to the elastic equation propagates precisely along the geometric rays. Our main theorem is then
\begin{theorem}\label{thm:Diffraction2}

 Let $u \in H^{S-1}_{b,\tmcD}(M;TX)$, $q \in \dot{\mathcal{H}}^p$, and $\alpha \in \mathcal{H}^p$ projecting to $q$. Suppose there is a compact time interval $\I \subset (-\infty, t(\alpha))$ and $B_0 \in \Psi^0_c(M^o)$ that is elliptic on $\dot{\mcF}^p_{I,q}(\I)$, for which
 $u \in H^{s+1}_{b,\tmcD'}(\I \times X) + B_0v$ for a distribution $v$ such that $B_0 v \in  \mcA^k_p(H^{s+1}_{b,\tmcD'}(\I \times X))$
 when $u$ is restricted to  $I \times X$. Then,
 \begin{align*}
 (\mathcal{F}^p_{I,\alpha'})^o \cap \text{WF}^s(u) = \emptyset \ \forall \alpha' \text{ geometrically related } &\text{to }\alpha \Rightarrow \\
 &(\mathcal{F}^p_{O,\alpha})^o \cap \text{WF}^r(u) = \emptyset \ \forall r<s.
\end{align*}
(Again we always assume $\I$ is close enough to $t(\alpha)$ so that $\mcF^p_I$ is a smooth manifold over $\I$.)
\end{theorem}

\proof
With the hypothesis of the theorem, pick any $ w \in (\mcF_{O,\alpha}^p)^o$.
Here, we follow the proof of Theorem \ref{thm:propagation of coinvolutive} to separate $u$ into microlocalized pieces. Thus, let $\chi$, $T_0 < T:=t(\alpha) < T_1$, and $v = \chi u$ as in that proof. Also, one has
$$ Pv = f_0 + f_1 $$
with $f_0 \in H^{s}_{b,\tmcD'}(\I \times X)$ and $f_1$ coinvolutive of order $k$ on $\dot{\mcF}^p_{I,\alpha}(\I)$ w.r.t. $H^{s}_{b,\tmcD'}(\I \times X)$.
By Theorem \ref{thm:EnergyEstimate}, $P^{-1}_+ f_0 \in H^{s-1}_{b,\tmcD}(M;TX)$ which in particular, since $w$ is away from $\p M$,
$$ w \notin \WF^s( P^{-1}_+f_0 ).$$
Thus, all we need to show is that $w \notin \WF^r( P^{-1}_+f_1 )$ $\forall r<s$ since $u$ agrees with $v$ for $t>T$. It will be convenient to denote $I':=\supp(d\chi) \subset \I$, where by construction $I' \pm \d \subset \I$ for some $\d >0$. Now, the assumption that $\WF^s(u) \cap \bigcup_{ \{\alpha': \alpha' \sim_G \alpha \}} (\mcF^p_{I,\alpha'})^o = \emptyset $ implies that
 \begin{equation}\label{eq:coinv f_1 is H^(s-1) near alpha' curves}
 \WF^{s-1}(f_1) \cap \bigcup_{ \{\alpha' \in \mcH^p_I: \alpha' \sim_G \alpha \}} (\mcF^p_{I,\alpha'})^o = \emptyset
 \end{equation}
when we write $f_1 = [P,\chi]u - f_0$ since $f_0$ satisfies this condition, $[P,\chi]$ is a first order differential operator, and the microlocality of $[P,\chi]$.
Thus, let $\ti{B} \in \Psi_b^{0}(M)$ scalar, with Schwartz kernel supported in $(\I \times X^o)^2$, elliptic on $\mcF^p_{I,\alpha'}(I')$ for each $\alpha'$ geometrically related to $\alpha$  such that
\begin{equation}\label{eq:coinv I-ti(B) microsupported away from alpha'}
\WF'(\Id - \ti{B}) \cap \mcF^p_{I,\alpha'}(I') = \emptyset \text{ for } \alpha' \sim_G \alpha
\end{equation}
(this is possible exactly because $I' \pm \d \subset \I$). We also make $\ti{B}$ microsupported in a small enough neighborhood of $\dot{\mcF}^p_{I,\alpha}$ so that
$$ \ti{B}f_1 \in H^{s-1} $$
using (\ref{eq:coinv f_1 is H^(s-1) near alpha' curves}).
So again using Theorem \ref{thm:EnergyEstimate} shows $P^{-1}_+ \ti{B}f_1$ is microlocally in $H^{s}$ at $w$. Hence, we are left with understanding $P^{-1}_+(\Id - \ti{B})f_1$.

Let $\alpha'$ be such that $\pi_0(\alpha') = q$ and is geometrically related to $\alpha$. Since, $f_1$ is supported in $I' \times X$, it vanishes on $ (\I \setminus I') \times X$. Also $(\Id-\ti{B})f_1= f_1-\ti{B}f_1$ is supported in $\I \times X$, so by microlocality and (\ref{eq:coinv I-ti(B) microsupported away from alpha'}), one has
\begin{align*}
&\WF^{\infty}((\Id - \ti{B})f_1) \cap (\mcF^p_{I,\alpha'})^o
 =\emptyset.
\end{align*}
Thus, if we set
$$ u_1 = P^{-1}_+(\Id - \ti{B})f_1, $$ then using propagation of coisotropic regularity over $M^o$ Corollary \ref{cor:module interior coiso for forward solution}(with $k=\infty$ and $s=S$ in the statement of the corollary) shows that
$$ \WF^S(\mcA^{\ti{k}} u_1 ) \cap (\mathcal{F}^p_{I,\alpha'})^o = \emptyset \quad
\forall \ti{k} \in \mathbb{N}, \ \alpha' \sim_G \alpha. $$

\indent Next, $(\Id - \ti{B})f_1$ vanishes for $t > T^-$ for some $T^- <T$ implies that $\WF^{\infty}(f_1)$ is disjoint from $(\dot{\mcF}^p_{O,\alpha})^o$, and $f_1 \in H^{\infty}_{b,\mcD'}$ microlocally in a neighborhood of $\p \dot{\mcF}^p_{O,\alpha}$. Thus, Corollary \ref{Cor:PropOfCoisotropy} applies to show that microlocally at $w$ one has
\begin{equation}\label{eq:coinv u_1 infinitely coisotropic}
 A^{\gamma}u_1 \in H^{S-0} \ \forall \gamma \text{ multiindex}.
 \end{equation}
Analogously, we may use Corollary \ref{cor:module interior coiso for forward solution} ( with $s$ and $k$ in the notation of that corollary being the same $s$ and $k$ here ) over $M^o$ to conclude that $u_1$ has coinvolutive regularity of order $k$ on $ (\dot{\mcF}^p_{I,q})^o $ w.r.t. $H^s$. Then the regularity properties just described for $f_1$ near $\dot{\mcF}^p_{O,\alpha}$ allow us to apply Theorem \ref{thm:propagation of coinvolutive} to conclude that microlocally near $w$
\begin{equation}\label{eq:coinv u_1 is coinv order k' w.r.t. H^r}
 u_1 \in \sum_{|\gamma|\leq k'} A^{\gamma}(H^r).
 \end{equation}

The rest of the proof is an interpolation argument, using an FIO to turn elements of $\mcA$ into a model form and then interpolating between infinite coisotropic regularity and low Sobolev regularity. The exact argument may be found in the proof of \cite[Theorem 9.6]{MVWCorners} to show $u_1 \in H^{r-\epsilon}$ microlocally at $w$. Since $r<s$ and $\epsilon >0$ was arbitrary, this gives the desired result.
$\Box$

\section{Application to the fundamental solution}
Now consider the fundamental solution
$$ u = P^{-1}_+\delta=P^{-1}_+ (\delta_{(t_0,o)} \otimes \Id) \text{ for }(t_0,o) \in M^o.$$
\begin{cor}
For all $(t_0,o) \in M^o$ let $\mathcal{L}^p_{t_0,o}$ (resp. $\mathcal{L}^s_{t_0,o}$) denote the flowout of $SN^*(\{o\})$ along $p$-bicharacteristics, resp. $s$-bicharacteristics, lying over $M^o$, which lie over $o$ and $t=t_0$. If $o$ is sufficiently close to $\partial X$, then for some short time beyond when the first $p$-bicharacteristic lying over $(t_0,o)$ reaches the boundary, the forward fundamental solution $u = u_{t_0,o}$ is a Lagrangian distribution along $\mathcal{L}_{t_0,o} := \mathcal{L}^p_{t_0,o} \cup
\mathcal{L}^s_{t_0,o}$ lying in $H^s$ for all $s < -n/2 +1$ together with diffracted waves, singular only at $\mcF_O^p \cup \mcF_O^s$, that lie in $H^r$ for all $r < -n/2 +1 +f/2$, away from its intersection with $\mathcal{L}_{t_0,o}$.
More precisely, if we consider the first incoming $p$-wave transverse to the boundary, i.e. $u \in H^{-n/2+1-0}$ along $\mcF^p_{I,\alpha}$ but $\WF^{\infty}(u) \cap (\mcF^p_{I,\alpha})^o \neq \emptyset$; then each of the outgoing diffracted $p$ and $s$ waves are weaker in the sense that
$u \in H^r$ along the diffracted $p$ and $s$ bicharacteristics generated by $\mcF^p_{I,\alpha}$ for all $r < -n/2 + 1 + f/2$. Similarly, if we consider the first incoming $s$-wave transverse to the boundary, i.e. $u \in H^{-n/2+1-0}$ along $\mcF^s_{I,\alpha}$, then each of the outgoing diffracted $p$ and $s$ waves are weaker in the sense that
$u \in H^r$ along the diffracted $p$ and $s$ bicharacteristics generated by $\mcF^s_{I,\alpha}$ for all $r < -n/2 + 1 + f/2.$
\end{cor}

\proof
Let $w \in \mcF^p_{O} \cup \mcF^s_{O}$ not geometrically related to any point in $\mcF^p_{I,\alpha}$ (by definition, no point in $\mcF^s_O$ is geometrically related to $\mcF^p_{I,\alpha}$, so we get the diffractive improvement on all such points).
In order to consider $p$-waves and $s$-waves separately, using that $\Sigma_p \cap \Sigma_s = \emptyset$, pick $Q \in \Psi^0_c(M^o)$ whose wavefront set is disjoint from $\Sigma_s$ and is elliptic on $S^*_{(t_0,o)}M \cap \Sigma_p$. Then write
$$ \delta = Q\delta + (I-Q)\delta = \delta_p + \delta_s$$
and denote
$$ u_p = P^{-1}_+\delta_p, \ \ u_s = P^{-1}_+ \delta_s.$$

We will show $u_p$ and $u_s$ each belong to $H^{-n/2+1+f/2-0}$ microlocally at $w$.

{\bf \underline{Step 1: Obtaining regularity for $u_p$} }

 By construction of $u_p$, using also that it is a forward solution, we obtain by propagation of singularities in the form of Corollary \ref{cor:module interior coiso for forward solution} that
\begin{equation}\label{u_p has no s-singularities for short time}
 \WF(u_p) \cap \mcL^s = \emptyset
 \end{equation}
$$ \Rightarrow \WF(u_p) \cap (\mcF^s_I)^o = \emptyset.$$

Denote $\bar{t}= t(\alpha)$. We are assuming that this is the first time a $p$-bicharacteristic strikes the boundary.
We will show the diffractive improvement for $u_p$ first, as $u_s$ involves a separate argument.

 We want to further break up $\u_p$ into a piece microsupported along the incoming bicharacteristic segment $\mcF^p_{I,\alpha}$ and a piece microsupported away from this ray. To this end, let $Q_0 \in \Psi_c^0(M)$ have microsupport in a small neighborhood (to be determined later) of $(\mcF^p_{I,\alpha})^o$ which is elliptic on a smaller neighborhood of $(\mcF^p_{I,\alpha})^o$. Consider the decomposition
$$ \delta_p = Q_0 \delta_p + (I-Q_0)\delta_p \equiv \delta_{p,0} + \delta_{p,1},$$
and
$$ u_p = P^{-1}_+(\delta_{p,0}) + P^{-1}_+(\delta_{p,1})
\equiv u_0 + u_1.$$

{\bf \underline{Step 1.1: Getting regularity of $u_0$ piece, where $u_p = u_0 + u_1$}}

By the microlocality property of PsiDO's, one has
\begin{equation}\label{eq:fund delta_(p,0) WF set}
\WF(\d_{p,0}) \subset S^*_{(t_0,o)}M \cap \WF'(Q) \cap \WF'(Q_0).
\end{equation}
Notice that if $w \in \mcF^s_{O,\beta}$ for some $\beta \in \mcH^s_O$, then (\ref{u_p has no s-singularities for short time}) already implies that
$$ \WF(u_0) \cap (\mcF^s_{I,\beta'})^o = \emptyset \text{ for }\beta' \sim_G \beta.$$
 On the other hand suppose $w \in \mcF^p_{O,\beta}$ for some $\beta \in \mcH^p_O$ that lies in the same fiber as $\alpha$.
Notice that if $\beta'\in \mcH^p_I$ is geometrically related to $\beta$ so that  $ \mcF^p_{I,\beta'}$ are points geometrically related to $w$, then $\mcF^p_{I,\beta'}$ does not lie over $(\bar{t},o)$ by the assumption that $w$ is not geometrically related to any point in $\mcF^p_{I,\alpha}$.
Thus,
$$ \WF(\d_{p,0}) \cap \mcF^p_{I,\beta'} =\emptyset \text{ for } \beta' \sim_G \beta. $$
As before, this implies using Corollary \ref{cor:module interior coiso for forward solution} and since $u_0$ is smooth for $t< t_0-\d$ that
\begin{equation}\label{eq:fund u_o nice on geom related to beta'}
\WF(u_0) \cap (\mcF^p_{I,\beta'})^o \text{ for } \beta' \sim_G \beta.
\end{equation}
Hence, taking any $q \in \dmcH^p$ such that $\beta$ lies in the same equivalence class as $q$, we need to show that $u_0$ satisfies the coinvolutive hypothesis of Theorem \ref{thm:Diffraction2} on $\dot{\mcF}^p_{I,q}$ in order to apply that theorem to $u_0$ to get improved regularity at $w$.

To proceed, it is well known that inside $(\infty,\bar{t}) \times M^o$ (i.e., before the first $p$-bicharacteristic hits $\p M$) $u_0$ is a Lagrangian distribution associated to the Lagrangian
\begin{align}\label{eq:fund Lagrangian for u_{p,0}}
 \mcL^p_0 &:= \mcL( \WF(\d_{p,0}) \cap \Sigma )
 \nonumber \\
 &= \mcL^p( \WF(\d_{p,0}) \cap \Sigma_p )
 \nonumber \\
 &\subset \mcL^p(S^*_{(t_0,o)}M^o \cap \WF'(Q_0)\cap \Sigma_p)
 \end{align}
 where $\mcL(K)$ refers to the flowout from $K \subset S^*M^o$ of both $p$ and $s$ bicharacteristics. \footnote{This can be proven by the argument we used to prove coisotropic regularity in Section \ref{sec: cosio regularity} by reducing $P$ to a standard wave equation, and then invoking a well-known wave equation result, which states that the forward fundamental solution is a Lagrangian distribution associated to the flowout of the Hamilton vector field.} The equality above follows from (\ref{eq:fund delta_(p,0) WF set}) and since $\WF'(Q) \cap \Sigma_s = \emptyset.$ Notice, $\mcL^p_0$ may be visualized as a conic spray of $p$-bicharacteristics that are close to $\mcF^p_{I,\alpha}$. In fact, $u_p$ is a Lagrangian distribution associated to $\mcL^p_0$ of order $s' = n/4-5/4$ (this is just determined by the Sobolev space $\d_{p,o}$ lies in, and its relation to the order of a Lagrangian distribution). By picking $o$ close enough to $\p M$, the intersection of $\mcL^p_0$ and $\mcF^p_I$ is transverse at $\dot{\mcF}^p_{I,q}$ as shown in \cite[Section 9]{MVWCorners}.\footnote{These are facts from symplectic geometry, unrelated to any particulars of the elastic equation.}
 Hence, by the analogous proof of \cite[Corollary 9.7]{MVWCorners}, which first brings $\mcL^p_0$ and $(\mcF^p_{I})^o$ to respective normal forms, there is a compact interval $\mathbf{I} \subset (t_0, \bar{t})$ (i.e. before the first $p$-bicharacteristic hits the boundary, and such that $\mcF^p_I$ is still well-defined on this interval)
 
 such that $u_0$ is coinvolutive of some large order $N>0$  on a neighborhood of $\dot{\mcF}^p_{I,q}(\I)$ with respect to $H^{-n/2+1 +f/2 - 0}$. Thus, there is a neighborhood $U_0$ of $\dot{\mcF}^p_{I,q}(\I)$ such that
   \begin{equation}\label{eq:fund B_0 coinv refularity}
   B_0 u_0 \in \mcA_p^N(H^{-n/2+1+f/2-0}) \text{ if $B_0 \in \Psi^0_c(M^o)$}, \ \WF'(B_0) \subset U_0.
   \end{equation}

By a standard microlocal partition of unity, when restricted to $\I \times X^o$, $u_0$ is the sum of a distribution in $H^{\infty}_{b,\mcD'}$ and a distribution in $\mcA^N_p(H^{-n/2+1+f/2-0})$, the former coming from the part of $u_0$ microsupported away from $\mcL^p_0.$

Hence, together with (\ref{eq:fund u_o nice on geom related to beta'}), $u_0$ satisfies the hypothesis of Theorem \ref{thm:Diffraction2}, and we may conclude by that theorem that $u_0$ lies in $H^{-n/2 +1 + f/2 -0}$ microlocally at $w$.

We may now turn to the piece $u_1=P^{-1}_+(\delta_{p,1})$.

{\bf \underline{Step 1.2: Getting regularity of $u_1$ piece}}

By construction of $u_1$ and propagation of singularities, we conclude (by an analogous argument as done for $u_0$) that $u_1$ is microsupported away from a neighborhood $V_{\alpha}$ of $\bar{\mcF}^p_{I,\alpha}$ for $t < \bar{t}$. Thus, letting $T_0 < \bar{t}$, since $\bar{t}$ is the \emph{first} time a $p$-bicharacteristic hits the boundary, then $u_1(T_0)$ is smooth in a neighborhood of $\partial X$. Let $B_{\epsilon_1}(\rho_X(\alpha))$ be a neighborhood contained in $\rho_X( \bar{\mcF}^p_{I,\alpha})$. Then assuming $o$ is close enough to $\partial X$ that no $p$-geodesics intersect each other up until time $\bar{t}$, then the microsupport property of $u_1$ implies $u_1$ is smooth on $B_{\epsilon_1}(\rho_X(\alpha))$ for $t\leq T_0$. Hence, by finite propagation speed via Theorem \ref{cor:speedFiniteSpeedOfProp}, $u_1$ is smooth on some interval of time past $\bar{t}$, but in a smaller ball. In particular, it is microlocally smooth on some point of $\mcF^{p}_{O,\beta}$ and hence at $w$ by propagation of singularities.

 {\bf \underline{Step 2: Getting regularity of the $u_s$ piece}}

This piece is even easier. By the same argument used to obtain (\ref{eq:fund Lagrangian for u_{p,0}}), if we let $\bar{t}_s$ denote the first time an $s$-bicharacteristic lying over $S^*_{(t_0,o)}M^o$ hits the boundary, then for times $t<\bar{t}_s$, $u_s$ is a Lagrangian distribution associated to the Lagrangian
\begin{equation}
\mcL^s = \mcL(\WF'(\d_s) \cap \Sigma) = \mcL^s( S^*_{(t_0,o)}M^o \cap \Sigma_s).
\end{equation}
Since $c_s < c_p$, then $\bar{t}_s < \bar{t}$ and so $u_s$ is smooth for times inside $(\bar{t}_s, \bar{t})$. Thus, we may find a geodesic ball $B$ centered at $\rho(\alpha)$, in which $u_s$ is in $H^{\infty}_{b,\tmcD}$. So by the same argument used for $u_{p,1}$ invoking Theorem \ref{cor:speedFiniteSpeedOfProp}, we conclude that $u_s$ lies in $H^{-n/2+1+f/2-0}$ microlocally at $w$ as well (notice that here, it is irrelevant whether $w \in \mcF^p$ or $w \in \mcF^s$, which only mattered when we spoke about $u_0$, since the finite propagation speed corollary makes no reference to $p$ or $s$ bicharacteristics).

The analogous argument then works when analyzing the first $s$-bicharacteristic lying over $S^*_{(t_0,o)}M^o$ to hit the boundary.
$\Box$

\section{Acknowledgments}
I would like to thank my advisor Andras V\'{a}sy for inspiring me to study microlocal analysis and for introducing to me the interesting problem of diffraction for elastic waves. I have had numerous helpful discussions with Andras that have forged my intuition for the topic. Without Andras's constant encouragement and patience, the completion of this project would not have been possible.

 \bibliography{researchStatement}
    \bibliographystyle{plain}
    \nocite{*}
    
\end{document}